\documentclass[12pt]{article}
\usepackage{amsmath}
\usepackage{epsfig}
\usepackage{graphicx}
\usepackage{amsthm}
\usepackage{amsfonts}
\usepackage{color}
\usepackage{amssymb}
\usepackage{amsmath}

\newcommand{\field}[1]{\mathbb{#1}}
\newcommand{\R}{\field{R}}

\newcommand{\N}{\field{N}}
\newcommand{\Z}{\field{Z}}

\newcommand{\E}{\field{E}}

\def\qed{\hfill$\diamondsuit$}

\theoremstyle{Conjecture} \theoremstyle{example}
\theoremstyle{remark} \theoremstyle{lemma}
\theoremstyle{definition} \theoremstyle{corol}
\theoremstyle{proposition} \theoremstyle{condition}
\newtheorem{theorem}{Theorem}

\newtheorem{remark}{Remark}

\def\cov{{\mbox{cov}}}
\def\var{{\mbox{var}}}

\begin{document}

\bibliographystyle{plain}

\begin{center}
{\Large  Fixed-$b$  Subsampling and Block Bootstrap: Improved Confidence Sets Based on P-value Calibration
}
\end{center}

\bigskip

\centerline{\sc By Xiaofeng Shao and Dimitris N. Politis\footnote{Xiaofeng Shao is Associate Professor, Department of Statistics, University
of Illinois at Urbana-Champaign, Champaign, IL, 61820, USA. Dimitris N. Politis is Professor, Department of Mathematics, University of
California, San Diego, La Jolla, CA, 92093, USA. Shao's research is supported in part by
NSF grant DMS-11-04545 and Politis's research is partially supported by NSF Grant DMS-10-07513.
Emails: xshao@uiuc.edu, politis@math.ucsd.edu. The authors would like
to thank Xianyang Zhang for helpful comments and Jason Schweinsberg
for suggesting some useful references on Gaussian processes. We are also grateful to two referees for constructive comments, which led to substantial improvements.  }}


\centerline{\today}

\begin{abstract}

Subsampling and block-based bootstrap methods have been used in a wide range of inference problems for time series.
To accommodate the dependence, these resampling methods involve a bandwidth parameter, such as  subsampling
window width and block size in the block-based bootstrap. In empirical work, using
different bandwidth parameters could lead to different inference results, but the traditional first order asymptotic theory does not capture the
choice of the bandwidth. In this article, we propose to adopt the fixed-$b$ approach, as advocated by Kiefer and Vogelsang (2005) in the heteroscedasticity-autocorrelation robust testing context, to account for the influence of the bandwidth on the inference. Under the fixed-$b$ asymptotic framework, we derive the asymptotic null distribution of the p-values for subsampling and the moving block bootstrap, and further propose a calibration of the traditional small-$b$ based confidence intervals (regions, bands) and tests. Our treatment is fairly general as it includes both finite dimensional parameters and infinite dimensional parameters, such as  marginal distribution function and normalized spectral distribution function.    Simulation results show that the fixed-$b$ approach is more accurate than the traditional small-$b$ approach in terms of approximating the finite sample distribution, and that the calibrated confidence sets tend to have smaller coverage errors  than the uncalibrated counterparts.

\bigskip

\end{abstract}

\newpage

\section{Introduction}

Subsampling and block-based bootstrap methods have been widely used in inference problems for time series; see Politis et al. (1999a) and
Lahiri (2003) for book-length treatments of these important resampling methods. To accommodate the unknown time series dependence nonparametrically, these methods introduce a bandwidth parameter $l_n$, such as the block size in the block-based bootstrap and the subsampling window width in subsampling.   The bandwidth $l_n$ plays an important role in the finite sample performance of subsampling or block bootstrap based inference. Intuitively, if the bandwidth (or block size) is too small, it may not capture the dependence in a time series sufficiently, whereas if it is too large, the number of blocks for subsampling/resampling is too small to lead to a good approximation of finite sample distribution. Statistically speaking, the bandwidth $l_n$ is a smoothing parameter as it usually leads to a bias-variance tradeoff in variance estimation or size-power tradeoff in testing on the basis of subsampling and block bootstrap. In the traditional asymptotic theory, $l_n$ goes to infinity as sample size $n$ goes to infinity and the fraction $b=l_n/n$ goes to zero, which is a necessary condition for the general consistency of subsampling and block-based bootstrap methods without additional assumptions. Therefore, the role of $l_n$ (or $b$) does not show up in the conventional first order asymptotics, although in practice the choice of $l_n$ does affect the subsampling/block bootstrap distribution estimator and related operating characteristics.

In this paper, we aim to offer a new perspective on the use
of these smoothing parameter dependent resampling methods  based on the so-called fixed-$b$ approach, which
 was first proposed by Kiefer and Vogelsang (2005) in the context of heteroscedasticity-autocorrelation robust (HAR) testing. It was found that the asymptotic distribution obtained under the fixed-$b$ framework (i.e. $b\in (0,1]$ is held fixed in the asymptotics) provides a better approximation to the sampling distribution of the studentized test statistic than its counterpart obtained under the small-$b$ framework (i.e., $b\rightarrow 0$
as $n\rightarrow\infty$). See Jansson (2004) and Sun et al. (2008) for rigorous theoretical justifications.
The fixed-$b$ approach has the advantage of accounting for the effect of the bandwidth, as different bandwidth parameters correspond to
different limiting (null) distributions.
The literature on inference using the fixed-$b$ approach and its variants has been growing steadily; see Hashimzade and Vogelsang (2008),  Sun et al. (2008), Shao (2010a), Goncalves and Vogelsang (2011), and Sayginsoy and Vogelsang (2011)  among others for recent contributions.

In this paper,  we adopt Kiefer and Vogelsang's fixed-$b$ approach and
investigate its possible gain in the context of subsampling [Politis and Romano (1994)] and the moving block bootstrap [K\"unsch (1989) and Liu and Singh (1992)].
 The extension to other bandwidth-dependent bootstrap methods, such as the tapered block bootstrap [Paparoditis and Politis (2001, 2002), Shao (2010b)] and the dependent wild bootstrap [Shao (2010c)] are possible but are not pursued here. Under the fixed-$b$ asymptotics, Lahiri (2001) showed that the subsampling and the moving block bootstrap approximations are no longer consistent in the case of sample mean, which seems to suggest that a direct application of the fixed-$b$ approach is fruitless. A novel feature of our extension is that we study the limiting null distribution of the p-value, which is  $U(0,1)$ (i.e., uniform distribution on $[0,1]$) under the small-$b$ asymptotics, but is dependent upon $b$ and differs from $U(0,1)$ under the fixed-$b$ asymptotics. For a scalar parameter,  we calibrate the nominal coverage level  on the basis of the pivotal limiting null distribution of the p-value under the fixed-$b$ framework, and modify the small-$b$ based confidence interval by inverting the corresponding test. Thus the impact of the bandwidth parameter $l_n$ on the subsampling/block bootstrap distribution approximation is captured to the first order  using a p-value based adjustment.
   Simulation studies are conducted to demonstrate that the fixed-$b$ approach delivers  confidence intervals of better coverage in most situations and that the fixed-$b$ based intervals are slightly wider than the small-$b$ counterparts,  consistent with early findings associated with the fixed-b approach; see e.g. Kiefer and Vogelsang (2005).

   So far the use of the fixed-$b$ approach has been restricted to the inference of a finite dimensional parameter. Since the subsampling and moving block bootstrap have also been used in the inference of infinite dimensional parameters, such as marginal distribution function and (normalized) spectral distribution function of a stationary time series, we explore an extension of the fixed-$b$ idea to construct  confidence bands for these infinite dimensional parameters. Unlike the case of  a scalar parameter, the limiting null distribution of the subsampling-based p-value  is not pivotal under the fixed-$b$ asymptotics and it depends on
   the unknown  dependence structure of the underlying process. To alleviate the problem, we apply the subsampling method to approximate the sampling distribution of the p-value so inference becomes feasible. This double subsampling approach is also used in constructing the confidence region  for  a vector parameter.

The remainder of the article is organized as follows. Section~\ref{sec:mean} introduces an extension of the fixed-$b$ approach
to subsampling and the moving block bootstrap in the mean case. Section~\ref{sec:general} lays out a general framework and describes
the fixed-$b$ based confidence interval (region) for a finite dimensional parameter and some implementational issues. Section~\ref{sec:infinite} presents an extension of
the fixed-$b$ approach to confidence band construction for the marginal distribution function and normalized spectral distribution function. Simulation results  are reported in Section~\ref{sec:simulation}. Section~\ref{sec:conclusion} concludes and Section~\ref{sec:appendix}
contains some technical details.

\section{Inference for the mean}
\label{sec:mean}

To help the reader understand the essence of the fixed-$b$ approach, we shall focus on the simple problem:
inference for the mean of a stationary time series.  Suppose we want to test $H_0: \mu=\mu_0$ versus
$H_1:\mu\not=\mu_0$ based on the observations $\{X_t\}_{t=1}^{n}$ from a univariate stationary time
series with $\E(X_t)=\mu$. Under suitable moment and weak dependence conditions, we have $\sqrt{n}(\bar{X}_n-\mu)\rightarrow_{D} N(0,\sigma^2)$
where $\sigma^2=\sum_{k\in\Z}\gamma(k)$ is the long run variance with $\gamma(k)=\cov(X_0,X_k)$ and ``$\rightarrow_{D}$'' denotes convergence in distribution.
The scale parameter $\sigma^2$ can be consistently estimated by
  the so-called lag window estimator $\hat{\sigma}_n^2=\sum_{j=1-n}^{n-1}K(j/l)\hat{\gamma}_n(j)$, where
  $l=l_n$ is a bandwidth parameter, $K(\cdot)$ is a kernel function and $\hat{\gamma}_n(j)=n^{-1}\sum_{k=|j|+1}^{n}(X_k-\bar{X}_n)(X_{k-|j|}-\bar{X}_n)$ is the sample autocovariance at lag $j$.

  A natural test statistic is $G_n=n(\bar{X}_n-\mu_0)^2/\hat{\sigma}_n^2$.
 To ensure the consistency of $\hat{\sigma}_n^2$ as an estimator of $\sigma^2$, the bandwidth parameter $l=bn$, where $b\in (0,1]$, typically satisfies $1/l+l/n=o(1)$ (i.e., $b + n^{-1}/b =o(1)$) as $n\rightarrow\infty$. This is the so-called small-$b$ asymptotics, under which the limiting null distribution of $G_n$ is the distribution of $\chi_1^2$. Under the fixed-$b$ asymptotics, the ratio of bandwidth to sample size $b$ is
held fixed and  $G_n$ converges in distribution (under the null) to $U(b)$, whose detailed form can be found in  Kiefer and Vogelsang (2005).
 The distribution of $U(b)$ depends on the kernel $K$ and $b$, so different choices of the kernel and bandwidth lead
to different limiting null distributions. From both empirical and theoretical perspectives, the fixed-$b$ approach has been shown to
provide a more accurate approximation to the finite sample distribution of $G_n$ than the small-$b$ counterpart under the null, so it corresponds to better size in hypothesis testing. Owing to the duality between confidence interval construction and hypothesis testing, the interval delivered by the fixed-$b$ approach tends to have an empirical coverage closer to the nominal one.


In the next two subsections, we describe an extension of the fixed-$b$ approach to subsampling and the moving block bootstrap for the mean inference problem. A further extension to the inference of a finite dimensional parameter is made in Section~\ref{sec:general}. Throughout, we use
 $\lfloor a\rfloor$ to denote the integer part of $a\in\R$ and $\lceil a\rceil$ to denote the smallest integer larger than or equal to $a$. The symbol $N(\mu,\Sigma)$ denotes the normal distribution with mean $\mu$ and covariance matrix $\Sigma$.

\subsection{Subsampling}
\label{subsec:ss}

For the inference of the mean, the subsampling method approximates the sampling distribution
of $\sqrt{n}(\bar{X}_n-\mu)$ with the empirical distribution generated by its subsample counterpart $\sqrt{l}(\bar{X}_{j,j+l-1}-\bar{X}_n)$, where
$\bar{X}_{j,j+l-1}=l^{-1}\sum_{i=j}^{j+l-1}X_i$, $j=1,\cdots,N=n-l+1$. Let $L_{n,l}(x)=N^{-1}\sum_{j=1}^{N}{\bf 1}\{\sqrt{l}(\bar{X}_{j,j+l-1}-\bar{X}_n)\le x\}$ be the subsampling approximation, where ${\bf 1}(A)$ denotes the indicator function of the set $A$.
 For a given $\alpha\in (0,1)$ (say, $\alpha=0.05$ or $0.1$), we define the subsampling-based critical values as
 $c_{n,l}(1-\alpha)=\inf\{x: L_{n,l}(x)\ge 1-\alpha\}$. Then  a $100(1-\alpha)\%$ (one-sided) confidence
 interval is $[\bar{X}_n-n^{-1/2}c_{n,l}(1-\alpha),\infty)$ and the $100(1-\alpha)\%$ (two-sided) equal-tailed
confidence  interval is $[\bar{X}_n-n^{-1/2}c_{n,l}(1-\alpha/2),\bar{X}_n-n^{-1/2}c_{n,l}(\alpha/2)]$.
 In the testing context, if the alternative hypothesis is $H_1: \mu>\mu_0$, then we reject the null hypothesis at the significance level $\alpha$ if
 $\mu_0\notin [\bar{X}_n-n^{-1/2}c_{n,l}(1-\alpha),\infty)$; if the alternative hypothesis is $H_1:\mu<\mu_0$, then
 the null is rejected provided that $\mu_0\notin (-\infty,\bar{X}_n-n^{-1/2}c_{n,l}(\alpha)]$, which is also a one-sided confidence interval with nominal level $(1-\alpha)$; If the alternative hypothesis is
 $H_1:\mu\not=\mu_0$, then the null hypothesis is rejected when $\mu_0\notin [\bar{X}_n-n^{-1/2}c_{n,l}(1-\alpha/2),\bar{X}_n-n^{-1/2}c_{n,l}(\alpha/2)]$.
The above inference is based on the traditional small-$b$  based asymptotic theory, under which  $\sup_{x\in\R} | L_{n,l}(x)-\Phi(x/\sigma)|= o_p(1)$, where $\Phi$ is the distribution function for the standard normal distribution,  and $\sup_{x\in \R} |L_{n,l}(x)-P\{\sqrt{n}(\bar{X}_n-\mu)\le x\}|=o_p(1)$, i.e., the subsampling method provides a consistent approximation to
the sampling distribution of $\sqrt{n}(\bar{X}_n-\mu)$ and its limiting distribution. Note that we are using the data-centered subsampling distribution for testing as
recommended by Berg et al. (2010).

 Under the fixed-$b$ framework, $L_{n,l}$ does not converge to $N(0,\sigma^2)$ in distribution. Instead, Lahiri (2001) showed that the limit of $L_{n,l}(x)$ is
 \[ (1-b)^{-1}\int_0^{1-b}{\bf 1}[\{W(b+t)-W(t)-bW(1)\}\sigma/\sqrt{b}\in (-\infty,x]] dt,\]
where $W(t)$ is a standard Brownian motion. Since $L_{n,l}$ converges to a random measure,   the  subsampling-based inference is asymptotically invalid under the fixed-$b$ framework. To alleviate the problem, we shall modify the traditional subsampling-based inference procedure by considering the subsampling-based p-value and its limiting null distribution. For the one-sided alternative hypothesis $H_1:\mu>\mu_0$, we define the p-value as
\[pval_{n,l}^{\mbox{SUB}}:=\frac{1}{N}\sum_{j=1}^{N} {\bf 1}\{\sqrt{n}(\bar{X}_n-\mu_0)\le \sqrt{l}(\bar{X}_{j,j+l-1}-\bar{X}_n)\}.\]
Under the small-$b$ asymptotics, it can be shown that $pval_{n,l}^{\mbox{SUB}}$ converges to $U[0,1]$ in distribution under the null, whereas under the fixed-b
asymptotics, its limiting null distribution is the distribution of $G(b)$, where
\[G(b)=(1-b)^{-1}\int_0^{1-b} {\bf 1}[W(1)\le \{W(b+t)-W(t)-bW(1)\}/\sqrt{b}] dt.\]
Note that the nuisance parameter $\sigma$ is canceled out in $G(b)$, which is pivotal for a given $b$.
Let $G_{\alpha}(b)$ denote the $100\alpha\%$ quantile of the distribution $G(b)$. Then
at the significance level $\alpha$, we reject the null and favor the alternative $H_1:\mu>\mu_0$, if the (realized) p-value is smaller than $G_{\alpha}(b)$. Correspondingly, a one-sided confidence interval under the fixed-$b$ asymptotics can be obtained by inverting the test, i.e.,
$\{\mu: \frac{1}{N}\sum_{j=1}^{N} {\bf 1}\{\sqrt{n}(\bar{X}_n-\mu)\le \sqrt{l}(\bar{X}_{j,j+l-1}-\bar{X}_n)\}\ge G_{\alpha}(b)\}$, which is
$$\{\mu: \sqrt{n}(\bar{X}_n-\mu)\le c_{n,l}(1-G_{\alpha}(b))\}=[\bar{X}_n-n^{-1/2}c_{n,l}(1-G_{\alpha}(b)),\infty).$$
 Compared to the conventional subsampling-based confidence interval, the difference lies in the replacement of $\alpha$ by $G_{\alpha}(b)$ in $c_{n,l}$.  Note that $\alpha$ is the $100\alpha\%$ quantile of $U(0,1)$, which is the limiting null distribution of the p-value under the small-$b$ asymptotics. In a similar manner, we can obtain the $100(1-\alpha)\%$
two sided equal-tailed confidence interval for $\mu$ as $[\bar{X}_n-n^{-1/2}c_{n,l}(1-G_{\alpha/2}(b)),\bar{X}_n-n^{-1/2}c_{n,l}(G_{\alpha/2}(b))]$ and another one-sided confidence interval  $(-\infty,\bar{X}_n-n^{-1/2}c_{n,l}(G_{\alpha}(b))]$ under the fixed-$b$ asymptotics. One can view the fixed-$b$ based inference  as a way of calibrating the small-$b$ counterpart with the level $\alpha$ adjusted by $G_{\alpha}(b)$, so the effect of $b$ on the inference is taken into account.

If one wants to construct a  symmetric two sided confidence interval for $\mu$, then one can approximate the sampling distribution of
$\sqrt{n}|\bar{X}_n-\mu|$ by $\widetilde{L}_{n,l}(x)=N^{-1}\sum_{j=1}^{N}{\bf 1}(\sqrt{l}|\bar{X}_{j,j+l-1}-\bar{X}_n|\le x)$. Letting $\widetilde{c}_{n,l}(1-\alpha)=\inf\{x: \widetilde{L}_{n,l}(x)\ge 1-\alpha\}$, then the $100(1-\alpha)\%$ symmetric confidence interval for $\mu$ is $[\bar{X}_n-n^{-1/2}\widetilde{c}_{n,l}(1-\alpha),\bar{X}_n+n^{-1/2}\widetilde{c}_{n,l}(1-\alpha)]$ under the small-$b$ asymptotic theory. The p-value is defined as
$$\widetilde{pval}_{n,l}^{SUB}=\frac{1}{N}\sum_{j=1}^{N} {\bf 1}\{\sqrt{n}|\bar{X}_n-\mu_0|\le \sqrt{l}|\bar{X}_{j,j+l-1}-\bar{X}_n|\}.$$
Under the fixed-$b$ asymptotics, the limiting null distribution of $\widetilde{pval}_{n,l}^{SUB}$ is the distribution of $\widetilde{G}(b)$, where
\[\widetilde{G}(b)=(1-b)^{-1}\int_0^{1-b} {\bf 1}\{|W(1)|\le |W(b+t)-W(t)-bW(1)|/\sqrt{b}\} dt.\]
Let $\widetilde{G}_{\alpha}(b)$ denote the $100\alpha\%$ quantile of the distribution $\widetilde{G}(b)$. Then  the fixed-$b$ based  $100(1-\alpha)\%$ symmetric confidence interval is
$\{\mu: N^{-1}\sum_{j=1}^{N}{\bf 1}(\sqrt{n}|\bar{X}_n-\mu|\le \sqrt{l}|\bar{X}_{j,j+l-1}-\bar{X}_n|) \ge \widetilde{G}_{\alpha}(b)\}$, i.e.,
\begin{eqnarray}
\label{eq:ci1}
[\bar{X}_n-n^{-1/2}\widetilde{c}_{n,l}(1-\widetilde{G}_{\alpha}(b)),\bar{X}_n+n^{-1/2}\widetilde{c}_{n,l}(1-\widetilde{G}_{\alpha}(b))].
\end{eqnarray}

\subsection{Moving block bootstrap}
\label{subsec:mbb}

In this subsection, we shall consider the approximation of the sampling distribution of $\sqrt{n}(\bar{X}_n-\mu)$ by the moving block bootstrap (MBB).
Denote the MBB sample by $\{X_1^*(l),\cdots,X_n^*(l)\}$ with the dependence on the block size $l$ being explicit. Then we approximate the sampling distribution of $\sqrt{n}(\bar{X}_n-\mu)$ by the conditional distribution of $\sqrt{n}\{\bar{X}_n^*(l)-\bar{X}_n\}$ or $\sqrt{n}[\bar{X}_n^*(l)-\E^*\{\bar{X}_n^*(l)\}]$ given the data, where  $\bar{X}_n^*(l)=n^{-1}\sum_{t=1}^{n}X_t^*(l)$ is the sample mean for the bootstrap sample, $\E^*$ and $\var^*$ are used to denote the conditional expectation and variance, respectively. To avoid the issue of centering, we could use the circular bootstrap [Politis and Romano (1992)], which is asymptotically equivalent to the moving block bootstrap [Lahiri (2003)]. For simplicity, we shall focus on the bootstrap approximation based on $\sqrt{n}\{\bar{X}_n^*(l)-\bar{X}_n\}$. The same idea can be applied to the other bootstrap approximation. We define the MBB-based p-value as
\[pval_{n,l}^{\mbox{MBB}}:=E^*[{\bf 1}\{\sqrt{n}(\bar{X}_n-\mu_0)\le \sqrt{n}(\bar{X}_{n}^*(l)-\bar{X}_n)\}],\]
which corresponds to the alternative $H_1:\mu>\mu_0$.
Under the small-$b$ asymptotics, the p-value $pval_{n,l}^{\mbox{MBB}}$ is expected to converge to $U(0,1)$ in distribution under the null, although we are unaware of a formal proof.
Under the fixed-$b$ asymptotics, we assume $R_b=n/l=1/b$ (i.e. reciprocal of $b$) to be an integer for the ease of our discussion. Then  $\bar{X}_{n}^*(l)=n^{-1}\sum_{j=1}^{R_b l}X_j^*(l)=n^{-1}\sum_{j=1}^{R_b} v_j^*$, where, conditional on the data,  $\{v_j^*\}_{j=1}^{R_b}$ are iid (independent and identically distributed)
 with a discrete uniform distribution, $P(v_1^*=\sum_{t=j}^{j+l-1}X_t)=1/N$, $j=1,\cdots,N$. Hence the above p-value is equal to
\[\frac{1}{N^{R_b}} \sum_{j_1,j_2,\cdots,j_{R_b}=1}^{N} {\bf 1}\left\{ n^{-1/2}\sum_{h=1}^{R_b}\sum_{s=j_h}^{j_h+l-1} (X_s-\mu_0)-n^{1/2}(\bar{X}_n-\mu_0) \ge n^{1/2} (\bar{X}_n-\mu_0) \right\}. \]
Under the fixed-$b$ asymptotics and the null, it converges in distribution to
\[H(b):=(1-b)^{-R_b} \int_0^{1-b}\cdots\int_{0}^{1-b} {\bf 1}\left[\sum_{h=1}^{R_b} \{W(t_h+b)-W(t_h)\} \ge 2 W(1)\right] dt_1\cdots dt_{R_b}\]
Let $H_{\alpha}(b)$ denote the $100\alpha\%$ quantile of $H(b)$. In practice, we usually further approximate the distribution of $\sqrt{n}\{\bar{X}_n^*(l)-\bar{X}_n\}$ by taking a finite number of bootstrap samples, say,
$\{X_t^{*,j}(l)\}_{t=1}^{n}$, $j=1,\cdots,B$. We approximate the sampling distribution of $\sqrt{n}(\bar{X}_n-\mu)$ by
$M_{n,l,B}^*(x)=\frac{1}{B}\sum_{j=1}^{B} {\bf 1}[\sqrt{n}\{\bar{X}_n^{*,j}(l)-\bar{X}_n\}\le x]$, where $\bar{X}_n^{*,j}(l)=n^{-1}\sum_{t=1}^{n}X_t^{*,j}(l)$.
Let $c_{n,l,B}^*(1-\alpha)=\inf\{x:M_{n,l,B}^*(x)\ge 1-\alpha\}$. The corresponding  fixed-$b$ based two sided equal tailed confidence interval for $\mu$ is then
$$[\bar{X}_n-n^{-1/2}c_{n,l,B}^*(1-H_{\alpha/2}(b)),\bar{X}_n-n^{-1/2}c_{n,l,B}^*(H_{\alpha/2}(b))]$$
 and the one-sided confidence intervals can be formed analogous to those developed for the subsampling method. The details are omitted. The above discussion is based on the assumption that $R_b=1/b$ is an integer.  When $R_b$ is not an integer, we use a fraction of the last resampled block to make the bootstrap sample size equal to original sample size. Then the p-value is
\begin{eqnarray*}
&&\frac{1}{N^{\lfloor R_b\rfloor }} \frac{1}{l\lfloor R_b\rfloor+1} \sum_{j_1,j_2,\cdots,j_{\lfloor R_b\rfloor}=1}^{N} \sum_{j_{\lfloor R_b\rfloor+1}=1}^{l\lfloor R_b\rfloor+1} {\bf 1}\left\{ n^{-1/2}\left(\sum_{h=1}^{\lfloor R_b\rfloor}\sum_{s=j_h}^{j_h+l-1} (X_s-\bar{X}_n) \right.\right.\\
&&\left.\left.+\sum_{s=j_{\lfloor R_b\rfloor+1}}^{j_{\lfloor R_b\rfloor+1}+n-l\lfloor R_b\rfloor-1}(X_s-\bar{X}_n)\right) \ge n^{1/2} (\bar{X}_n-\mu_0) \right\}.
\end{eqnarray*}
and its limiting null distribution can be derived similarly. Below we shall focus our discussion on the case $1/b$ is an integer for simplicity.

In a similar fashion, if we want to construct an MBB-based symmetric confidence interval for $\mu$, we consider the approximation of the sampling distribution of $\sqrt{n}|\bar{X}_n-\mu|$ by the conditional distribution of $\sqrt{n}|\bar{X}_n^*(l)-\bar{X}_n|$ given the data. The corresponding p-value is
\begin{eqnarray*}
&&\widetilde{pval}_{n,l}^{\mbox{MBB}}:=E^*\{{\bf 1}(\sqrt{n}|\bar{X}_n-\mu_0|\le \sqrt{n}|\bar{X}_{n}^*(l)-\bar{X}_n|)\}\\
&&\hspace{0.5cm}=\frac{1}{N^{R_b}}\sum_{j_1,j_2,\cdots,j_{R_b}=1}^{N} {\bf 1}\left( \left|n^{-1/2} \sum_{h=1}^{R_b}\sum_{s=j_h}^{j_h+l-1} (X_s-\mu_0)- \sqrt{n} (\bar{X}_n-\mu_0)\right|\ge \sqrt{n}|\bar{X}_n-\mu_0| \right)
\end{eqnarray*}
and it converges in distribution to
\[\widetilde{H}(b)=(1-b)^{-R_b} \int_0^{1-b}\cdots\int_{0}^{1-b} {\bf 1}\left(\left|\sum_{h=1}^{R_b} \{W(t_h+b)-W(t_h)\} -W(1)\right| \ge  |W(1)| \right) dt_1\cdots dt_{R_b}\]
under the null. Define $\widetilde{M}_{n,l,B}^*(x)=\frac{1}{B}\sum_{j=1}^{B} {\bf 1}\{\sqrt{n}|\bar{X}_n^{*,j}(l)-\bar{X}_n|\le x\}$ and $\widetilde{c}_{n,l,B}^*(1-\alpha)=\inf\{x:\widetilde{M}_{n,l,B}^*(x)\ge 1-\alpha\}$.
Then the fixed-$b$ based $100(1-\alpha)\%$ symmetric confidence interval for $\mu$ is $[\bar{X}_n-n^{-1/2}\widetilde{c}_{n,l,B}^{*}(1-\widetilde{H}_{\alpha}(b)),\bar{X}_n+n^{-1/2}\widetilde{c}_{n,l,B}^*(1-\widetilde{H}_{\alpha}(b))]$.

\section{Finite dimensional parameter}
\label{sec:general}

We first introduce some notation. Let $D[0,1]$ be the space of
functions on $[0,1]$ which are right continuous and have left limits, endowed with the Skorokhod topology (Billingsley 1968). Denote by ``$\Rightarrow$" weak convergence in $D[0,1]$ or more generally in the $\R^k$-valued function space $D^k[0,1]$, where $k\in\N$. Later in Section~\ref{sec:infinite}, we also use ``$\Rightarrow$" to denote
the weak convergence in $D[0,\pi]$, $D([0,1]\times [0,\pi])$ and $D([-\infty,\infty]\times [0,1])$.

\subsection{Subsampling}
\label{subsec:subsampling}

Following Politis et al. (1999a), we assume that the parameter of interest is $\theta(P)\in\R^k$, where $P$ is the joint probability law that governs
the $p$-dimensional stationary sequence $\{X_t\}_{t\in\Z}$.
Let $\hat{\theta}_n=\hat{\theta}_n(X_1,\cdots,X_n)$ be an estimator of $\theta=\theta(P)$ based on the observations
$(X_1,\cdots,X_n)$. Further we define the subsampling estimator of $\theta(P)$ by $\hat{\theta}_{j,j+l-1}=\hat{\theta}_l(X_j,\cdots,X_{j+l-1})$ on the basis of
the subsample $(X_j,\cdots,X_{j+l-1})$, $j=1,\cdots,N$.
Let $\|\cdot\|$ be a norm in $\R^k$. The subsampling-based distribution estimator of $\|\sqrt{n}\{\hat{\theta}_n-\theta(P)\}\|$ is denoted as
$\widetilde{L}_{n,l}(x)=N^{-1}\sum_{j=1}^{N}{\bf 1}(\|\sqrt{l}(\hat{\theta}_{j,j+l-1} -\hat{\theta}_n)\|\le x)$.   In the testing context (say $H_0:\theta=\theta_0$ versus $H_1:\theta\not=\theta_0$), we define the subsampling based p-value as
\begin{eqnarray}
\label{eq:pvaluesub}
\widetilde{pval}_{n,l}^{SUB}=N^{-1}\sum_{j=1}^{N}{\bf 1}(\|\sqrt{n}(\hat{\theta}_n-\theta)\|\le \|\sqrt{l}(\hat{\theta}_{j,j+l-1}-\hat{\theta}_n)\|),
\end{eqnarray}
where we do not distinguish $\theta$ and $\theta_0$ for notational convenience
because they are the same under the null.


To obtain the limiting null distribution of the p-value under the fixed-$b$ framework, we further assume
$\theta(P)=T(F)$, where $F$ is the marginal distribution of $X_1\in \R^p$, and $T$ is a functional that takes value in $\R^k$. Then a natural estimator of $T(F)$ is $\hat{\theta}_n=T(\rho_{1,n})$, where $\rho_{1,n}=n^{-1}\sum_{t=1}^{n}\delta_{X_t}$ is the empirical distribution. Here $\delta_x$ stands for the point mass at $x$. Similarly, $\hat{\theta}_{j,j+l-1}=T(\rho_{j,j+l-1})$, where $\rho_{j,j+l-1}=l^{-1}\sum_{h=j}^{j+l-1}\delta_{X_h}$.
 Assume that there is an expansion of $T(\rho_{1,n})$ in the neighborhood of $F$, i.e., $$T(\rho_{1,n})=T(F)+n^{-1}\sum_{t=1}^{n}IF(X_t; F)+R_{1,n},$$
where $IF(X_t;F)$ is the influence function of $T$ (Hampel, Ronchetti, Rousseeuw and Stahel, 1986) defined by
$IF(x;F)=\lim_{\epsilon\downarrow 0}\frac{T((1-\epsilon)F+\epsilon\delta_x)-T(F)}{\epsilon}$ and $R_{1,n}$ is the remainder term. Similarly,
$T(\rho_{j,j+l-1})=T(F)+l^{-1}\sum_{h=j}^{j+l-1}IF(X_h;F)+R_{j,j+l-1}$. Below are the two key assumptions we need.

(A.1) Assume that $\E \{IF(X_j;F)\} = 0$ and $ n^{-1/2}\sum_{j=1}^{\lfloor nr\rfloor} IF(X_j;F)\Rightarrow \Sigma(P)^{1/2} W_k(r)$, where $\Sigma(P)$ is a positive definite matrix and $W_k(\cdot)$ denotes the $k$-th dimensional vector of independent Brownian motions.

(A.2) Assume that $\sqrt{n}\|R_{1,n}\|=o_p(1)$ and $\sqrt{l}\sup_{j=1,\cdots,N}\|R_{j,j+l-1}\|=o_p(1)$.

Note that (A.1) is Assumption 1 in Shao (2010a) and its verification has been discussed in Remark 1 therein.
  The assumption (A.2) is to ensure the negligibility of remainder terms.  In the sample mean case, $IF(X_t;F)=(X_t-\mu)$ and the remainder terms vanish, so (A.2) is automatically satisfied and (A.1)   reduces to a functional central limit theorem for the partial sum process of $X_t$.

\begin{theorem}
\label{th:finite}
Suppose the assumptions (A.1) and (A.2) hold and $b\in (0,1]$ is held fixed as $n\rightarrow\infty$. The limiting null distribution of
$\widetilde{pval}_{n,l}^{SUB}$ is the distribution of $\widetilde{G}(b;k)$, where
\[\widetilde{G}(b;k)=(1-b)^{-1}\int_0^{1-b} {\bf 1}[\|\Sigma(P)^{1/2} W_k(1)\|\le \|\Sigma(P)^{1/2}\{W_k(b+t)-W_k(t)-bW_k(1)\}\|/\sqrt{b}]dt.\]
In the special case $k=1$, $\widetilde{G}(b;1)=\widetilde{G}(b)$.
\end{theorem}

 Thus for a scalar parameter, the limiting null distribution of the p-value is pivotal for a given $b$
 and the $100(1-\alpha)\%$ symmetric  confidence interval for $\theta$ is
\[ [\hat{\theta}_n-n^{-1/2}\widetilde{c}_{n,l}(1-\widetilde{G}_{\alpha}(b)),\hat{\theta}_n+n^{-1/2}\widetilde{c}_{n,l}(1-\widetilde{G}_{\alpha}(b))],\]
 which reduces to (\ref{eq:ci1}) in the mean case. To conduct the inference for the case $k=1$, we need to know $G_{\alpha}(b)$, $\widetilde{G}_{\alpha}(b)$, $H_{\alpha}(b)$ and $\widetilde{H}_{\alpha}(b)$. Following the practice of Kiefer and Vogelsang (2005), we first generate the simulated values for $\alpha=0.05,0.1$ and $b=0.01,0.02,\cdots,0.2$, then
 fit  the quadratic equation $cv(b)=a_0+a_1b + a_2b^2$  to the simulated values by ordinary least squares. The intercept $a_0$
was set to be equal to $\alpha$, so that $cv(0)=\alpha$.  Table~\ref{tb:cv} reports the estimated coefficients and $R^2$ from the regressions (ranging from
0.9584 to 0.9997), which suggests quite satisfactory fits. For $H_{\alpha}(b)$ and $\widetilde{H}_{\alpha}(b)$, fitting higher order polynomials does not lead to substantial higher $R^2$.  To simulate $G_{\alpha}(b)$ and $\widetilde{G}_{\alpha}(b)$ for a given $\alpha$ and $b\in (0,0.2)$, we generate 5000 iid $N(0,1)$
random variables, and use its normalized partial sum to approximate the standard Brownian motion. For $H_{\alpha}(b)$ and $
\widetilde{H}_{\alpha}(b)$, we approximate $\E^*$ in the definition of p-value by performing bootstrap 50000 times. 50000 monte carlo replications were used for all the cases. For small $\alpha$ (say $\alpha=0.01$) and relatively large $b$, say $b=0.15, \cdots, 0.2$, our simulated critical values are mostly zero, so we are unable to
provide a fitted quadratic equation when $\alpha$ is very small. Nevertheless, if the goal is to construct a $90\%$ or $95\%$ symmetric confidence interval, or a $90\%$ equal-tailed confidence interval, or a one-sided confidence interval of nominal coverage $90\%$ or $95\%$,   Table~\ref{tb:cv} is useful when $b\in (0,0.2]$.

 \centerline{Please insert Table~\ref{tb:cv} about here!}

 For a vector parameter (i.e. $k\ge 2$), the limiting null distribution of the p-value depends on the unknown
 long run variance matrix, so is not pivotal in general. One way out is to approximate the limiting null distribution by subsampling. Denote by $n'$ the subsampling width at the second stage. Let $l'=\lceil n'b\rceil$ and $N'=n'-l'+1$. For each subsample $\{X_t,\cdots,X_{t+n'-1}\}$, we define the subsampling
 counterpart of $\widetilde{pval}_{n,l}^{SUB}$ as
 \[q_{n',t}=(N')^{-1} \sum_{j=t}^{t+N'-1} {\bf 1}\left\{\|\sqrt{l'}(\hat{\theta}_{j,j+l'-1}-\hat{\theta}_{t,t+n'-1})\|\ge \|\sqrt{n'}(\hat{\theta}_{t,t+n'-1}-\hat{\theta}_n)\| \right\}\]
 for $t=1,\cdots,n-n'+1$. Denote the empirical distribution function of $\{q_{n',t}\}_{t=1}^{n-n'+1}$ by $Q_{n,n'}(x)=(n-n'+1)^{-1}\sum_{j=1}^{n-n'+1}{\bf 1}(q_{n',j}\le x)$, which can be used to approximate the sampling distribution or the limiting null distribution of $\widetilde{pval}_{n,l}^{SUB}$. Let $c_{n,n',l}(1-\alpha)=\inf\{x: Q_{n,n'}(x)\ge 1-\alpha\}$. Then the calibrated $100(1-\alpha)\%$ subsampling-based confidence region for $\theta$ is
 \begin{eqnarray}
 \label{eq:calibregion}
 \{\theta\in \R^k: ~\widetilde{pva}_{n,l}^{SUB}~\mbox{in}~(\ref{eq:pvaluesub})~\ge c_{n,n',l}(\alpha)\},
 \end{eqnarray}
whereas the  traditional subsampling-based confidence region is $\{\theta\in \R^k: ~\widetilde{pva}_{n,l}^{SUB}~\mbox{in}~(\ref{eq:pvaluesub})~\ge \alpha\}$.

\begin{theorem}
\label{th:theoremregion}
Assume that $1/n'+n'/n=o(1)$ and $b\in (0,1]$ is fixed. Suppose that the process $X_t$ is $\alpha$-mixing and $\widetilde{G}(b;k)$ is a continuous random variable. Then we have $$\sup_{x\in\R} |Q_{n,n'}(x)-P(\widetilde{G}(b;k) \le x)|=o_p(1).$$ Consequently, the asymptotic coverage probability of
confidence region in~(\ref{eq:calibregion}) is $(1-\alpha)$.
\end{theorem}

\begin{remark}{\rm
As we have done subsampling twice, this procedure is naturally called double subsampling in the spirit of double bootstrap.
The use of subsampling at the second stage is mainly to approximate the sampling distribution or the limiting null distribution of the p-value, which is unknown under the fixed-$b$ asymptotic framework. Of course, the approximation error depends on  the subsampling window size $n'$ at the second stage. If we view $n'/n$ as a fixed constant in the above asymptotics, then the asymptotic coverage of the calibrated confidence region is still different from the nominal level. One can perform further calibration by subsampling, which leads to iterative subsampling, similar to iterative bootstrap in Beran (1987, 1988).
In practice, however, the selection of the subsampling window size at each stage usually involves quite expensive computation, and the (finite sample) improvement in coverage errors  is not guaranteed by doing subsampling iteratively.

As pointed out by a referee, a possible alternative approach is to simulate the asymptotic null distribution of the p-value, i.e. the distribution of $\widetilde{G}(b;k)$ after plugging in a consistent estimator of long run variance matrix. Note that in general consistent estimation of  long run variance matrix also involves the bandwidth selection; see e.g. Politis (2011). Since the above-mentioned double subsampling approach is also applicable to the infinitely dimensional case [see Section~\ref{sec:infinite}], we shall not pursue this alternative approach.

}
\end{remark}

\subsection{Moving block bootstrap}

 For the moving block bootstrap, we approximate the sampling distribution of $\|\sqrt{n}\{\hat{\theta}_n-\theta\}\|$
by the conditional distribution of $\sqrt{n}(\hat{\theta}_n^*-\hat{\theta}_n)$, where $\hat{\theta}_n^*=\hat{\theta}_n\{X_1^*(l),\cdots,X_n^*(l)\}$.
Define the p-value as $\widetilde{pval}_{n,l}^{MBB}:= \E^*\{ {\bf 1}(\|\sqrt{n}(\hat{\theta}_n-\theta_0)\|\le  \|\sqrt{n}(\hat{\theta}_n^*-\hat{\theta}_n)\|)\}$.
It can be expected that under certain regularity conditions, the limiting null distribution of $\widetilde{pval}_{n,l}^{MBB}$ is the distribution of
\begin{eqnarray*}
&&\widetilde{H}(b;k)=\frac{1}{(1-b)^{R_b}} \int_0^{1-b}\cdots\int_{0}^{1-b} {\bf 1}\left[\left\|\Sigma(P)^{1/2}\left[\sum_{h=1}^{R_b}  \{W_k(t_h+b)-W_k(t_h)\} -W_k(1)\right]\right\|\right.\\
&& \hspace{3cm} \left. \ge \|\Sigma(P)^{1/2} W_k(1)\| \right] dt_1\cdots dt_{R_b},
\end{eqnarray*}
which coincides with $\widetilde{H}(b)$ when $k=1$. When $k\ge 2$, the p-value is not asymptotically pivotal under the
fixed-$b$ asymptotics, and its sampling distribution can be approximated by subsampling or the moving block bootstrap.
Since the idea is similar to the double subsampling procedure described above, we omit the details. We mention in passing that Lee and Lai (2009)
have recently studied the benefit of performing double block bootstrap for the smooth function model.


The p-value based calibration is closely related to the prepivoting method proposed by Beran (1987, 1988).
The p-value of a statistic is itself a statistic that has a pivotal limiting distribution or  tends to be more pivotal than
the original (unstudentized) statistic. In Beran (1987), the limiting null distribution of the p-value was assumed to be
 $U(0,1)$, and he focused on the refinement of the approximation error of
sampling distribution of the p-value to $U(0,1)$ by prepivoting and iterative bootstrap. His treatment is quite
general but is mainly focused on the iid setting. By contrast, we deal with time series with independent data as a special case and the limiting null distribution
of the p-value (under the fixed-$b$ asymptotics) is not $U(0,1)$. In addition, our calibration can be applied to the inference of infinite dimensional parameters [see Section~\ref{sec:infinite}],
which is not covered by Beran (1987, 1988). Another related calibration method in the bootstrap literature was proposed by
Loh (1987, 1991), who calibrated confidence coefficients using a consistent estimate of actual coverage probability.
For a given confidence interval, its estimated coverage probability is used to
alter the nominal level of the interval, and it is shown that the calibrated interval is asymptotically robust under iid assumptions
and some regularity conditions. Similar to Beran's work, Loh's discussion is limited to the iid setting and  his calibration method seems only applicable to the inference of  finite dimensional parameters. For a comprehensive account of bootstrap iteration and calibration, see Hall (1992).

For a finite dimensional parameter, another way of making the statistic more pivotal is to do studentization using a consistent estimate
of asymptotic variance of the original statistic. For dependent data, this typically involves the estimation of long run
variance using the lag window type estimate. Although theoretically possible,   consistent estimation can be
difficult to carry out  in practice for some statistics. For example, if $k=p=1$, $\theta=\mbox{median}(F)$ and
$\hat{\theta}_n=\mbox{median}(X_1,\cdots,X_n)$, then
\[\Sigma(P)=\{4 g^2(\theta)\}^{-1}\sum_{k=-\infty}^{\infty}\cov\{1-2{\bf 1}(X_0\le \theta),1-2{\bf 1}(X_k\le \theta)\}\]
 with $g(\cdot)$ being the density function of $X_1$. Consistent estimation of $\Sigma(P)$ involves kernel density estimation
 for $g(\theta)$ and long run variance estimation for the transformed series $1-2{\bf 1}(X_t\le \theta)$, both of which involve
 the choice of a bandwidth parameter. By contrast,  subsampling and the moving block bootstrap can be used to provide consistent
 variance estimate, which lead to a studentized statistic, or a p-value, which is more pivotal than the original unstudentized statistic.
 Both methods  are relatively easier to
 implement, although they also require the user to choose the subsampling window width or block size. The self-normalized approach of Shao (2010a), which uses recursive subsample estimates in its studentization, would be another good candidate when a direct consistent long run variance estimation is difficult, although there is an efficiency loss under certain loss functions.


\section{Infinite dimensional parameter}
\label{sec:infinite}

In previous sections, our discussion focuses on the inference of a finite dimensional parameter, for which a $\sqrt{n}$-consistent
estimator exists and the asymptotic normality holds. In general, the use of subsampling and block bootstrap methods are not limited to the inference
for finite dimensional parameters. In the time series setting, they have been used to provide an approximation of the nonpivotal limiting distribution
when the parameter of interest is of infinite dimension, such as marginal distribution function and spectral distribution function of a stationary time series.
 In what follows, we use $\|F-G\|_{\infty}$ to denote $\sup_{x\in D}|F(x)-G(x)|$, where $D=[-\infty,\infty]$ in Section~\ref{subsec:mdf}
and $D=[0,\pi]$ in Section~\ref{subsec:sdf}.

\subsection{Marginal distribution function}
\label{subsec:mdf}

Consider a stationary sequence $\{X_k, k\in\Z\}$ and let $m(s)=P(X_0\le s)$
be its marginal distribution. Given the observations $\{X_t\}_{t=1}^{n}$, the empirical process is defined
as $m_n(s)=n^{-1}\sum_{k=1}^{n}{\bf 1}(X_k\le s)$. More generally, we define the standardized recursive process
$$K_n(s,\lfloor nt\rfloor)=n^{-1/2}\sum_{k=1}^{\lfloor nt\rfloor }\{{\bf 1}(X_k\le s)-m(s)\}, ~t\in [0,1].$$
Under certain regularity conditions [see Berkes et al. (2009)], we have that
\begin{eqnarray}
\label{eq:weakconv}
K_n(s,\lfloor nt\rfloor) \Rightarrow K(s,t).
\end{eqnarray}
Here $K(s,t)$, $(s,t)\in [-\infty,\infty] \times [0,1]$ is a two-parameter
mean zero Gaussian process with
$$\cov\{K(s,t),K(s',t')\}= (t\wedge t')\Gamma(s,s'),$$
where $\Gamma(s,s')=\sum_{k=-\infty}^{\infty}\cov\{{\bf 1}(X_0\le s),{\bf 1}(X_k\le s')\}$.
To construct a confidence band for $m(\cdot)$, we note that by the continuous mapping theorem, (\ref{eq:weakconv}) implies that
$\sqrt{n} \|m_n-m\|_{\infty} \rightarrow_D \sup_{s\in \R} |K(s,1)|.$
Since $K(s,1)$ is a Gaussian process with mean zero
and unknown covariance $\cov\{K(s,1),K(s',1)\}=\Gamma(s,s')$, direct inference of $m(\cdot)$ is difficult.
To circumvent the difficulty, both the moving block bootstrap and subsampling  have been proposed
to approximate the limiting distribution $\sup_{s\in \R}  |K(s,1)|$ consistently; see B\"uhlmann (1994), Naik-Nimbalkar  and Rajarshi (1994),  and Politis et al. (1999b).
Below we shall focus our discussion on the subsampling method, and a similar argument applies to the moving block bootstrap approach in view of the argument in Section~\ref{subsec:mbb}.
Let $g_n(t,s)=l^{1/2}\{m_{t,t+l-1}(s)-m_n(s)\}$, $t=1,\cdots, N=n-l+1$
be the subsampling counterpart of $n^{1/2} \{m_n(s)-m(s)\}$, where $m_{t,t+l-1}(s)=l^{-1}\sum_{h=t}^{t+l-1} {\bf 1}(X_h\le s)$ .
Assuming $l/n+1/l=o(1)$ and other regularity conditions, Politis et al. (1999b) showed that
the subsampling approximation based on $\{g_n(t,s)\}_{t=1}^{N}$ is consistent in certain function space. This implies that
the sampling distribution of $\sqrt{n}\|m_n-m\|_{\infty}$ (or the distribution of $\sup_{s\in\R} |K(s,1)|$)
 can be consistently approximated by the empirical distribution of
$\{ \sqrt{l}\|m_{t,t+l-1}-m_n\|_{\infty} \}_{t=1}^{N}$.
The above result is obtained under the small-$b$ asymptotics. To introduce our calibration method, we again start with the p-value and study its limiting null distribution under the fixed-$b$ asymptotics. For notational simplicity, we do not distinguish the true marginal distribution function $m(x)$ and the hypothesized function $m_0(x)$, because they are
 identical under the null hypothesis.

Define the p-value
\begin{eqnarray}
\label{eq:pvalueband}
pval_{n,l}^{\mbox{E}}=N^{-1}\sum_{j=1}^{N} {\bf 1}\left\{l^{1/2}\|m_{j,j+l-1}-m_n\|_{\infty} \ge \sqrt{n}\|m_n-m\|_{\infty} \right\}.
\end{eqnarray}
Let $b=l/n$. Under the fixed-$b$ asymptotics, the limiting null distribution of  the p-value  is the distribution of $J(b)$, where
\begin{eqnarray*}
J(b):=(1-b)^{-1} \int_0^{1-b} {\bf 1} \left\{ \sup_{s\in\R} |K(s,r+b) - K(s,r) - bK(s,1)|/\sqrt{b} \ge  \sup_{s\in \R} |K(s,1)|\right\} dr.
\end{eqnarray*}
Note that the distribution of $J(b)$ is not pivotal for a given $b$, because it depends upon the Gaussian process $K(s,t)$, whose covariance structure
 is tied to the unknown dependence structure of $X_t$. So subsampling at the first stage is insufficient under the fixed-$b$ asymptotic framework. It is worth noting that in the iid setting, the
 quantity $\sqrt{n} \|m_n-m\|_{\infty}$ is pivotal provided that $m$ is continuous, so the inferential difficulty is mainly caused by the presence of unknown
 weak dependence.

 To make the inference feasible, we propose to approximate the sampling distribution of the p-value or its limiting null distribution by subsampling; see Section~\ref{subsec:subsampling}.
 Let $n'$ be the subsampling window size at the second stage, $l'=\lceil n'b\rceil$ and $N'=n'-l'+1$.
 For each subsample $\{X_{t},\cdots,X_{t+n'-1}\}$,
 the subsampling counterpart of $pval_{n,l}^{\mbox{E}}$ is defined as
 \[h_{n',t}=(N')^{-1}\sum_{j=t}^{t+N'-1} {\bf 1}\left\{  \sqrt{l'}\|m_{j,j+l'-1}-m_{t,t+n'-1}\|_{\infty} \ge \sqrt{n'}\|m_{t,t+n'-1}-m_n\|_{\infty} \right\}\]
 for $t=1,\cdots,n-n'+1$.
Then we can approximate the sampling distribution of $pval_{n,l}^{\mbox{E}}$ or its limit null distribution $J(b)$
 by the empirical distribution associated with $\{h_{n',t}\}_{t=1}^{n-n'+1}$, denoted as $J_{n,n'}(x)=(n-n'+1)^{-1}\sum_{t=1}^{n-n'+1}{\bf 1}(h_{n',t}\le x)$.
 For a given $\alpha\in (0,1)$, the $100(1-\alpha)\%$ traditional subsampling-based confidence band for $m(\cdot)$ is
 \begin{eqnarray*}
\label{eq:band0}
\{m:~m~\mbox{is a distribution function and}~pval_{n,l}^{\mbox{E}}~\mbox{in} ~(\ref{eq:pvalueband}) \ge \alpha\},
\end{eqnarray*}
and the calibrated confidence band is
\begin{eqnarray}
\label{eq:band1}
\{m:~m~\mbox{is a distribution function and}~pval_{n,l}^{\mbox{E}}~\mbox{in} ~(\ref{eq:pvalueband}) \ge \bar{c}_{n,n',l}(\alpha)\},
\end{eqnarray}
where $\bar{c}_{n,n',l}(1-\alpha)=\inf\{x: J_{n,n'}(x) \ge 1-\alpha\}$.
The following theorem states the consistency of
subsampling at the second stage, which implies that the  coverage for the calibrated
confidence band is asymptotically correct. Let
$$\widetilde{V}_{b}(r,\epsilon):=P\left\{\left|\sup_{s\in\R} |K(s,r+b) - K(s,r) - bK(s,1)|/\sqrt{b}  - \sup_{s\in \R} |K(s,1)|\right| = \epsilon\right\}.$$
\begin{theorem}
\label{th:theorem2}
Assume that $1/n'+n'/n=o(1)$, (\ref{eq:weakconv}) and $b\in (0,1]$ is fixed. (a) The limiting null distribution of the p-value in (\ref{eq:pvalueband}) is the distribution of $J(b)$ provided that
$\widetilde{V}_b(r,0)=0$ for every $r\in [0,1-b]$. (b) Suppose that the process $X_t$ is $\alpha$-mixing, $J(b)$ is a continuous random variable and
$\widetilde{V}_b(r,\epsilon)=0$ for every $r\in [0,1-b]$ and $\epsilon\ge 0$. Then we have $$\sup_{x\in\R} |J_{n,n'}(x)-P(J(b)\le x)|=o_p(1).$$ Consequently, the asymptotic coverage probability of
confidence band in~(\ref{eq:band1}) is $(1-\alpha)$.
\end{theorem}

The conditions on $J(b)$ and $\widetilde{V}_b(r,\epsilon)$ are technical ones that are not easy to verify. The verification seems related to
the regularity of the distribution of the maximum of Gaussian processes; see Diebolt and Posse (1996),  Aza\"is and Wschebor (2001) and references therein.   We conjecture that they hold
for a large class of Gaussian processes. Note that our calibration is based on the subsampling based approximation to sampling distribution of  the p-value, which is obtained by doing the subsampling in the first stage.
As the p-value is a prepivoted statistic, we are effectively combining the prepivoting idea with subsampling in the infinite
 dimensional parameter case, for which the usual studentizing technique in the finite dimensional parameter case does
 not seem to apply. The idea of prepivoting (using the p-value) in the infinite dimensional parameter case seems new and quite general.
 We can also use the moving block bootstrap in the first stage to obtain a p-value or in the second stage to approximate the sampling distribution of the
 p-value. But the implementation of the moving block bootstrap in this setting seems very computationally demanding, especially when the block size is chosen through some data
driven algorithms. For this reason, we shall focus on the subsampling method in simulation studies for the infinite dimensional case.

\subsection{Spectral distribution function}
\label{subsec:sdf}

Another infinite dimensional parameter of interest in time series analysis is the spectral distribution function $F(\lambda)=\int_0^{\lambda} f(w)dw$,  $\lambda\in [0,\pi]$, where
$f(\cdot)$ is the spectral density function of $\{X_t\}$. Let $I_n(w)=(2\pi n)^{-1}\left|\sum_{t=1}^{n} (X_t-\bar{X}_n) e^{itw}\right|^2$
be the periodogram. A commonly used estimator for $F(\lambda)$ is $F_n(\lambda)=\int_0^{\lambda}I_n(w)dw$ or its discretized version $F_n(\lambda)=(2\pi)/n \sum_{0<2\pi s/n\le \lambda} I_n(2\pi s/n)$. It has been shown that the two versions are asymptotically equivalent [Dahlhaus (1985a)]
 and we shall use the discrete version for the computational convenience. Under certain regularity conditions,
we have that $\sqrt{n}\{F_n(\lambda)-F(\lambda)\}\Rightarrow G(\lambda)$,
where $G(\lambda)$ is a mean zero Gaussian process with covariance
\[C(\lambda,\lambda')=\cov\{G(\lambda),G(\lambda')\}=2\pi \int_0^{\lambda\wedge\lambda'} f^2(w)dw+ 2\pi \int_0^{\lambda}\int_0^{\lambda'} f_4(w_1,-w_1,-w_2)dw_1 dw_2.\]
 Here $f_4(\cdot,\cdot,\cdot)$ is the fourth order cumulant spectrum. For various sets of conditions for this weak convergence to hold, see
 Brillinger (1975), Dahlhaus (1985b) and Anderson (1993).
Applying the continuous mapping theorem, we get
\begin{eqnarray}
\label{eq:limit2}
\sqrt{n}\|F_n-F\|_{\infty}\rightarrow_{D} \sup_{\lambda\in [0,\pi]} |G(\lambda)|.
\end{eqnarray}
 Since the covariance of $G(\lambda)$ depends on unknown second order and fourth order spectrum, the distribution of $\sup_{\lambda\in [0,\pi]}|G(\lambda)|$ is unknown and is usually difficult to estimate directly, which renders the confidence band construction for $F$ a hard task. To alleviate the problem,   Politis et al. (1999b) proposed to apply the subsampling method to approximate the limiting distribution in (\ref{eq:limit2}) and they proved the consistency. See Politis et al. (1993) for some related numerical work.
 Often in practice, the main interest is on the pattern of dependence described in terms of autocorrelations, then the normalized spectral distribution function $\widetilde{F}(\lambda)=F(\lambda)/F(\pi)$, $\lambda\in [0,\pi]$, which does not depend on the marginal variance of $X_t$, is  of more practical relevance. Politis et al. (1999b) mentioned that the subsampling method is still consistent in the approximation of the sampling distribution or the limiting distribution of $\sqrt{n}\|\widetilde{F}_n-\widetilde{F}\|_{\infty}$, where $\widetilde{F}_n(\lambda)=F_n(\lambda)/F_n(\pi)$.

To introduce our calibration method, we need to define the estimate of $\widetilde{F}(\lambda)$ based on the subsample
$(X_t,\cdots,X_{t'})$ for $1\le t<t'\le n$. In particular, we define the periodogram on the basis of the subsample
$\{X_t,\cdots,X_{t'}\}$ as $I_{t,t'}(w)=\{2\pi (t'-t+1)\}^{-1}|\sum_{j=t}^{t'}(X_j-\bar{X}_{t,t'})\exp(ijw)|^2$, where $\bar{X}_{t,t'}=(t'-t+1)^{-1}\sum_{j=t}^{t'}X_j$,
  $F_{t,t'}(\lambda)=\int_0^{\lambda}I_{t,t'}(w) dw$, and $\widetilde{F}_{t,t'}(\lambda)=F_{t,t'}(\lambda)/F_{t,t'}(\pi)$. The subsampling method approximates the sampling  distribution of $\sqrt{n} \|\widetilde{F}_n-\widetilde{F}\|_{\infty}$ by
  the empirical distribution generated from  $ \sqrt{l}\|\widetilde{F}_{t,t+l-1}-\widetilde{F}_{1,n}\|_{\infty}$, $t=1,\cdots,N$
and the corresponding p-value is
\begin{eqnarray}
\label{eq:pvalue2}
pval_{n,l}^{\mbox{S}}=N^{-1}\sum_{t=1}^{N} {\bf 1}\left\{\sqrt{l}\|\widetilde{F}_{t,t+l-1}-\widetilde{F}_{1,n}\|_{\infty} \ge  \sqrt{n}\|\widetilde{F}_n-\widetilde{F}\|_{\infty}\right\}.
\end{eqnarray}
Under the small-b asymptotics, the limiting null distribution of $pval_{n,l}^{\mbox{S}}$ is $U(0,1)$, but under the fixed-$b$
asymptotics, the limiting null distribution is expected to depend on $b$ and the intricate second and fourth order dependence structure of $X_t$. The derivation of the limiting distribution of the p-value relies on
the functional central limit theorem for $\sqrt{n}\{F_{1,\lfloor nr\rfloor}(\lambda)-F(\lambda)\}$, $(r,\lambda)\in [0,1]\times [0,\pi]$, which
seems unavailable in the literature.  In view of Theorem 1 in Shao (2009), Theorem 2 in Shao (2010a), and Theorem 3.3 in Dahlhaus (1985b), we conjecture that
\begin{eqnarray}
\label{eq:fclt1}\sqrt{n}\{F_{1,\lfloor nr\rfloor}(\lambda)-F(\lambda)\}\Rightarrow H(r,\lambda),
\end{eqnarray}
 where $H(r,\lambda)$ is a mean zero Gaussian process with covariance
$\cov\{H(r,\lambda),H(r',\lambda')\}=(r\wedge r') C(\lambda,\lambda')$. Let $\widetilde{H}(r,\lambda)=\{H(r,\lambda) F(\pi) -F(\lambda) H(r,\pi)\}/F^2(\pi)$.
Then (\ref{eq:fclt1}) implies that $\sqrt{n}\{ \widetilde{F}_{1,\lfloor nr\rfloor}(\lambda)-\widetilde{F}(\lambda)\}\Rightarrow \widetilde{H}(r,\lambda)$ by the continuous mapping theorem and that the limiting null distribution of the p-value is the distribution of
$$(1-b)^{-1}\int_0^{1-b} {\bf 1}\left(\sup_{\lambda\in [0,\pi]} |\widetilde{H}(r+b,\lambda)-\widetilde{H}(r,\lambda)-b\widetilde{H}(1,\lambda)|/\sqrt{b}\ge \sup_{\lambda\in [0,\pi]} |\widetilde{H}(1,\lambda)|\right) dr,$$
which is not pivotal. Following the calibration idea described in Section~\ref{subsec:mdf},  we apply the subsampling method to approximate the sampling distribution of the p-value. For a given subsampling block size $n'$ at the second stage, let $l'=\max(\lceil n' b\rceil,2)$ since a minimum sample size of 2 is needed to estimate the spectral distribution function. Let $N'=n'-l'+1$ and
\[\widetilde{h}_{n',t}=(N')^{-1}\sum_{j=t}^{t+N'-1} {\bf 1}\left\{ \sqrt{l'} \|\widetilde{F}_{j,j+l'-1}-\widetilde{F}_{t,t+n'-1}\|_{\infty} \ge \sqrt{n'}\|\widetilde{F}_{t,t+n'-1}-\widetilde{F}_n\|_{\infty} \right\}\]
for $t=1,\cdots,n-n'+1$. The traditional subsampling-based $100(1-\alpha)\%$ confidence band for $\widetilde{F}(\cdot)$ is
\begin{eqnarray*}
\label{eq:bandspec0}
\{\widetilde{F}: \widetilde{F}~\mbox{is a normalized spectral distribution function and}~pval_{n,l}^{\mbox{S}}~\mbox{in~(\ref{eq:pvalue2})}\ge \alpha\}.
\end{eqnarray*}
In contrast, the calibrated confidence band is
\begin{eqnarray}
\label{eq:bandspec1}
\{\widetilde{F}: \widetilde{F}~\mbox{is a normalized spectral distribution function and}~pval_{n,l}^{\mbox{S}}~\mbox{in~(\ref{eq:pvalue2})}\ge \tilde{c}_{n,n',l}(\alpha)\},
\end{eqnarray}
 where $\tilde{c}_{n,n',l}(1-\alpha)=\inf\{x: (n-n'+1)^{-1}\sum_{t=1}^{n-n'+1} {\bf 1}(\widetilde{h}_{n',t}\le x)\ge 1-\alpha\}$. If (\ref{eq:fclt1}) is true, then the confidence band in (\ref{eq:bandspec1})
 is expected to have $100(1-\alpha)\%$ coverage asymptotically under
 appropriate mixing and moment conditions and the assumptions
 that $b\in (0,1]$ is held fixed and $1/n'+n'/n\rightarrow 0$ as $n\rightarrow\infty$.

\section{Simulation results}
\label{sec:simulation}

\bigskip

In this section, we conduct simulation studies to evaluate the accuracy of the asymptotic approximations provided
by both small-$b$ and fixed-$b$ approaches to the finite sample distribution. Specifically, we examine the empirical coverage
probabilities and the volumes of confidence sets to see if the fixed-$b$ approach corresponds to smaller coverage errors.

\subsection{Finite sample performance of confidence intervals}
\label{subsec:interval}

In this subsection, we consider
a univariate stationary time series model with various types of  dependence structure. To be specific, we let
\[X_t=\mu+u_t,~~u_t=\rho u_{t-1}+\epsilon_t+\theta \epsilon_{t-1},~\epsilon_t\sim iid~ N(0,1).\]
We consider (i) AR(1)-$N(0,1)$ error: $(\rho,\theta)=(0,0)$, $(0.5,0)$ and $(0.8,0)$; (ii) MA(1)-$N(0,1)$ error: $(\rho,\theta)=(0,-0.5)$;
 and their corresponding AR(1)-EXP(1) and MA(1)-EXP(1) models, where $\epsilon_t\sim iid~\mbox{EXP}(1)-1$ has mean zero, unit variance but with an asymmetric distribution.
Following the suggestion of a referee, we also include two nonlinear time series models: Nonlinear 1, $X_t=0.6\sin(X_{t-1})+\epsilon_t$, where $\epsilon_t\sim iid~N(0,1)$. This model was used in the simulation work of
 Paparoditis and Politis (2001) and Shao (2010c); Nonlinear 2 (threshold autoregressive model of order 1), $X_t=0.3X_{t-1}{\bf 1}(X_{t-1}>0)+0.8X_{t-1}{\bf 1}(X_{t-1}\le 0)+\epsilon_t$, where $\epsilon_t\sim iid~N(0,1)$. Sample size $n=100$ and the number of bootstrap replications is $5000$. The bandwidth parameter $l$ varies from $3$ to $16$, i.e. $b=0.03,0.04,\cdots,0.16$. We examine the empirical coverages and average widths of symmetric confidence intervals for
$\mu=\E(X_1)$ and the $25\%$ trimmed mean based on 10000 replications.  Nominal coverage is set to be $95\%$.

\begin{figure}
\caption{The empirical coverage probabilities (left panel) and  the ratios of interval widths (calibrated fixed-$b$  over traditional small-$b$)
 (right panel) for the mean and for the models with normally distributed errors. Sample size $n=100$ and number of replications is 10000. }
\begin{center}
{\includegraphics[height=8cm,width=4.5cm,angle=270]{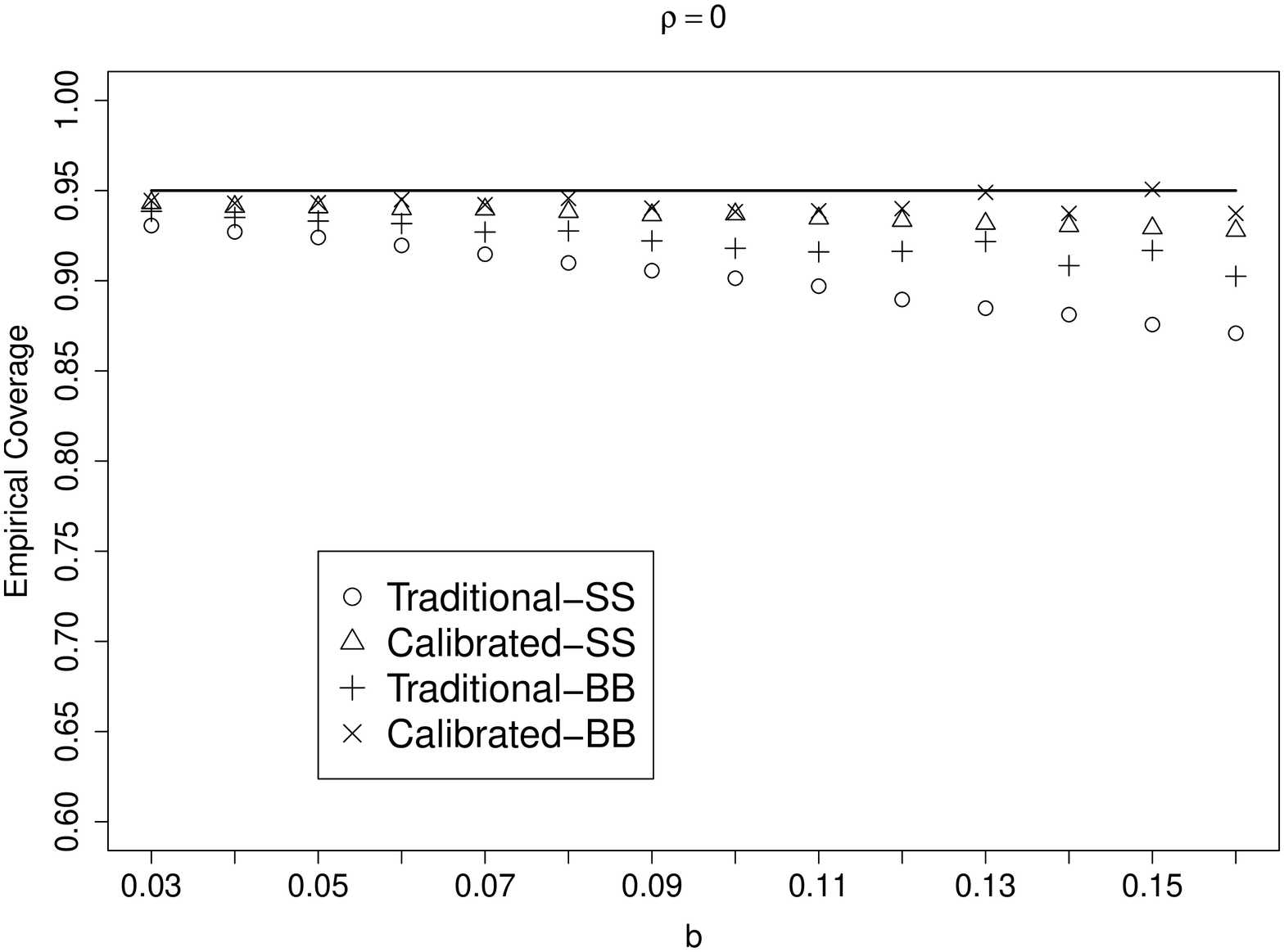}}
{\includegraphics[height=8cm,width=4.5cm,angle=270]{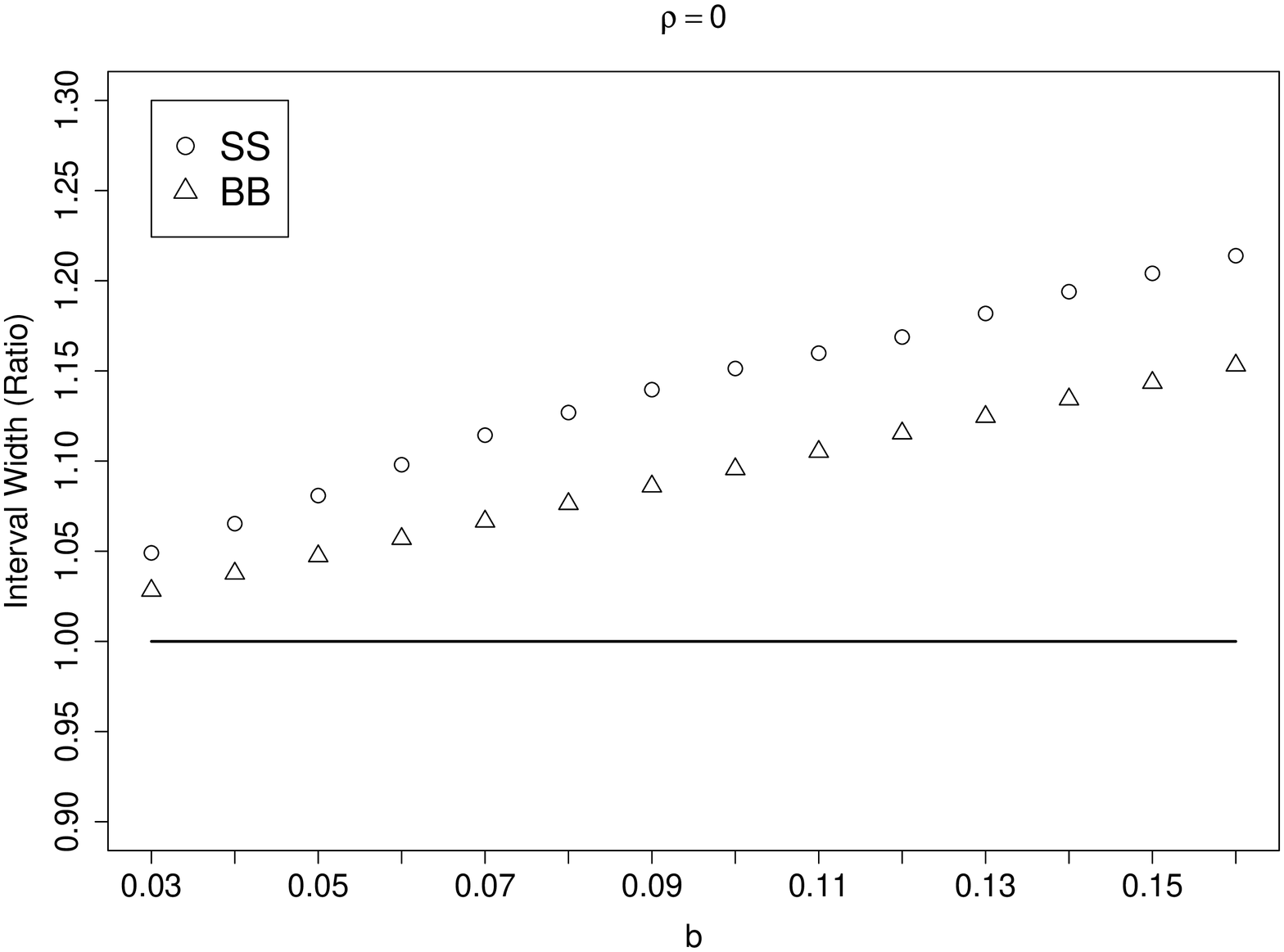}}
{\includegraphics[height=8cm,width=4.5cm,angle=270]{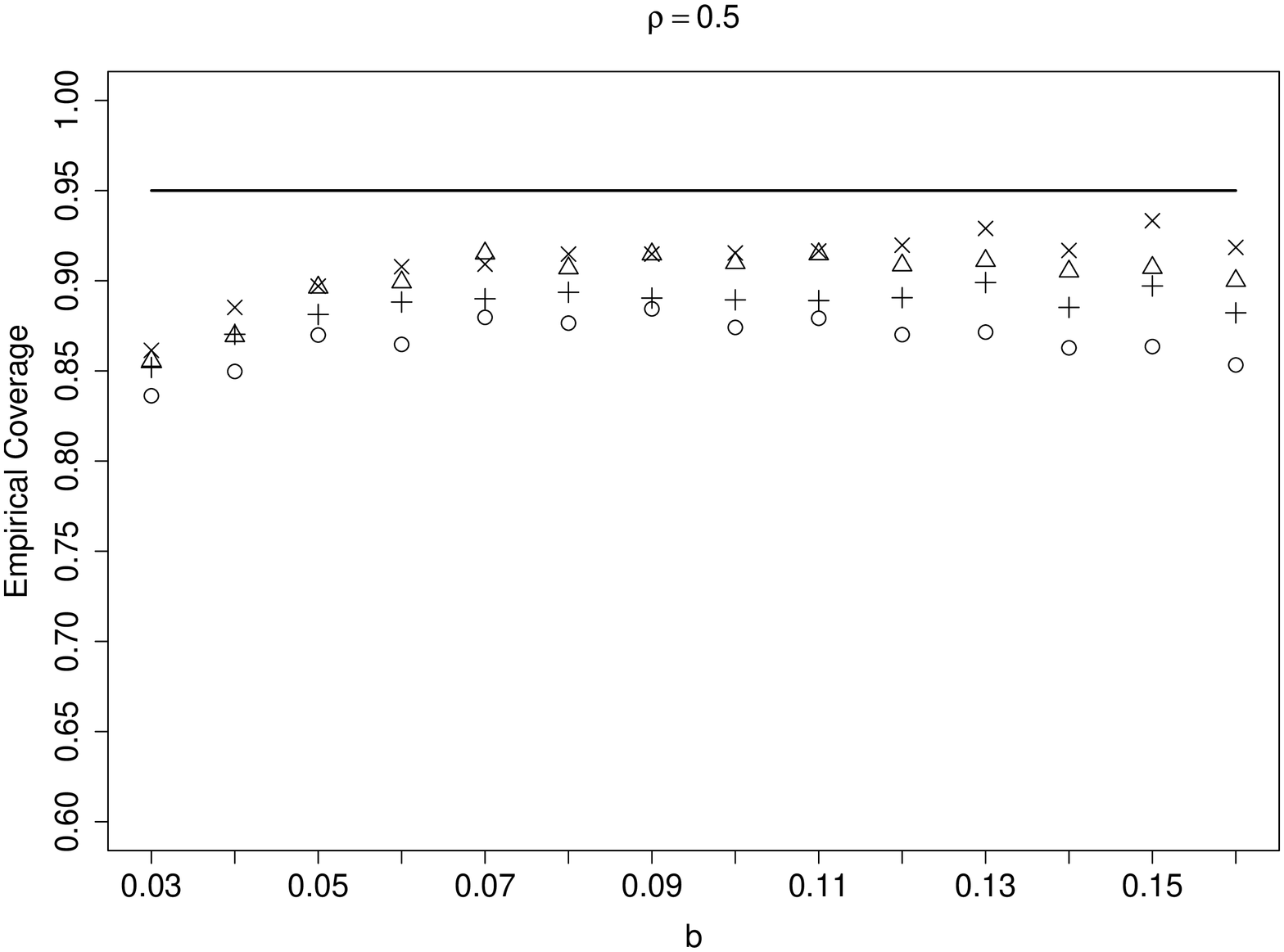}}
{\includegraphics[height=8cm,width=4.5cm,angle=270]{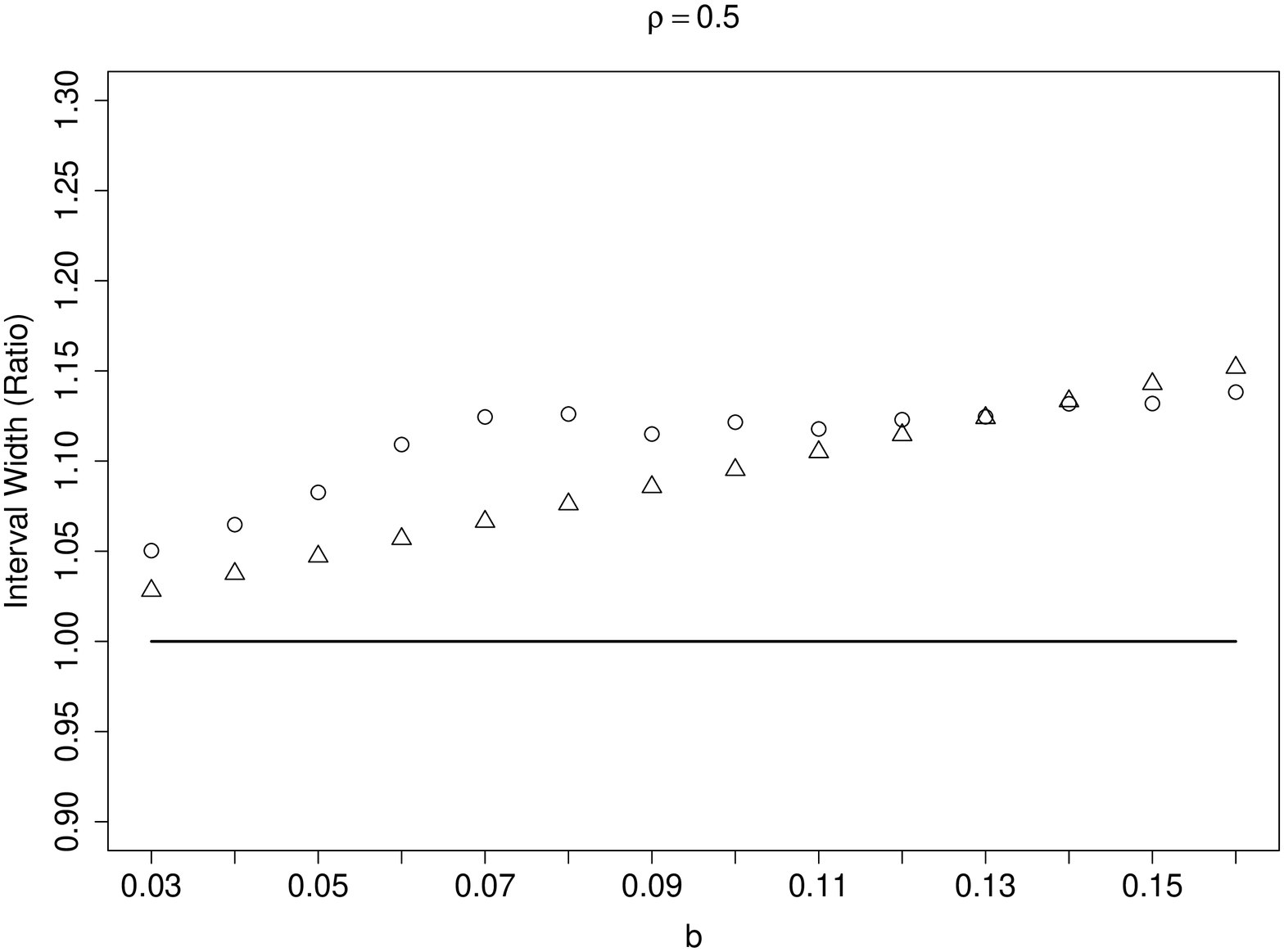}}
{\includegraphics[height=8cm,width=4.5cm,angle=270]{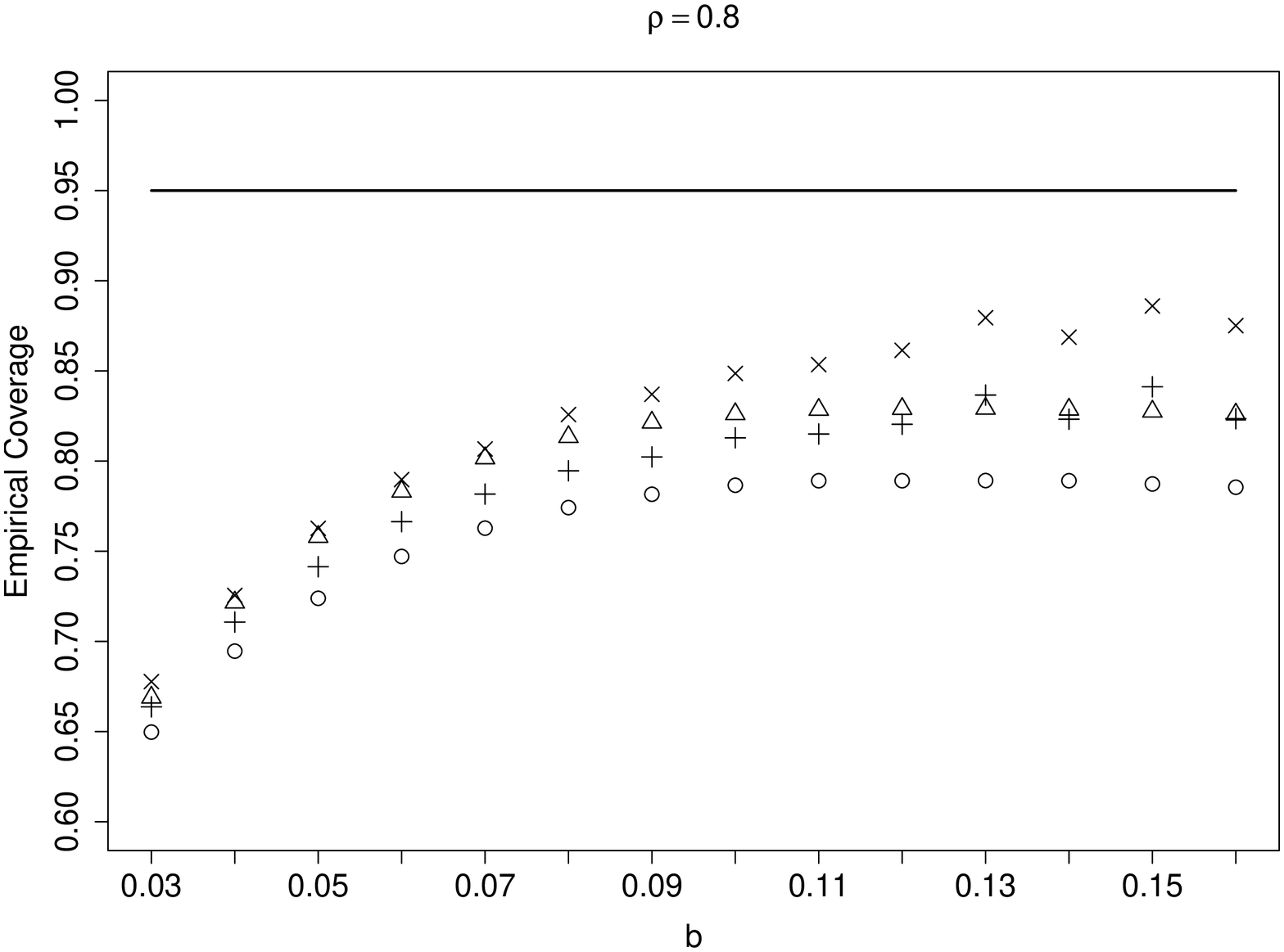}}
{\includegraphics[height=8cm,width=4.5cm,angle=270]{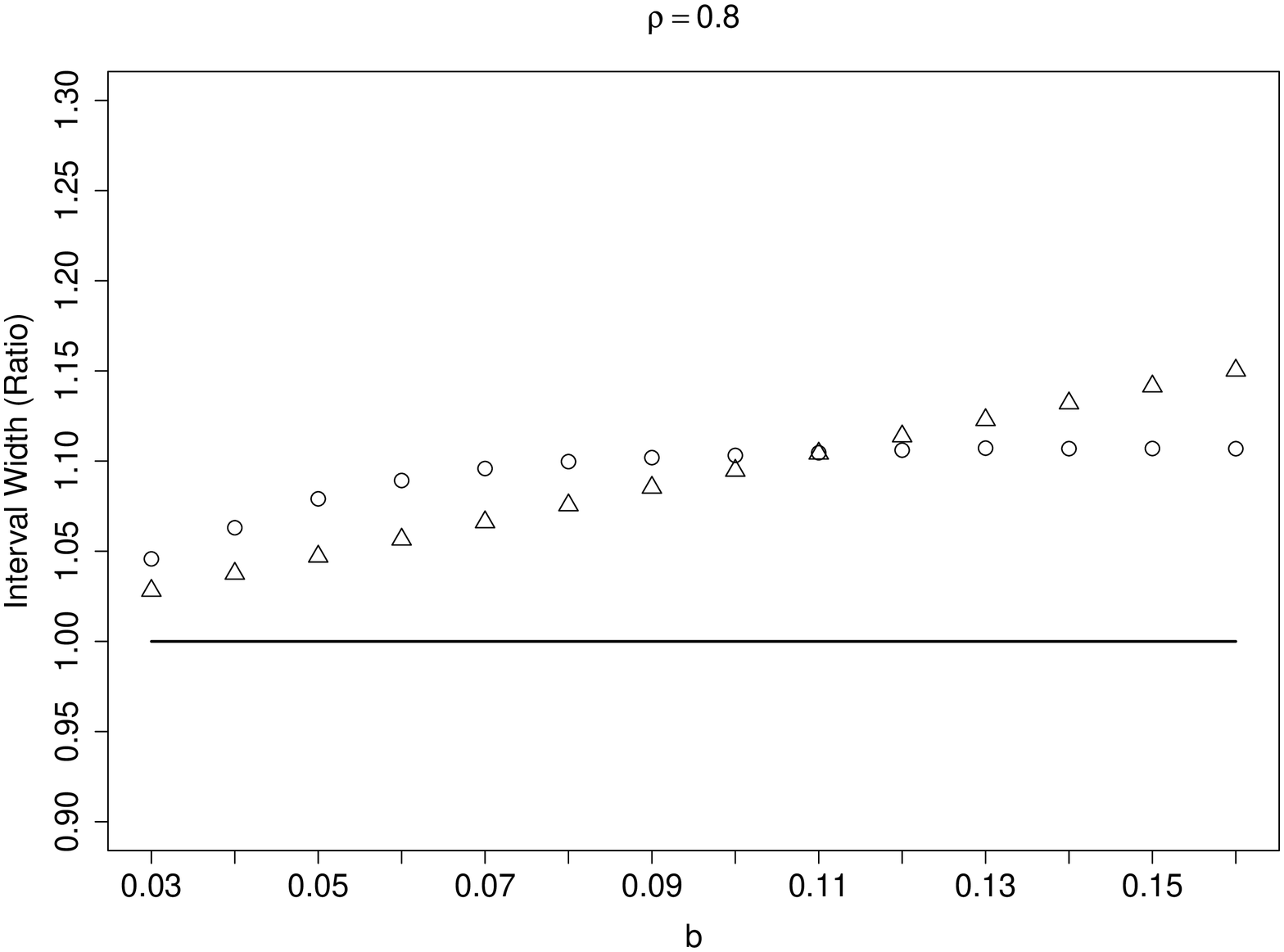}}
{\includegraphics[height=8cm,width=4.5cm,angle=270]{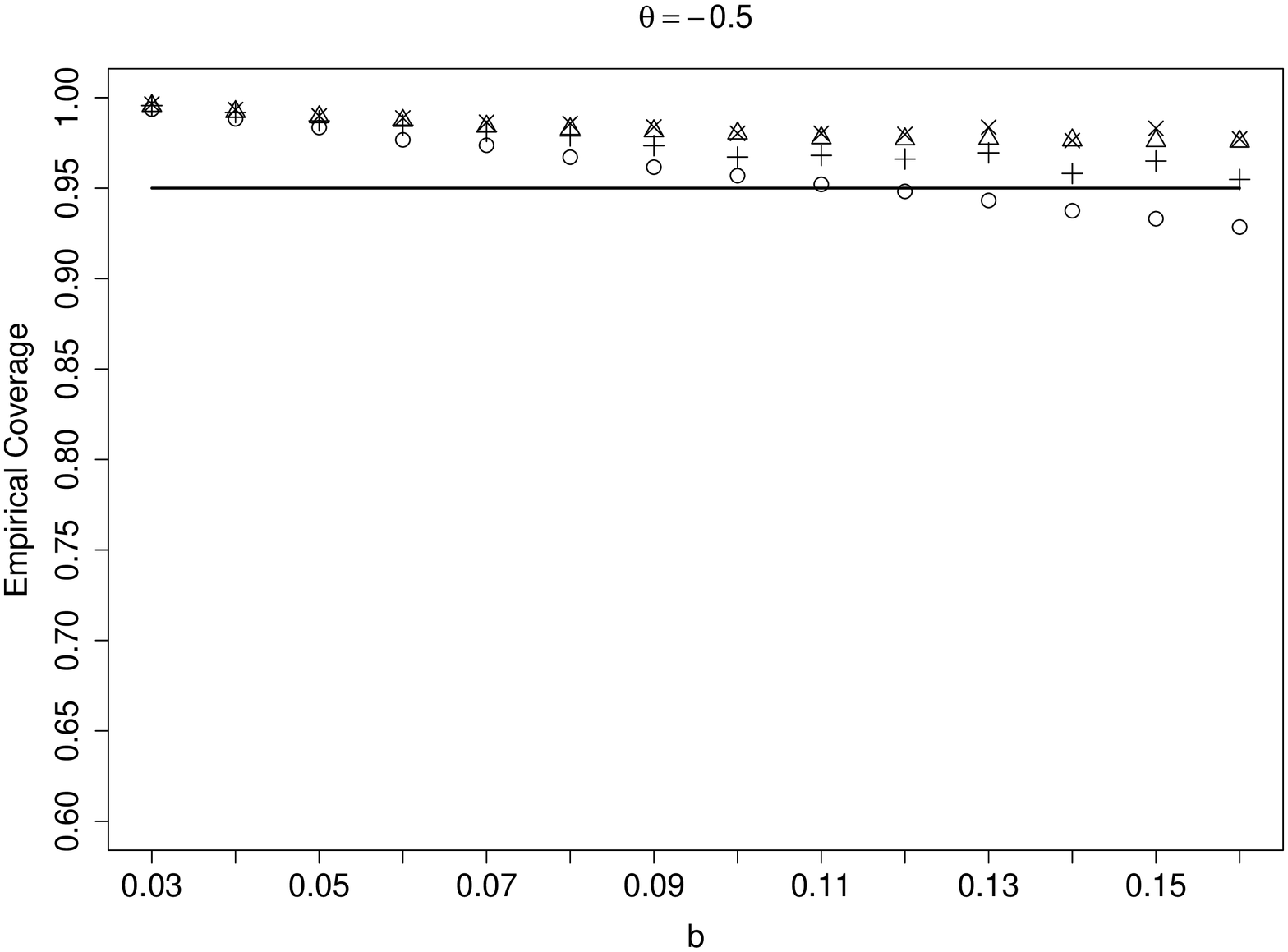}}
{\includegraphics[height=8cm,width=4.5cm,angle=270]{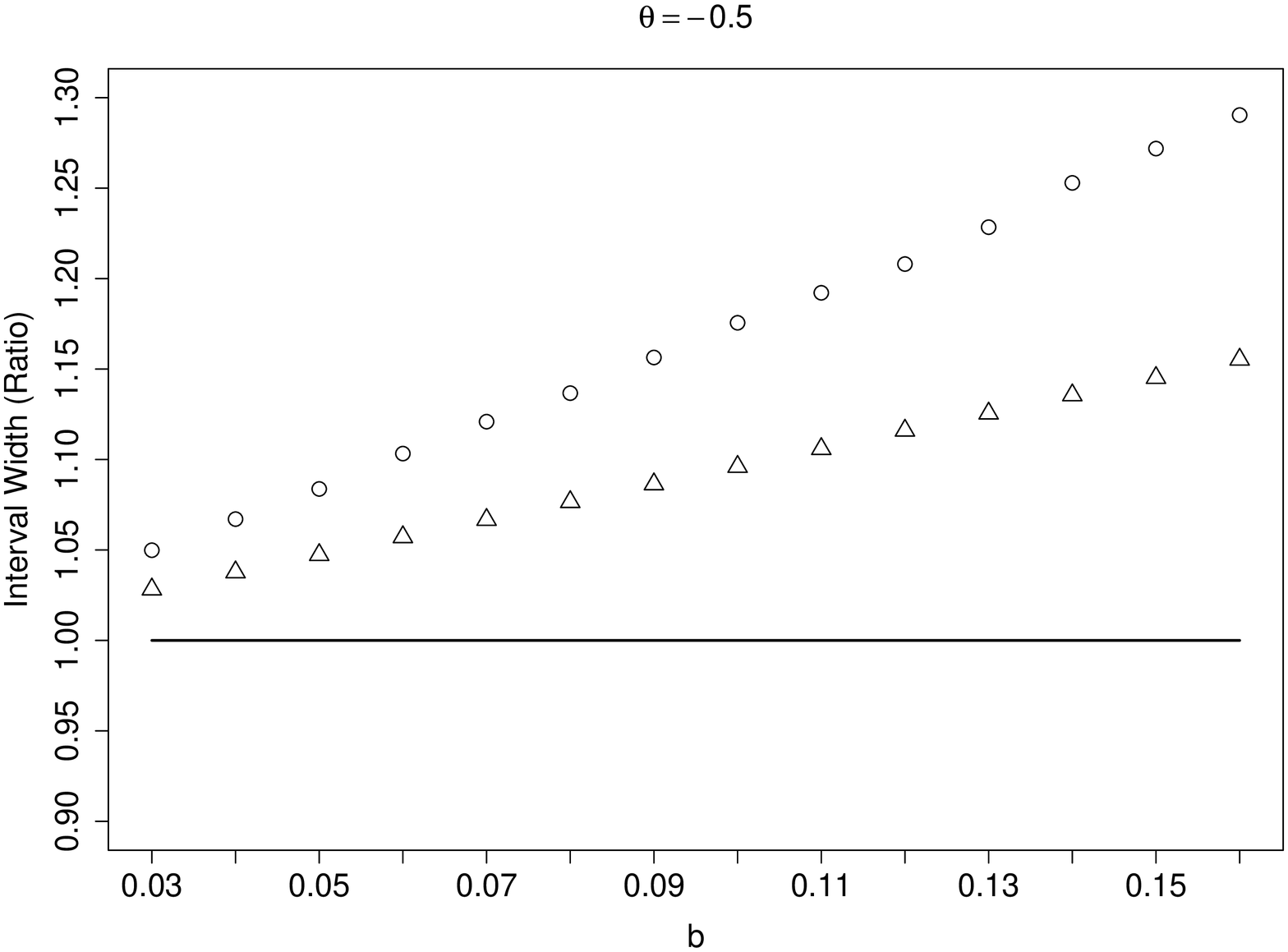}}
\label{fig:mean1}
\end{center}

\end{figure}

For the models with normally distributed errors, the results for the mean case are depicted in Figure~\ref{fig:mean1}, in which the left panel shows the empirical
coverages and the right panel shows the corresponding ratios of average interval widths (fixed-$b$ over small-$b$). The symbols ``SS" and ``BB" stand for  subsampling and the moving block bootstrap, respectively.  For both subsampling and the moving block bootstrap, the undercoverage occurs and it gets more severe as the dependence strengthens. The empirical coverages for the fixed-$b$ approach are closer to the nominal level than those for the small-$b$ approach, with the difference between two empirical coverages
 increasing as $b$ gets large. On the other hand,  the  fixed-$b$ based interval is slightly wider than its small-$b$ counterpart, with the ratio of widths increasing with respect to $b$
 in general. These findings are consistent with the intuition that the larger $b$ is, the more accurate the fixed-$b$ based approximation provides relative to its small-$b$ counterpart.
  The intervals constructed based on the moving block bootstrap have noticeably better coverage than the ones based on subsampling, especially for large $b$.
 For the MA(1) model with $\theta=-0.5$, there is an overcoverage for the fixed-$b$ based
 interval, which is usually slightly more conservative than the small-$b$ counterpart. 
 The overcoverage  in the case of negatively correlated time series corresponds to the underrejection for Kiefer and Vogelsang's (2005) studentized statistic when using normal approximation (i.e. small-$b$ approach) and $b$ is small (see Figure 1 therein), so our results are in a sense consistent with those in Kiefer and Vogelsang (2005). Practically speaking, the overcoverage is less harmful to the practitioner than the undercoverage, so are less concerned in practice. 
 
\begin{figure}
\caption{The empirical coverage probabilities (left panel) and  the ratios of interval widths (calibrated fixed-$b$  over traditional small-$b$)
 (right panel) for the mean and for the models with exponentially distributed errors. Sample size $n=100$ and number of replications is 10000. }
\begin{center}
{\includegraphics[height=8cm,width=4.5cm,angle=270]{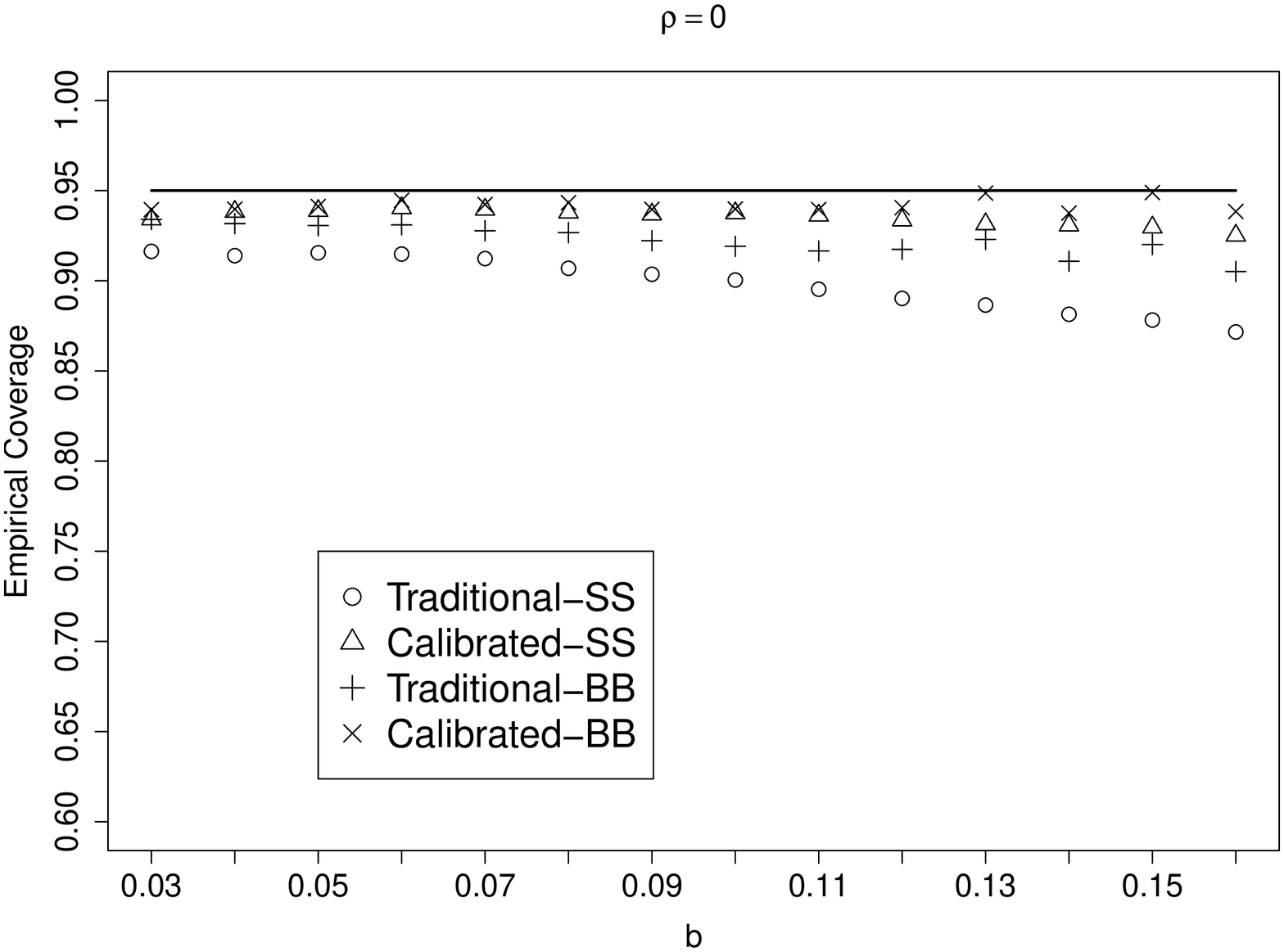}}
{\includegraphics[height=8cm,width=4.5cm,angle=270]{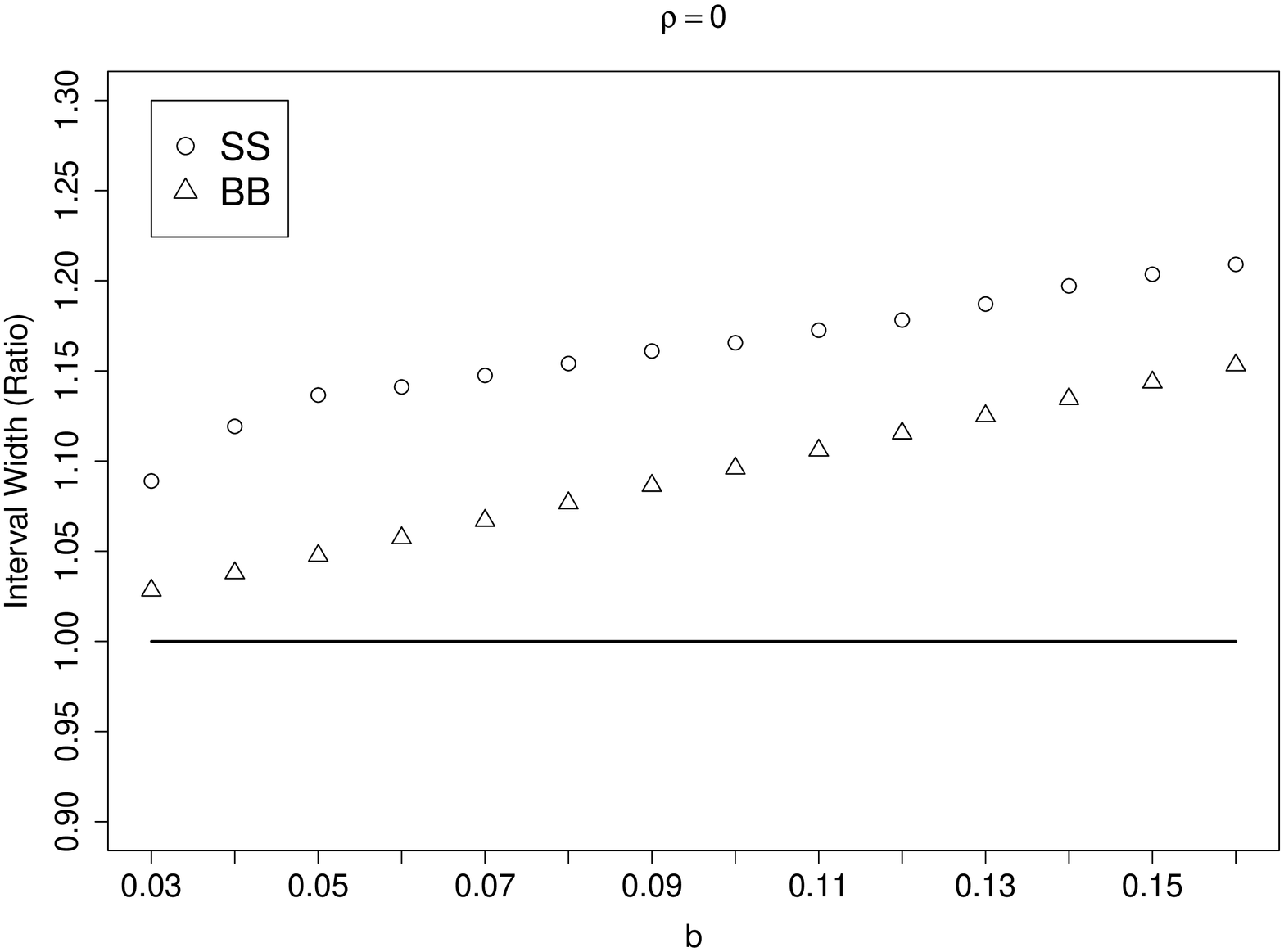}}
{\includegraphics[height=8cm,width=4.5cm,angle=270]{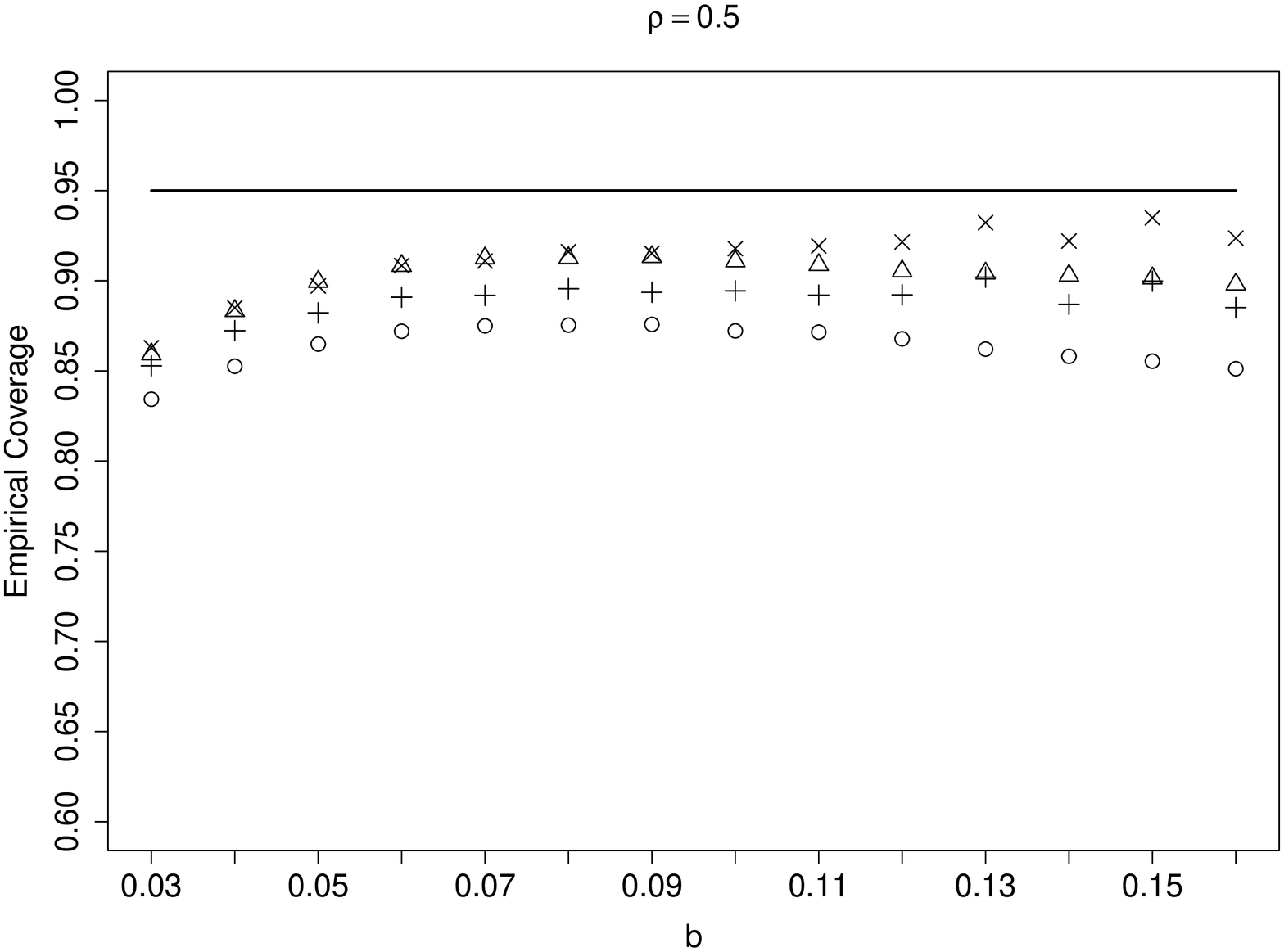}}
{\includegraphics[height=8cm,width=4.5cm,angle=270]{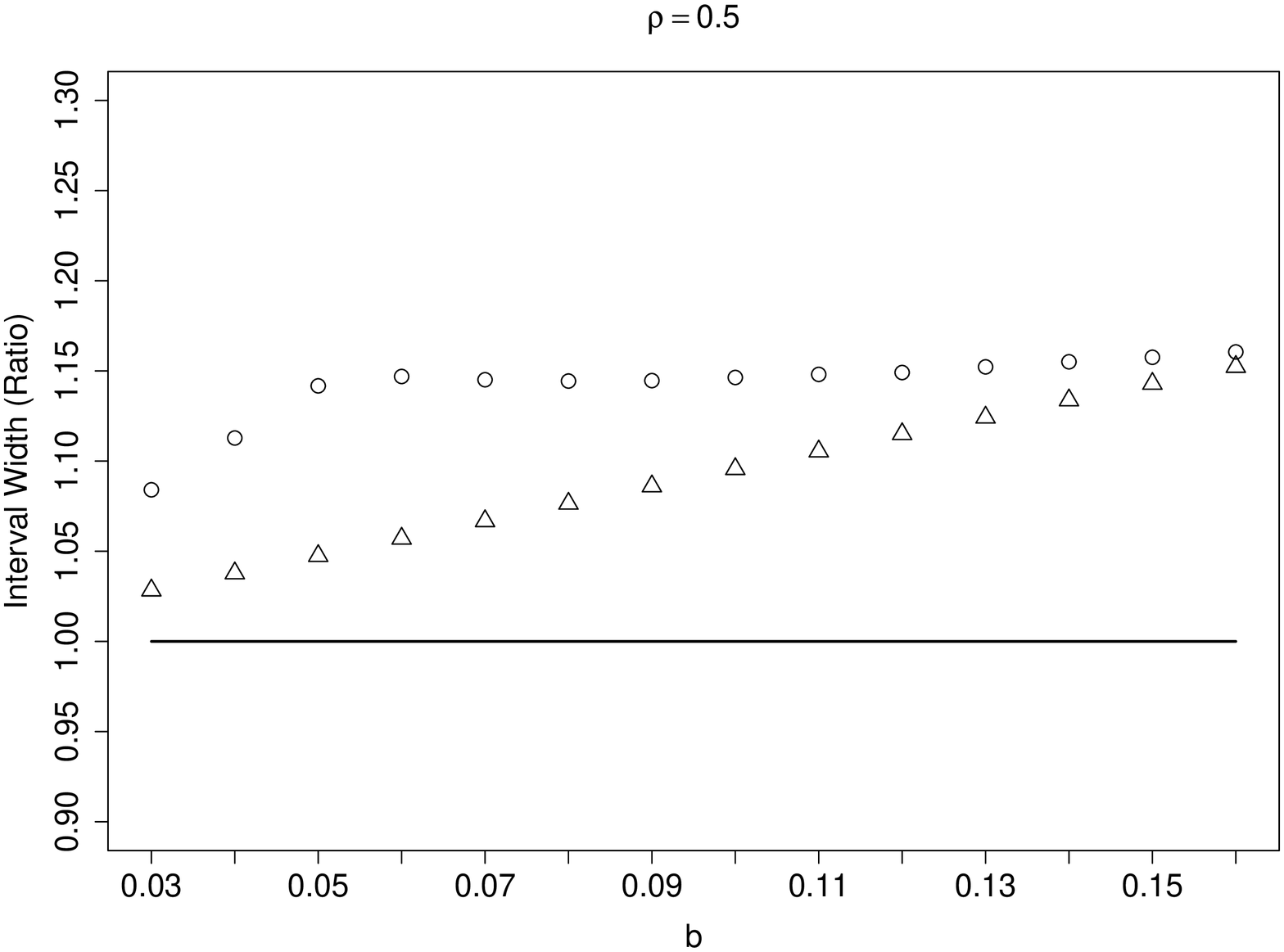}}
{\includegraphics[height=8cm,width=4.5cm,angle=270]{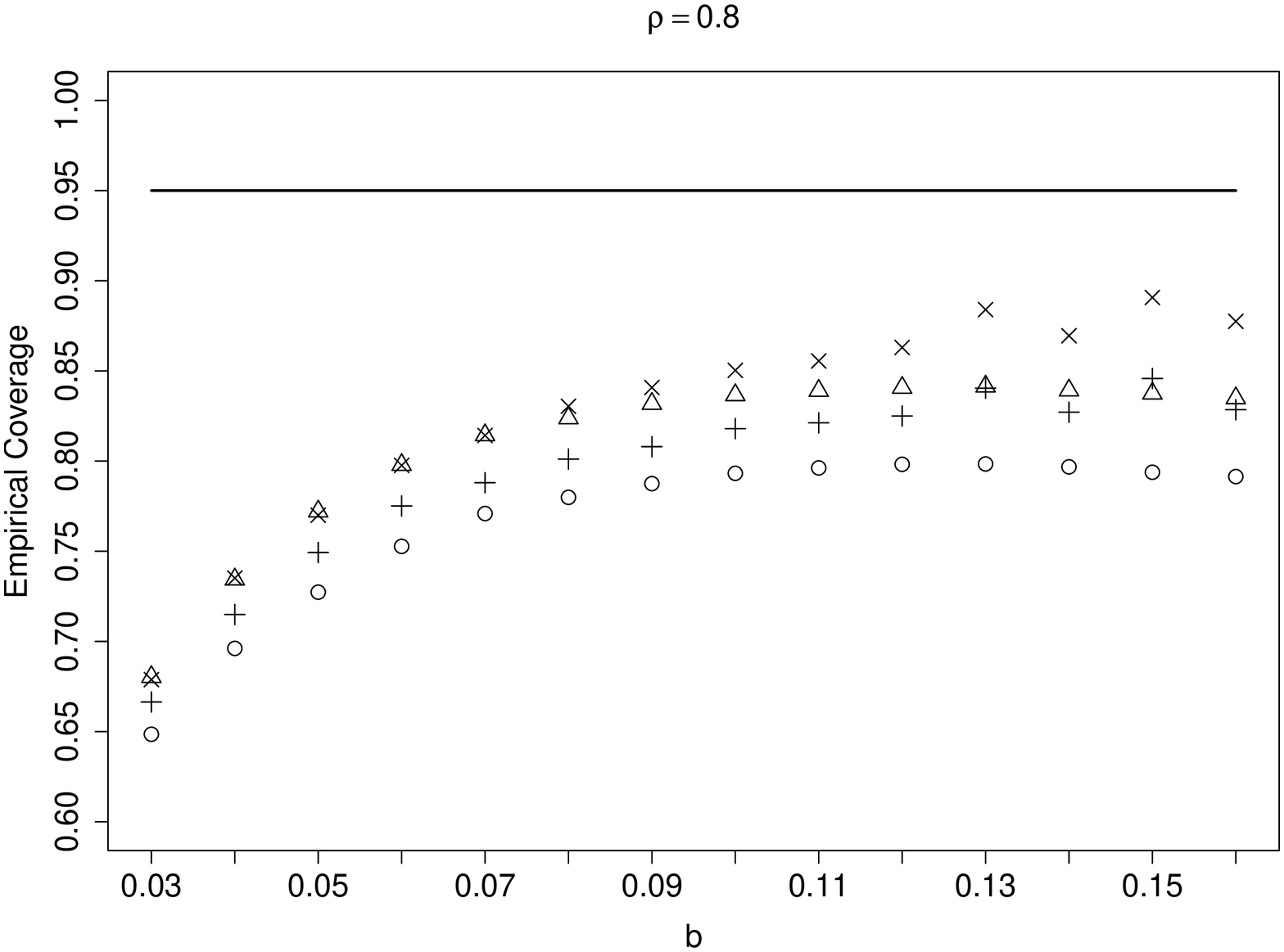}}
{\includegraphics[height=8cm,width=4.5cm,angle=270]{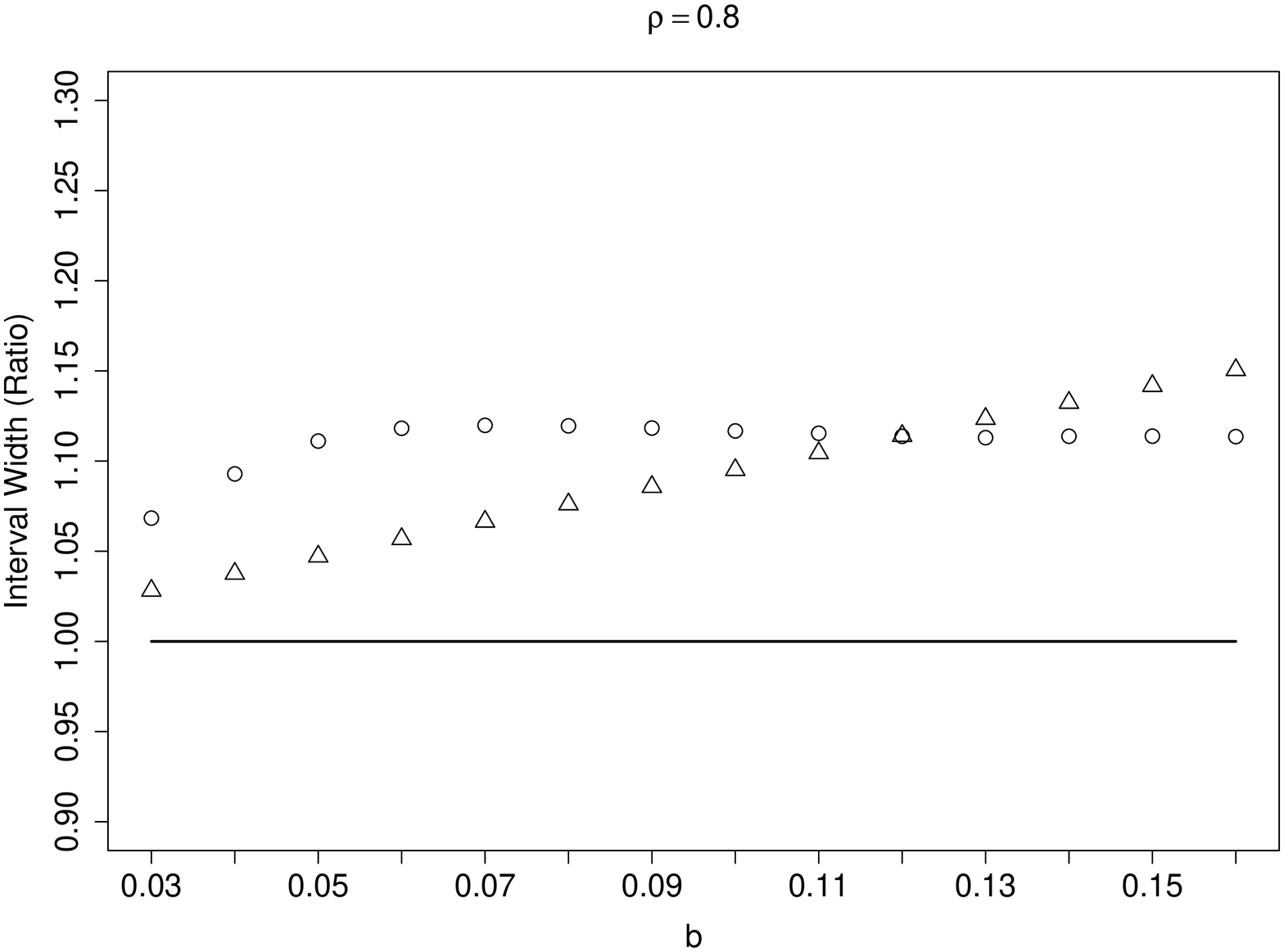}}
{\includegraphics[height=8cm,width=4.5cm,angle=270]{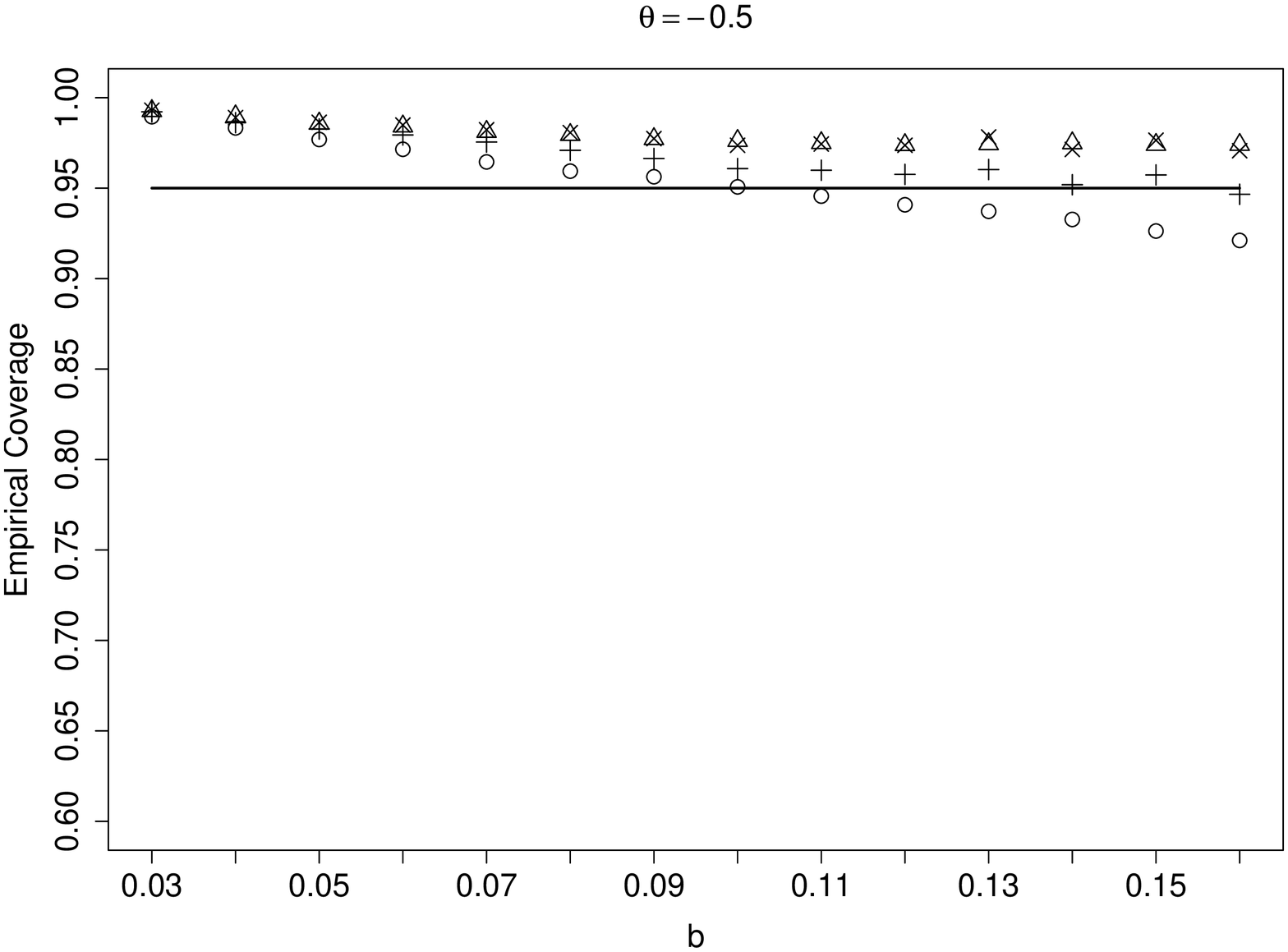}}
{\includegraphics[height=8cm,width=4.5cm,angle=270]{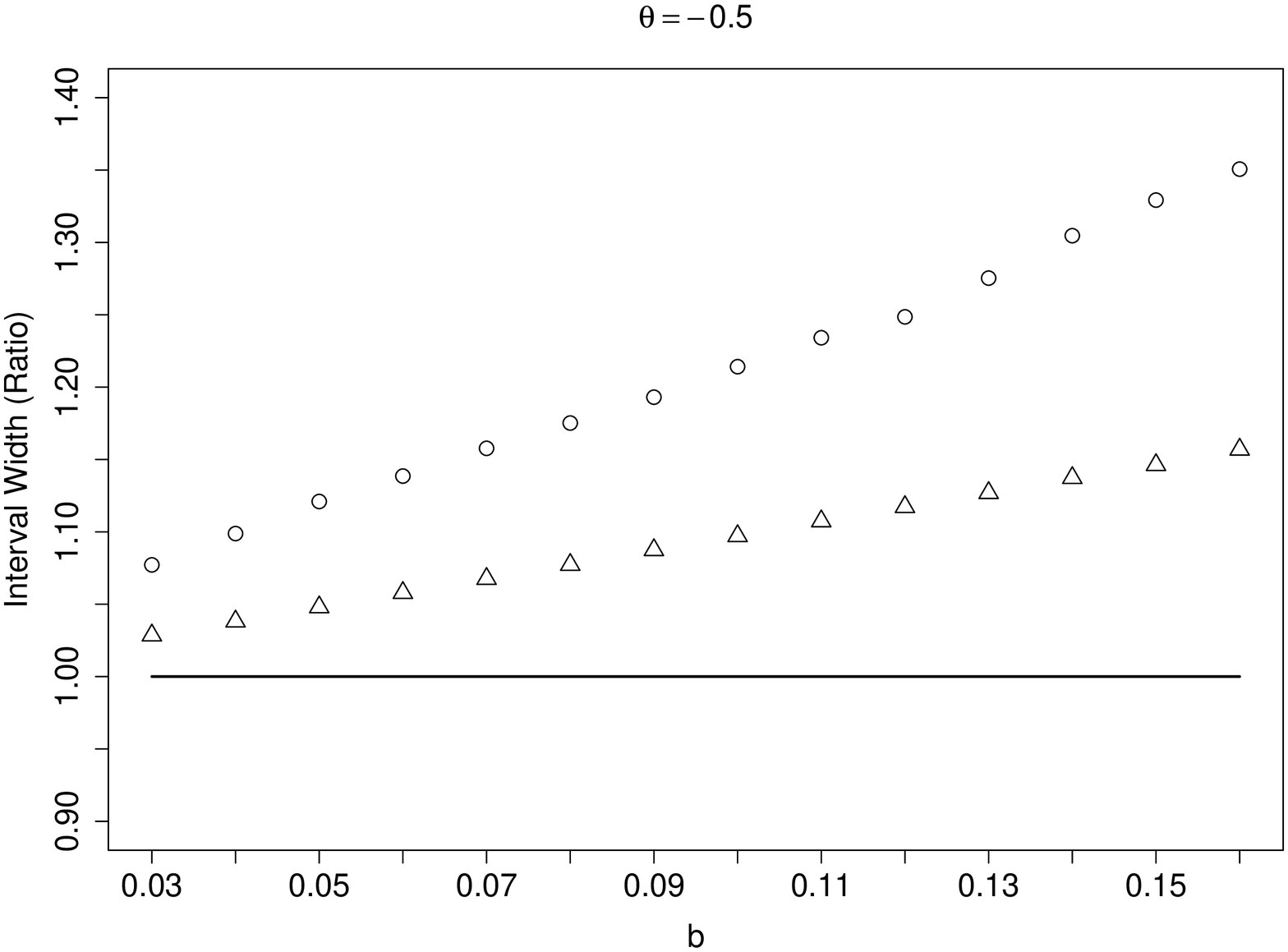}}
\label{fig:meanEXP1}
\end{center}

\end{figure}

\begin{figure}
\caption{The empirical coverage probabilities (left panel) and  the ratios of interval widths (calibrated fixed-$b$  over traditional small-$b$)
 (right panel) for the $25\%$ trimmed mean and for the models with normally distributed errors. Sample size $n=100$ and number of replications is 10000. }
\begin{center}
{\includegraphics[height=8cm,width=4.5cm,angle=270]{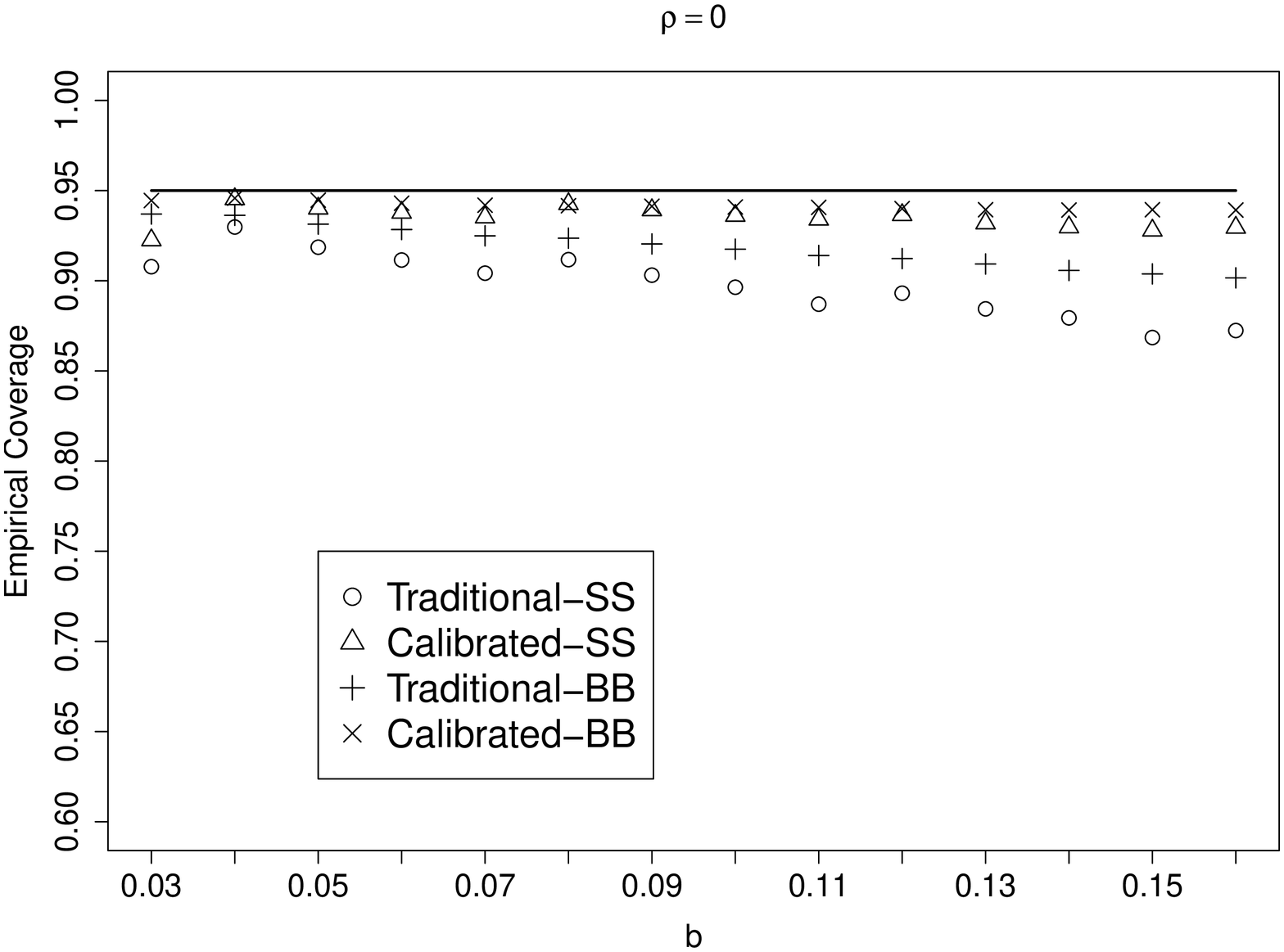}}
{\includegraphics[height=8cm,width=4.5cm,angle=270]{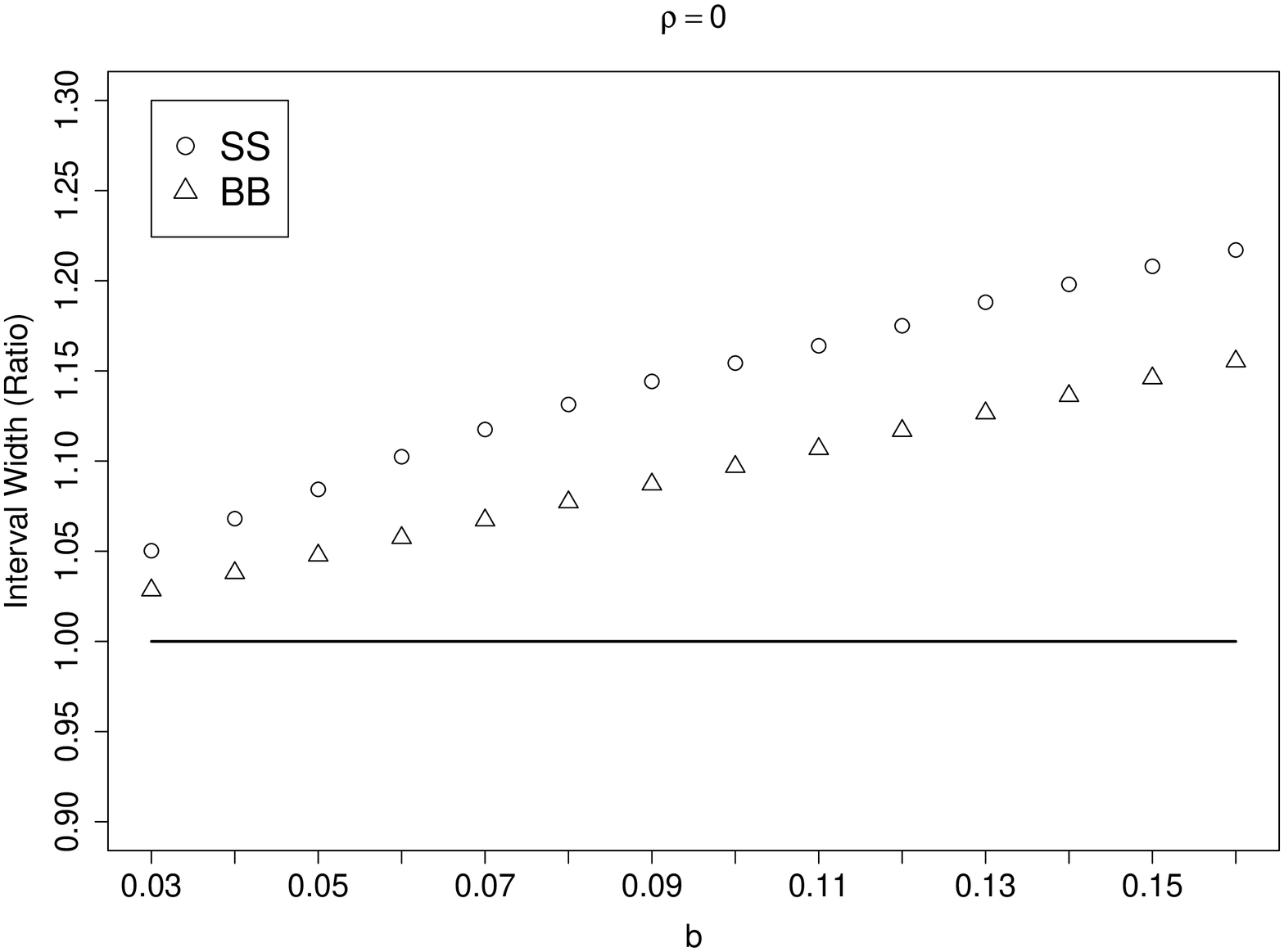}}
{\includegraphics[height=8cm,width=4.5cm,angle=270]{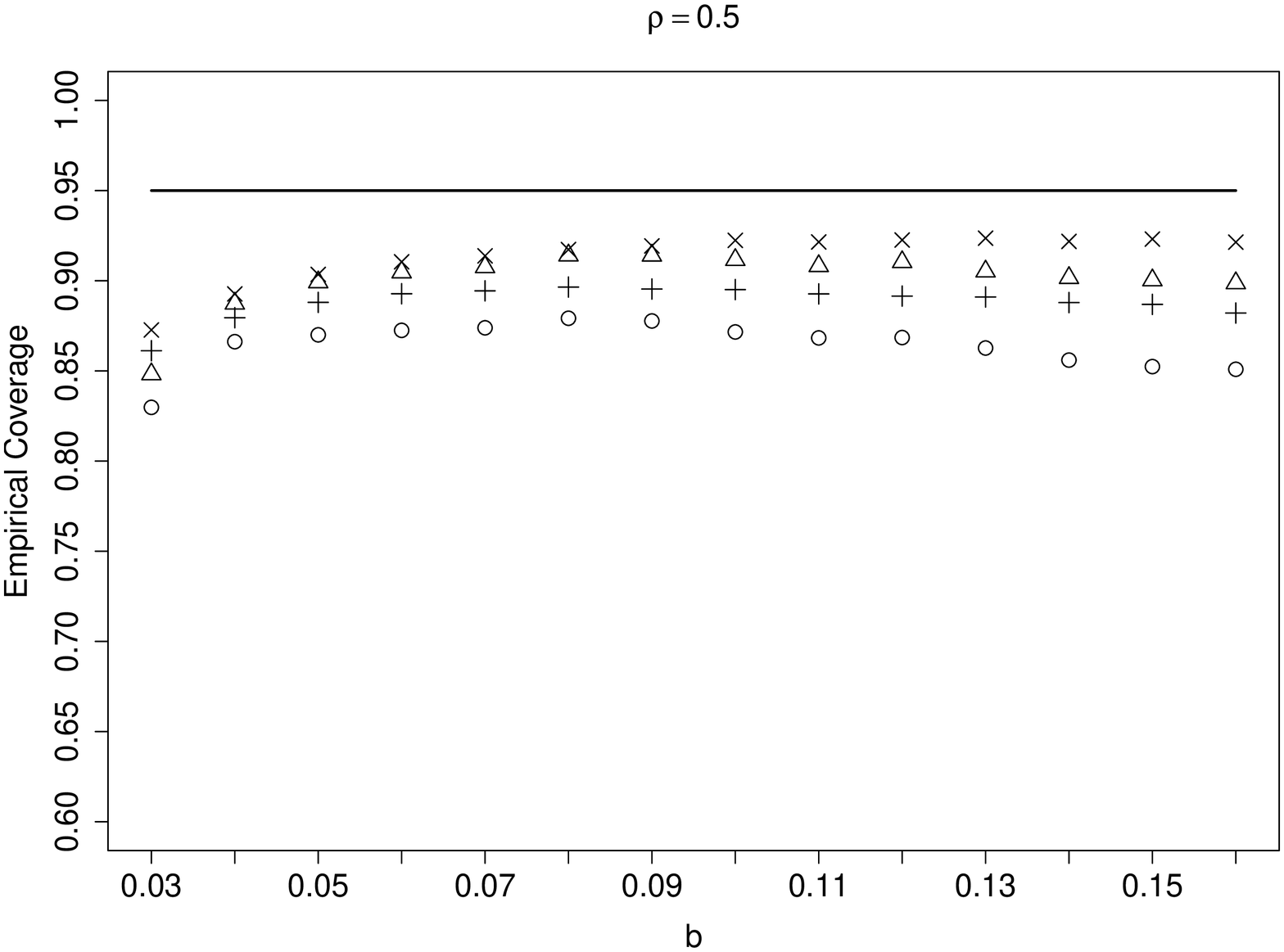}}
{\includegraphics[height=8cm,width=4.5cm,angle=270]{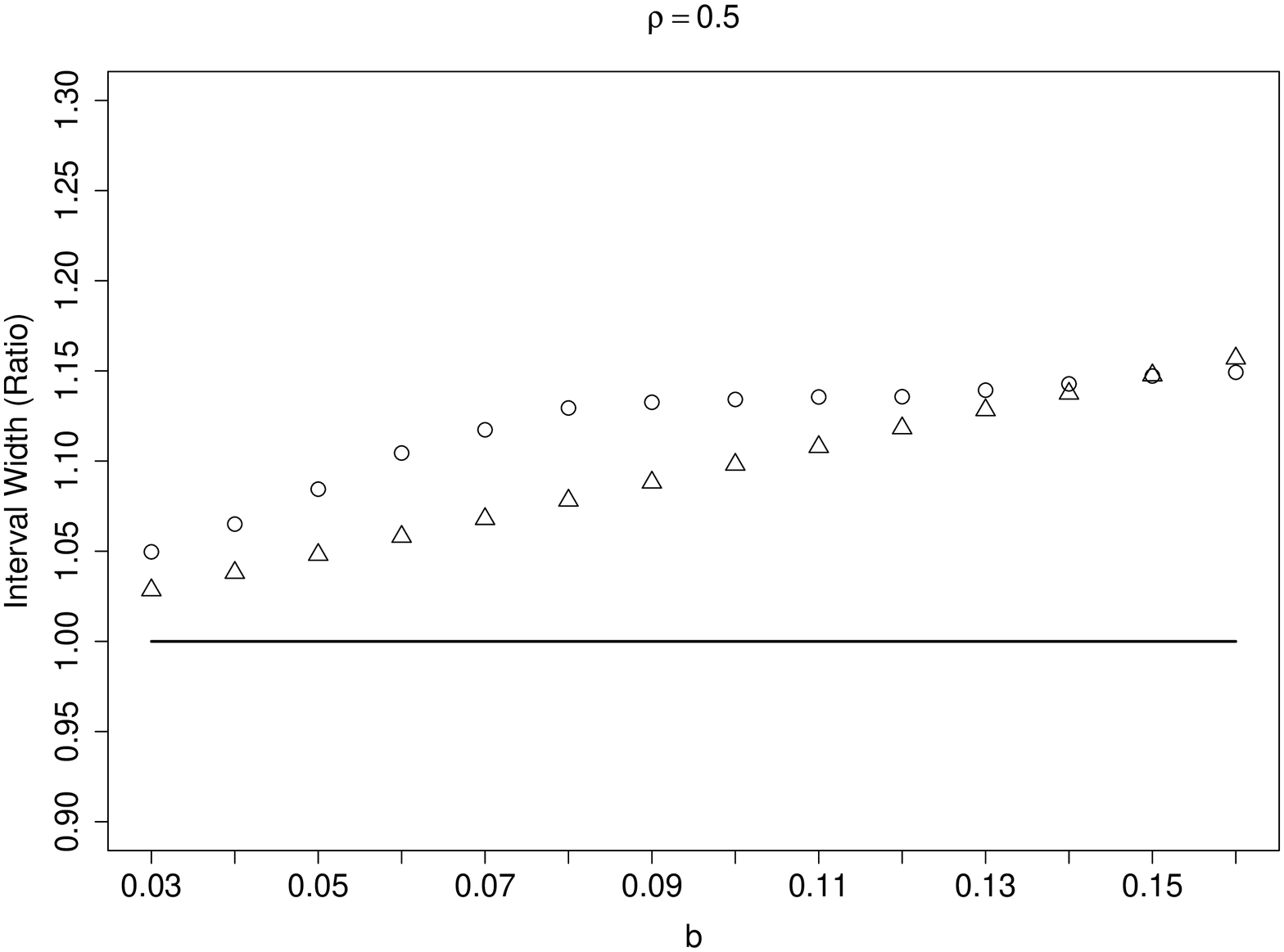}}
{\includegraphics[height=8cm,width=4.5cm,angle=270]{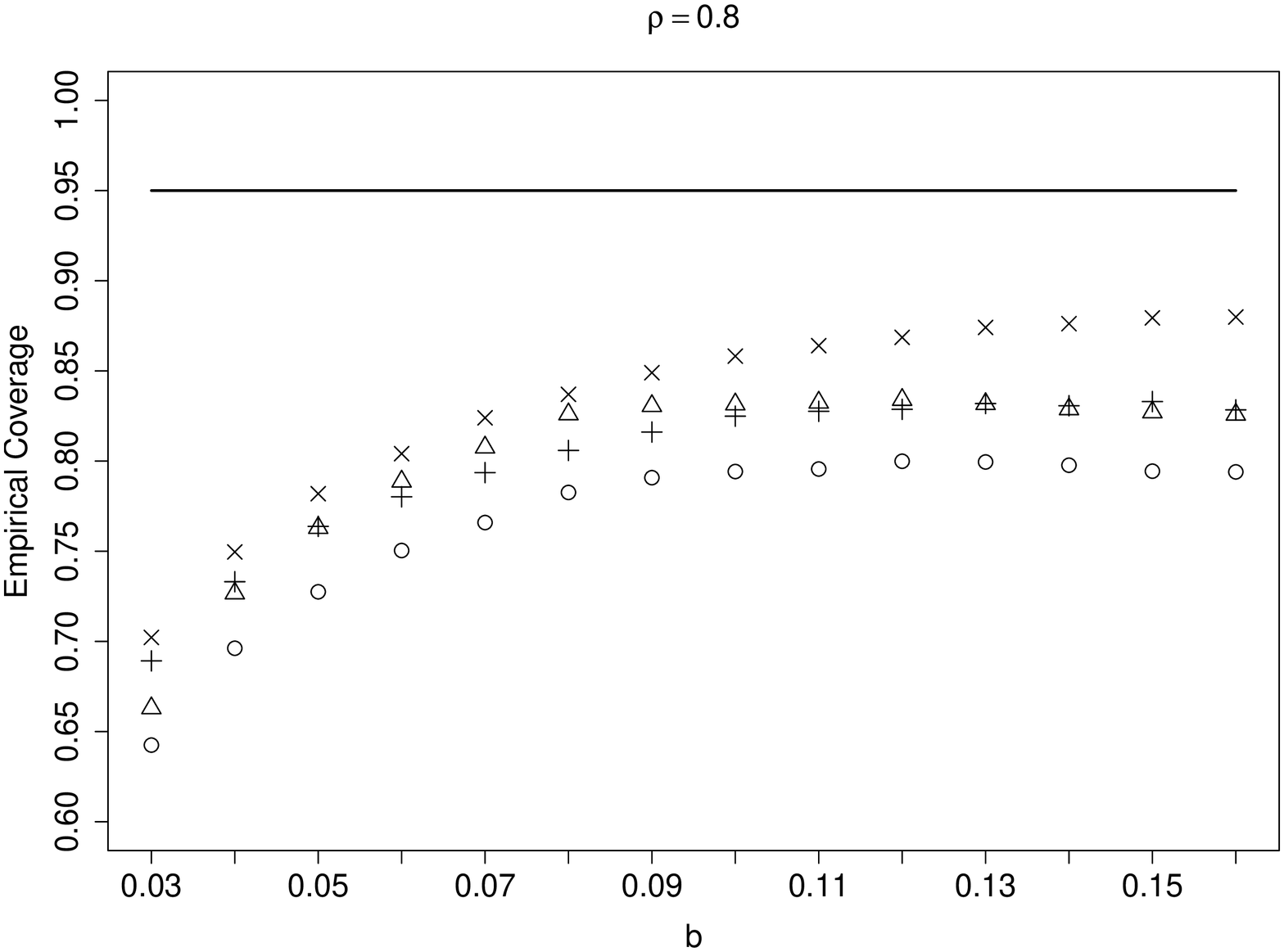}}
{\includegraphics[height=8cm,width=4.5cm,angle=270]{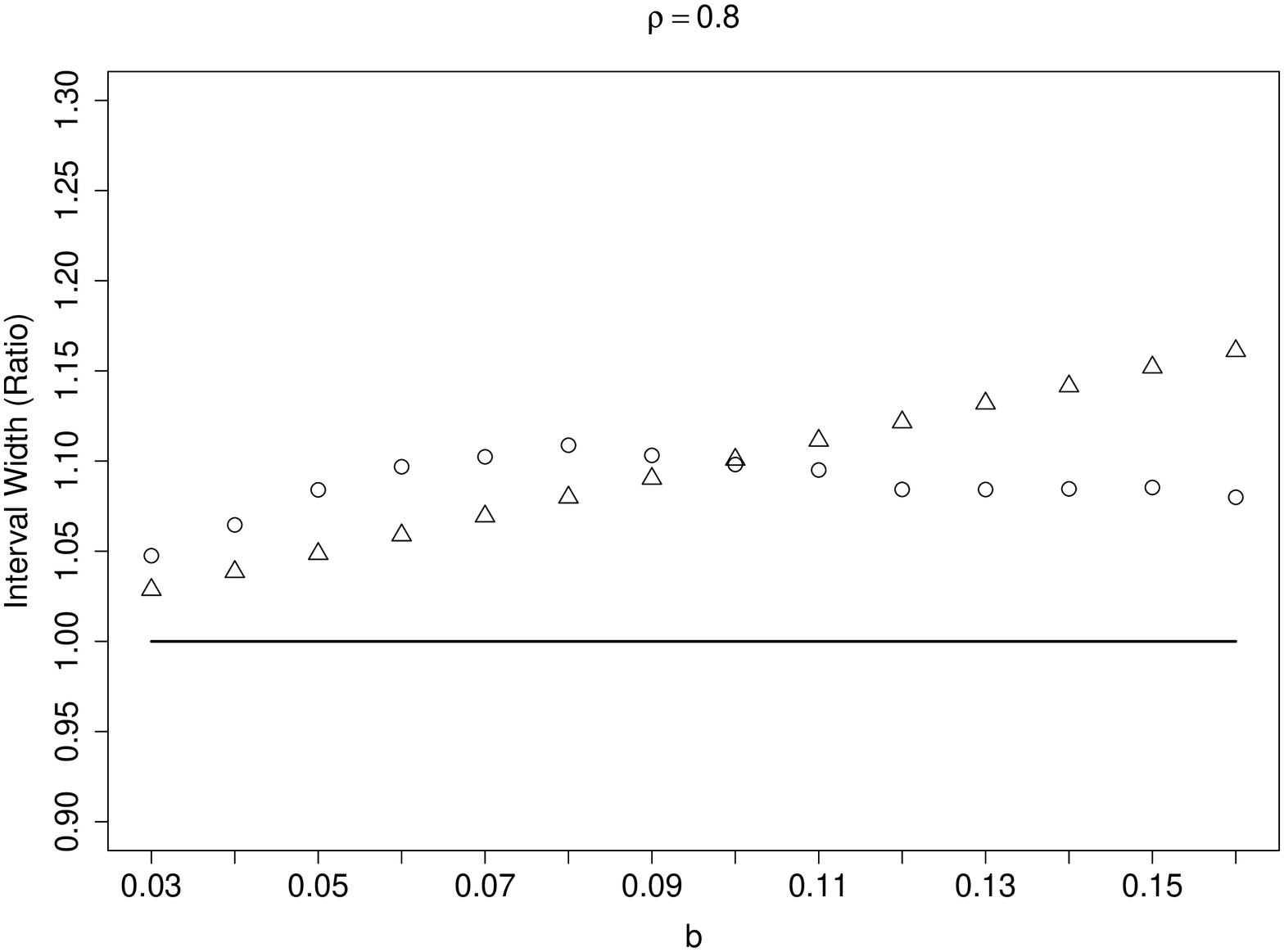}}
{\includegraphics[height=8cm,width=4.5cm,angle=270]{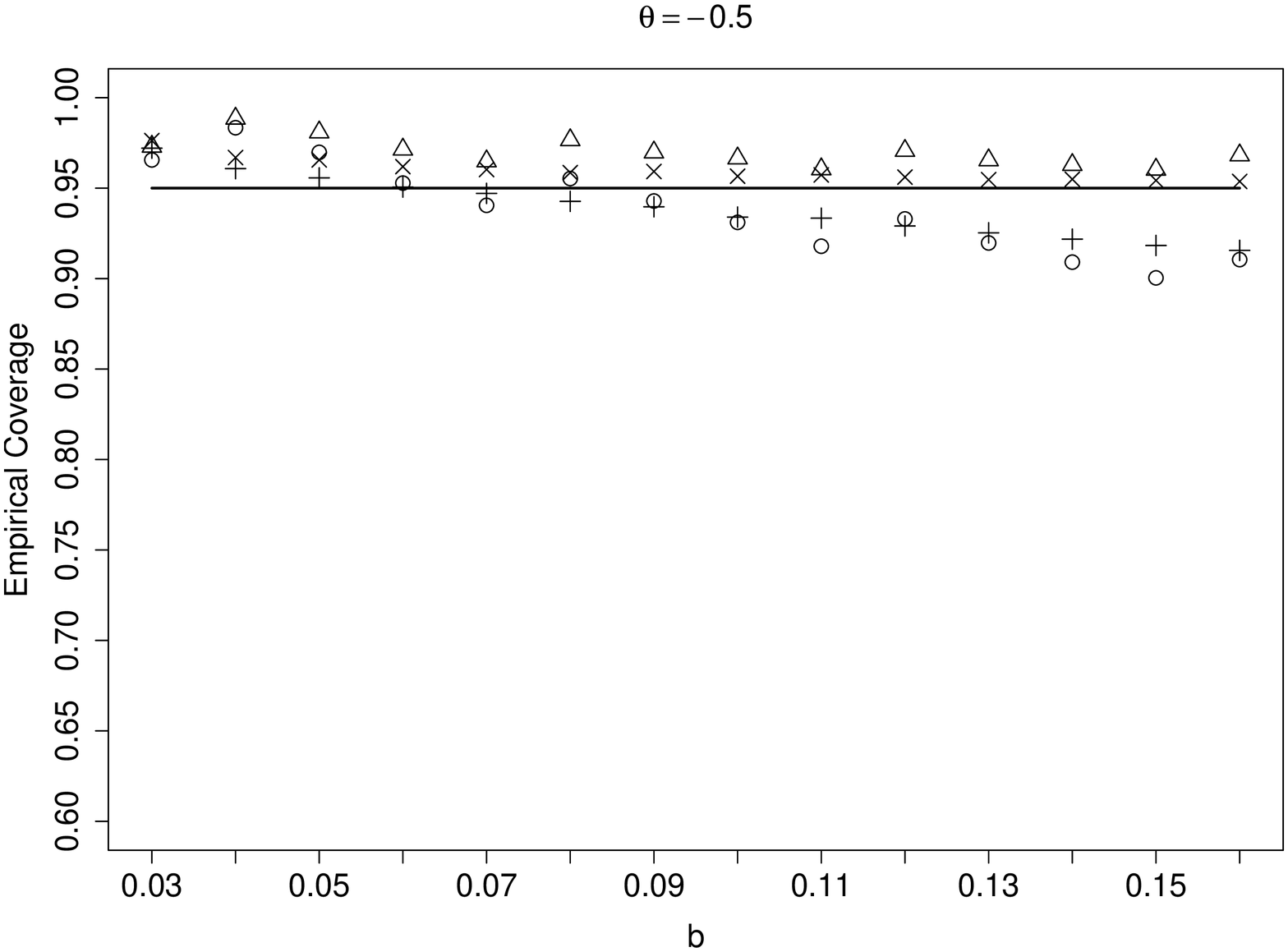}}
{\includegraphics[height=8cm,width=4.5cm,angle=270]{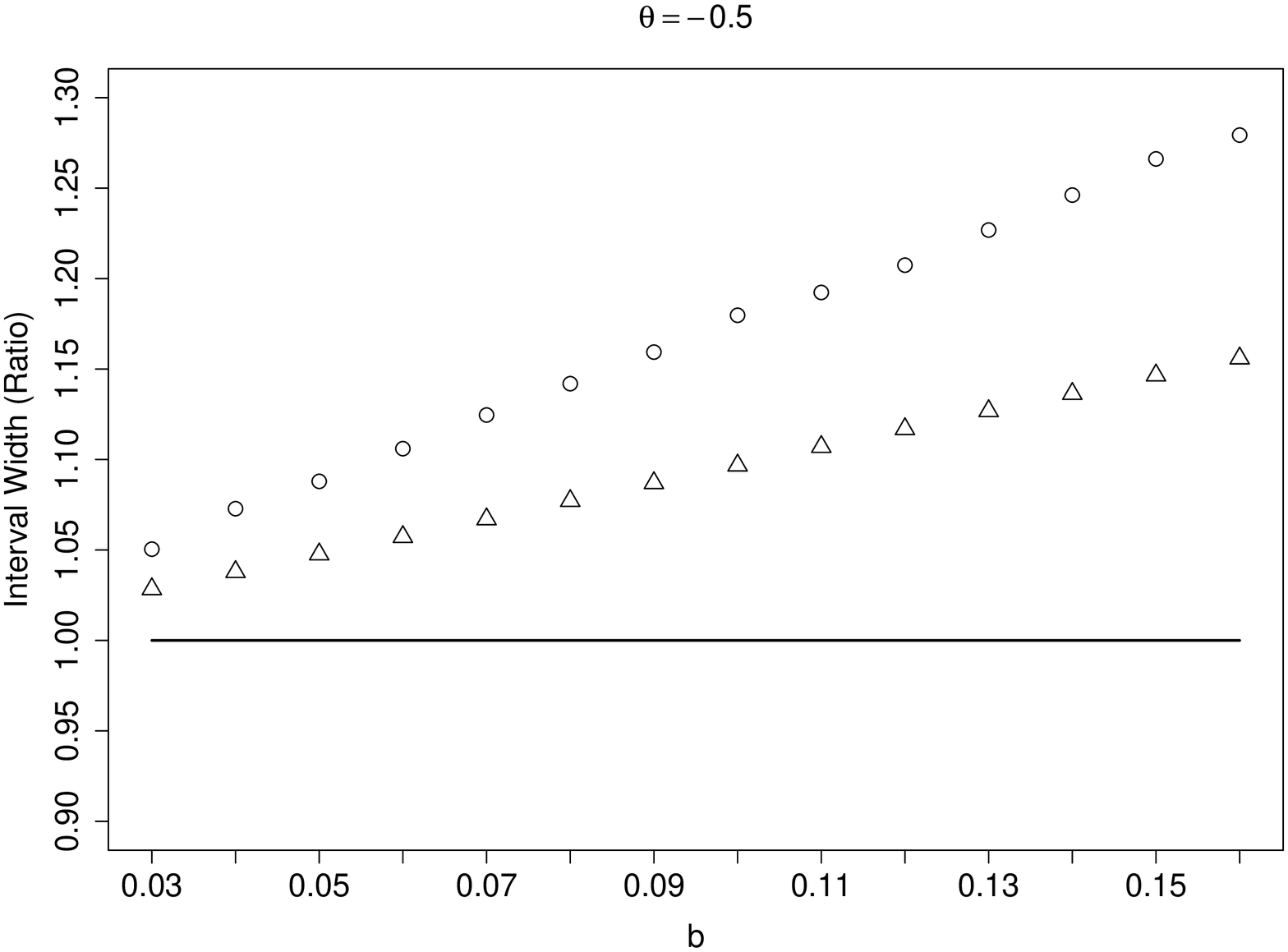}}
\label{fig:trmean1}
\end{center}

\end{figure}

\begin{figure}
\caption{The empirical coverage probabilities (left panel) and  the ratios of interval widths (calibrated fixed-$b$  over traditional small-$b$)
 (right panel) for the $25\%$ trimmed mean and for the models with exponentially distributed errors. Sample size $n=100$ and number of replications is 10000. }
\begin{center}
{\includegraphics[height=8cm,width=4.5cm,angle=270]{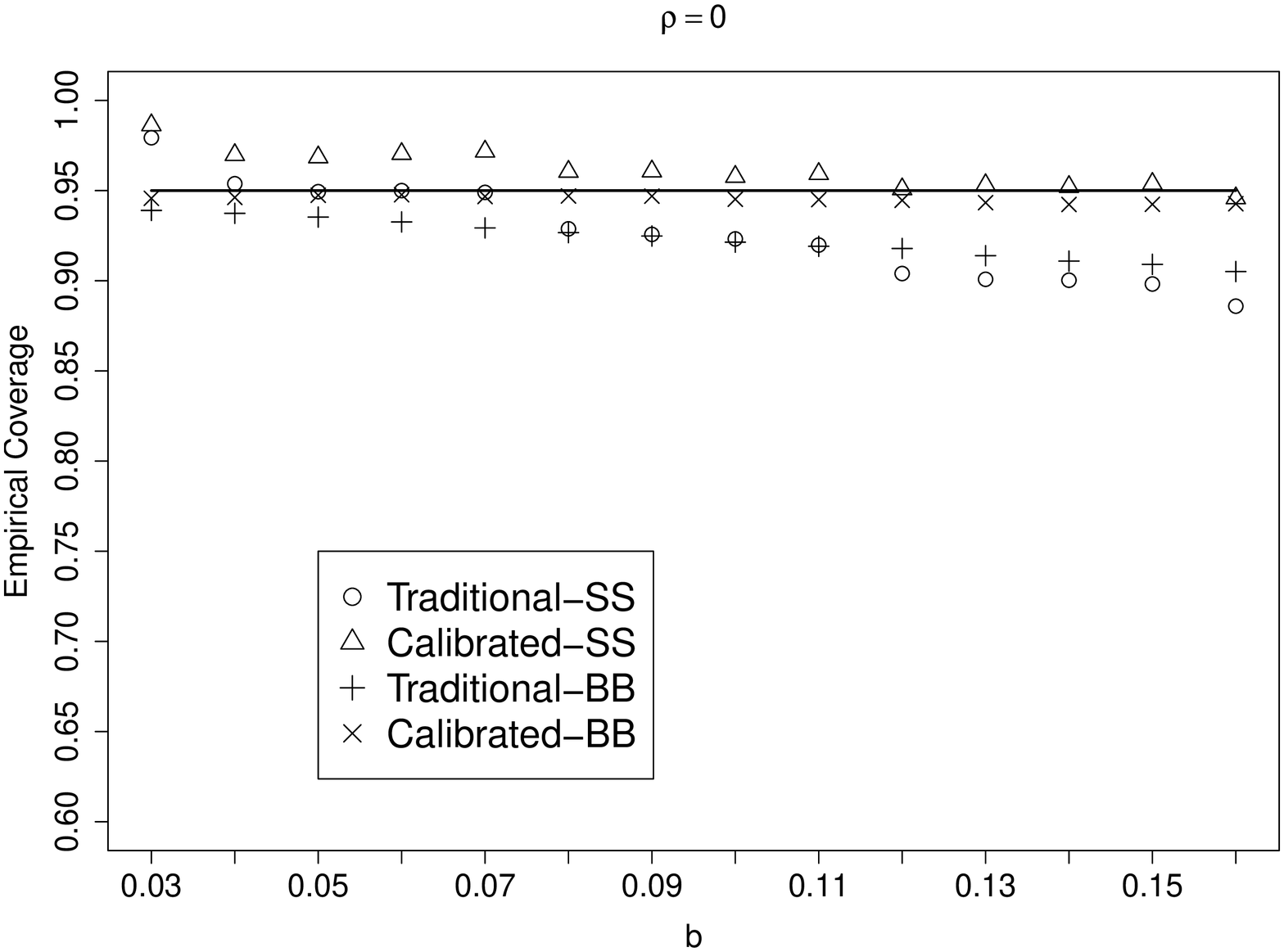}}
{\includegraphics[height=8cm,width=4.5cm,angle=270]{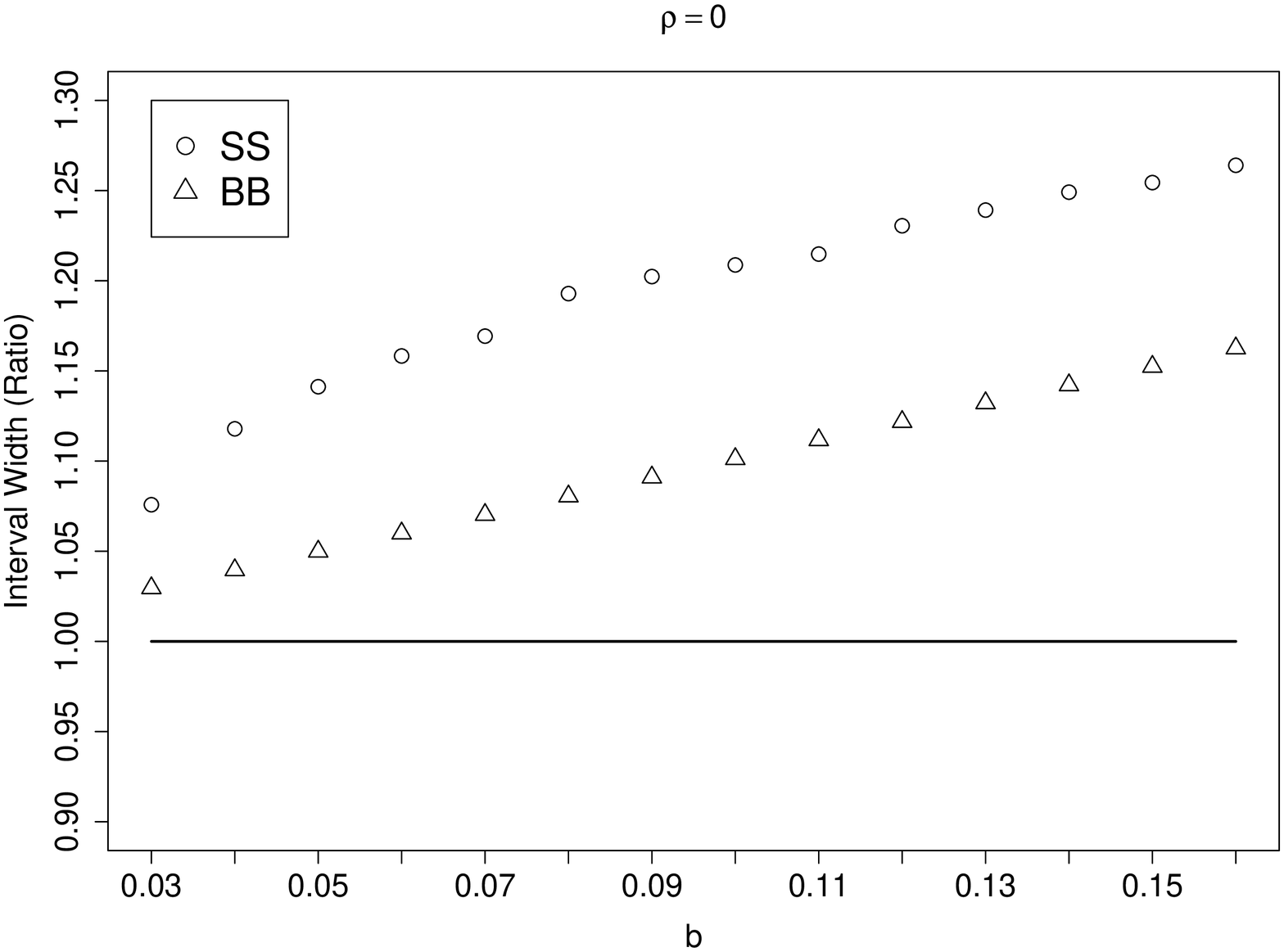}}
{\includegraphics[height=8cm,width=4.5cm,angle=270]{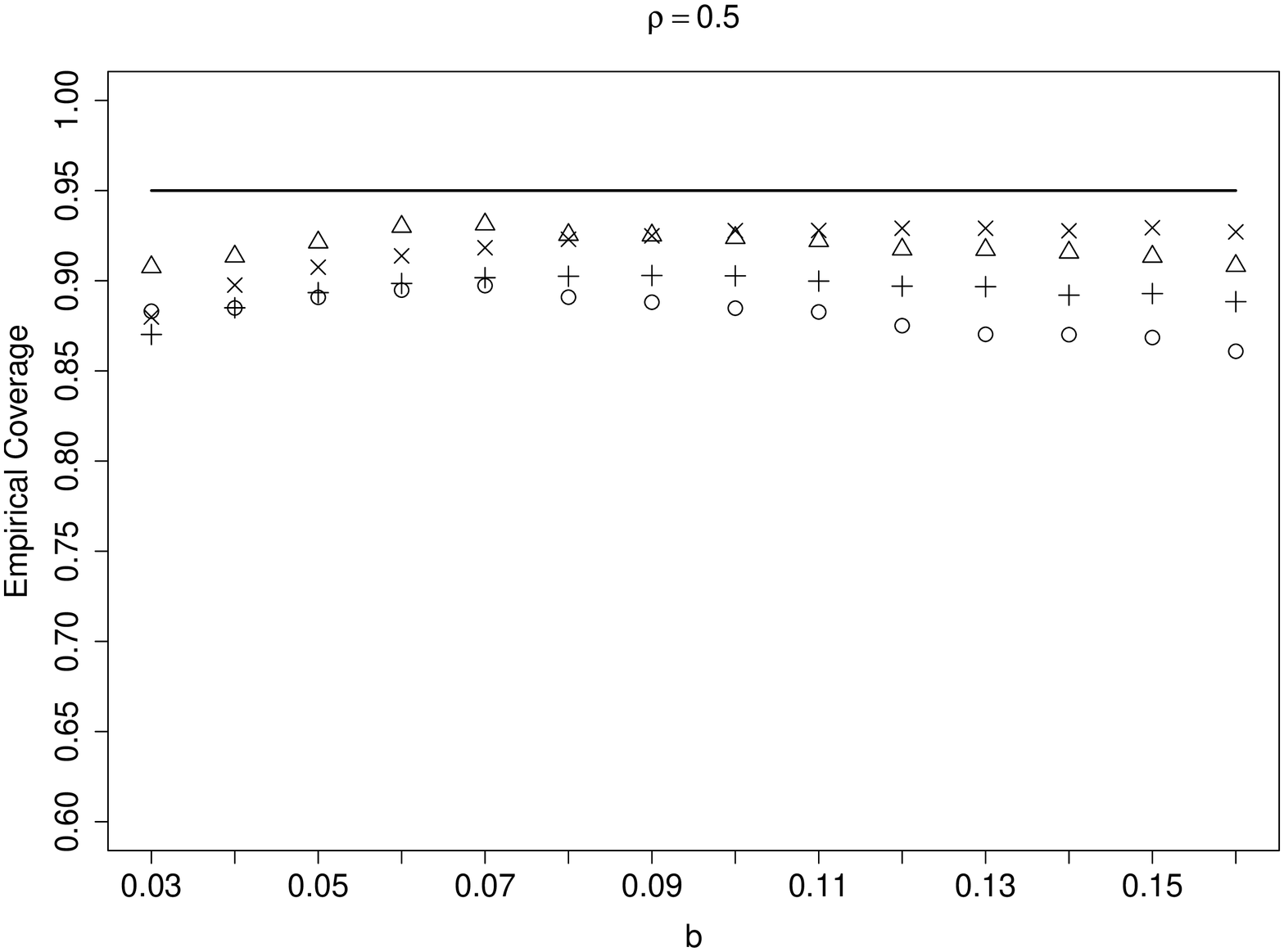}}
{\includegraphics[height=8cm,width=4.5cm,angle=270]{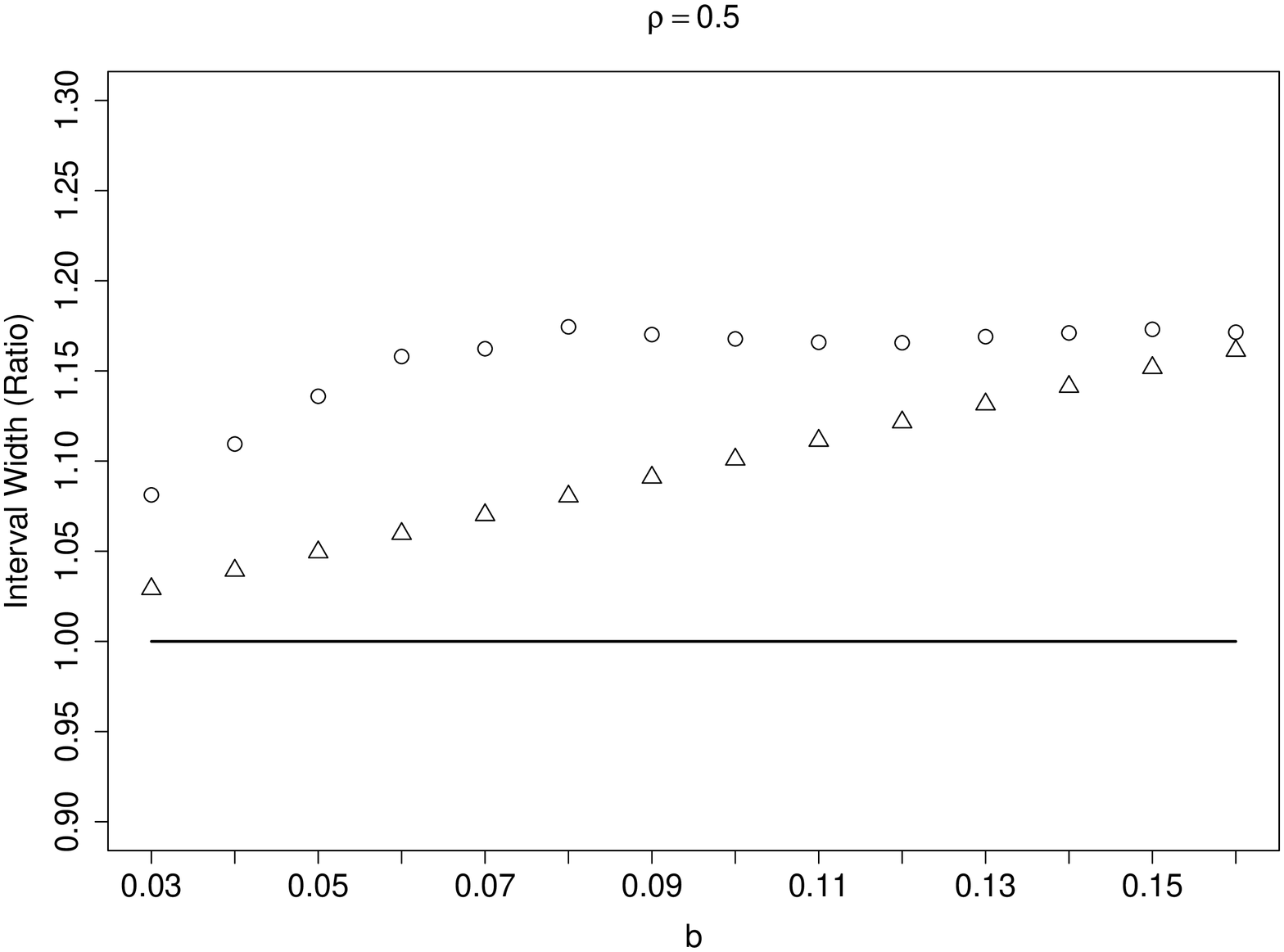}}
{\includegraphics[height=8cm,width=4.5cm,angle=270]{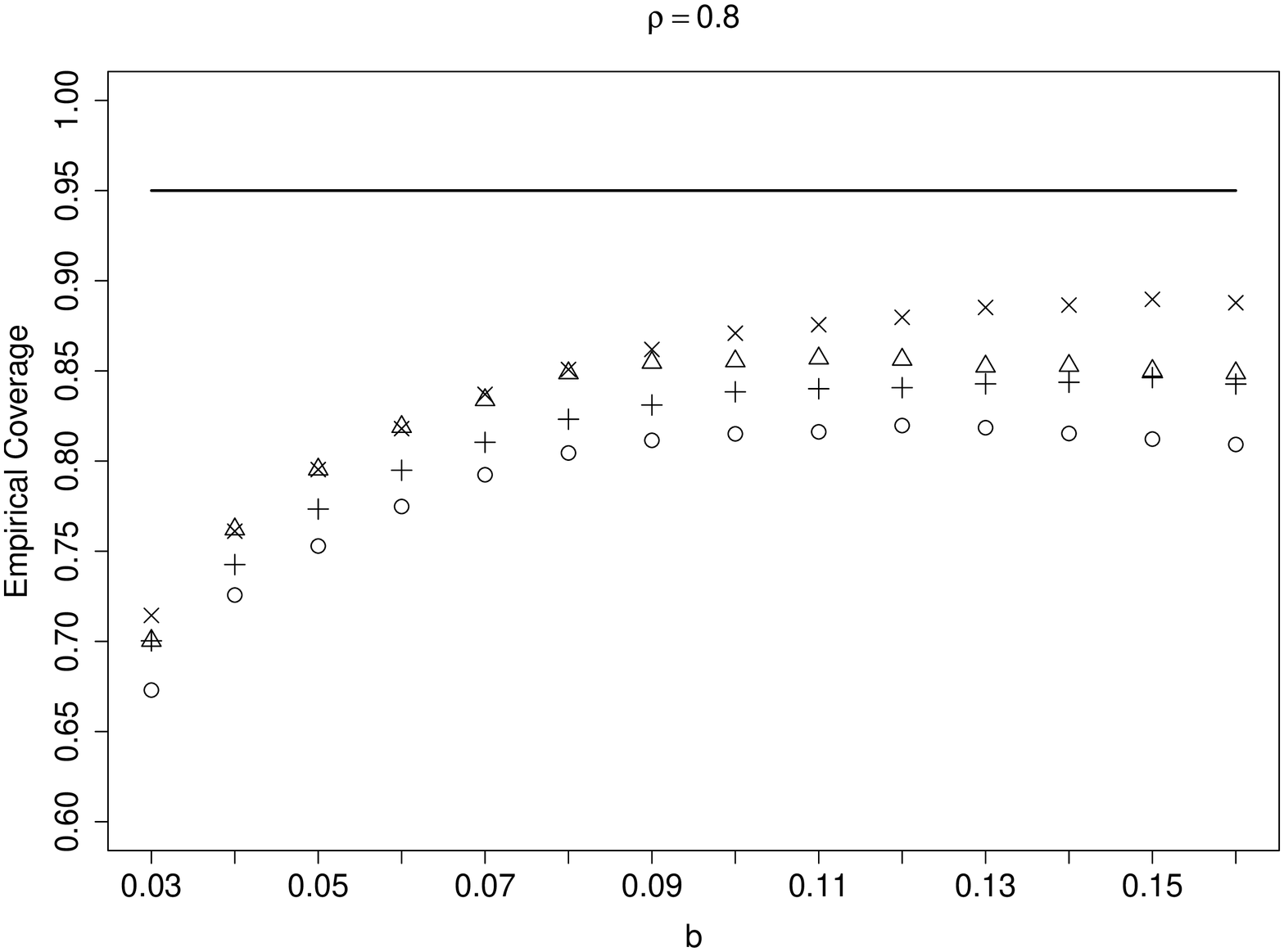}}
{\includegraphics[height=8cm,width=4.5cm,angle=270]{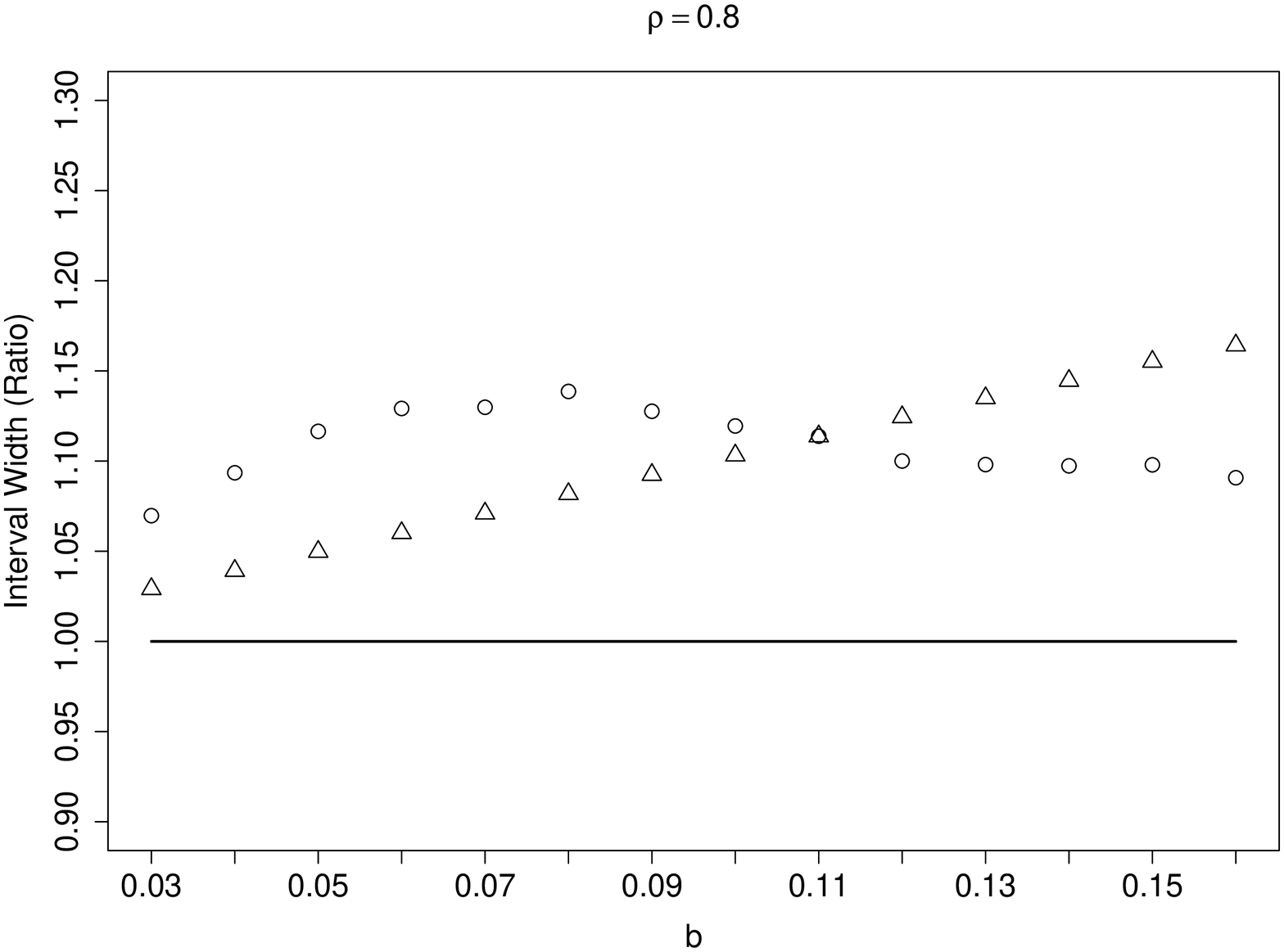}}
{\includegraphics[height=8cm,width=4.5cm,angle=270]{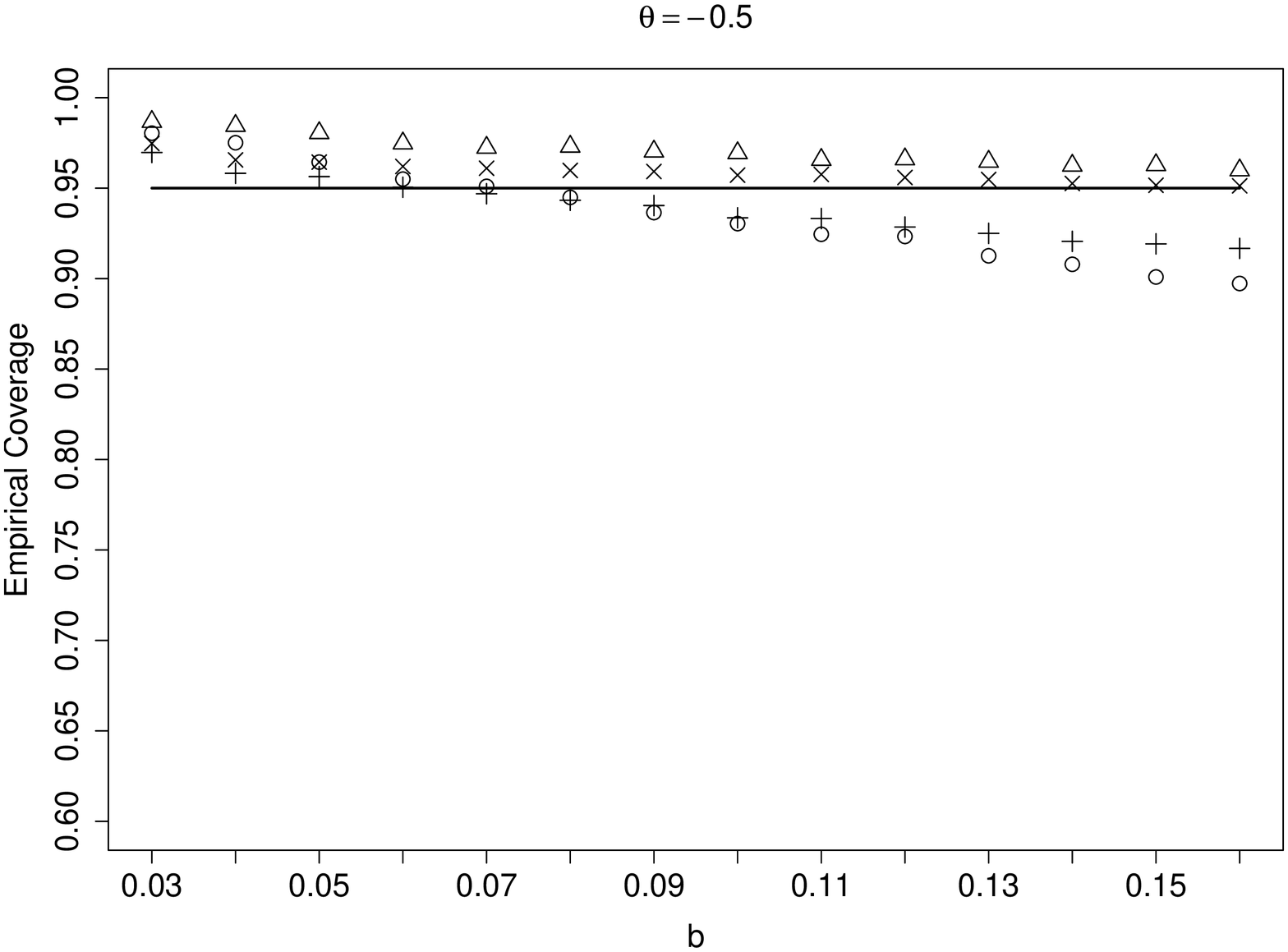}}
{\includegraphics[height=8cm,width=4.5cm,angle=270]{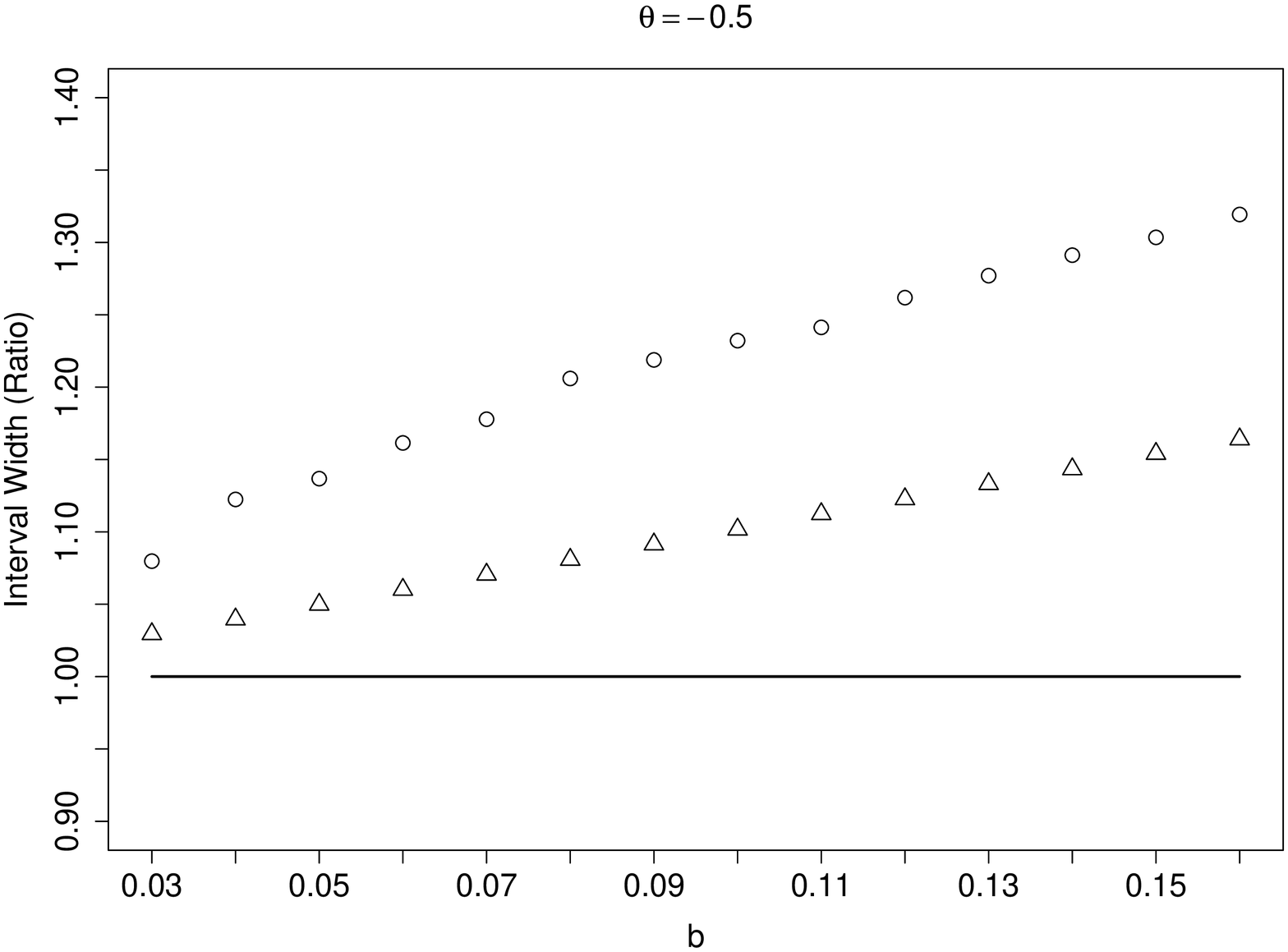}}
\label{fig:trmeanEXP1}
\end{center}

\end{figure}

\begin{figure}
\caption{The empirical coverage probabilities (left panel) and  the ratios of interval widths (calibrated fixed-$b$  over traditional small-$b$)
 (right panel) for the mean (top two plots) and the trimmed mean (bottom two plots) for the two nonlinear models. Sample size $n=100$ and number of replications is 10000. }
\begin{center}
{\includegraphics[height=8cm,width=4.5cm,angle=270]{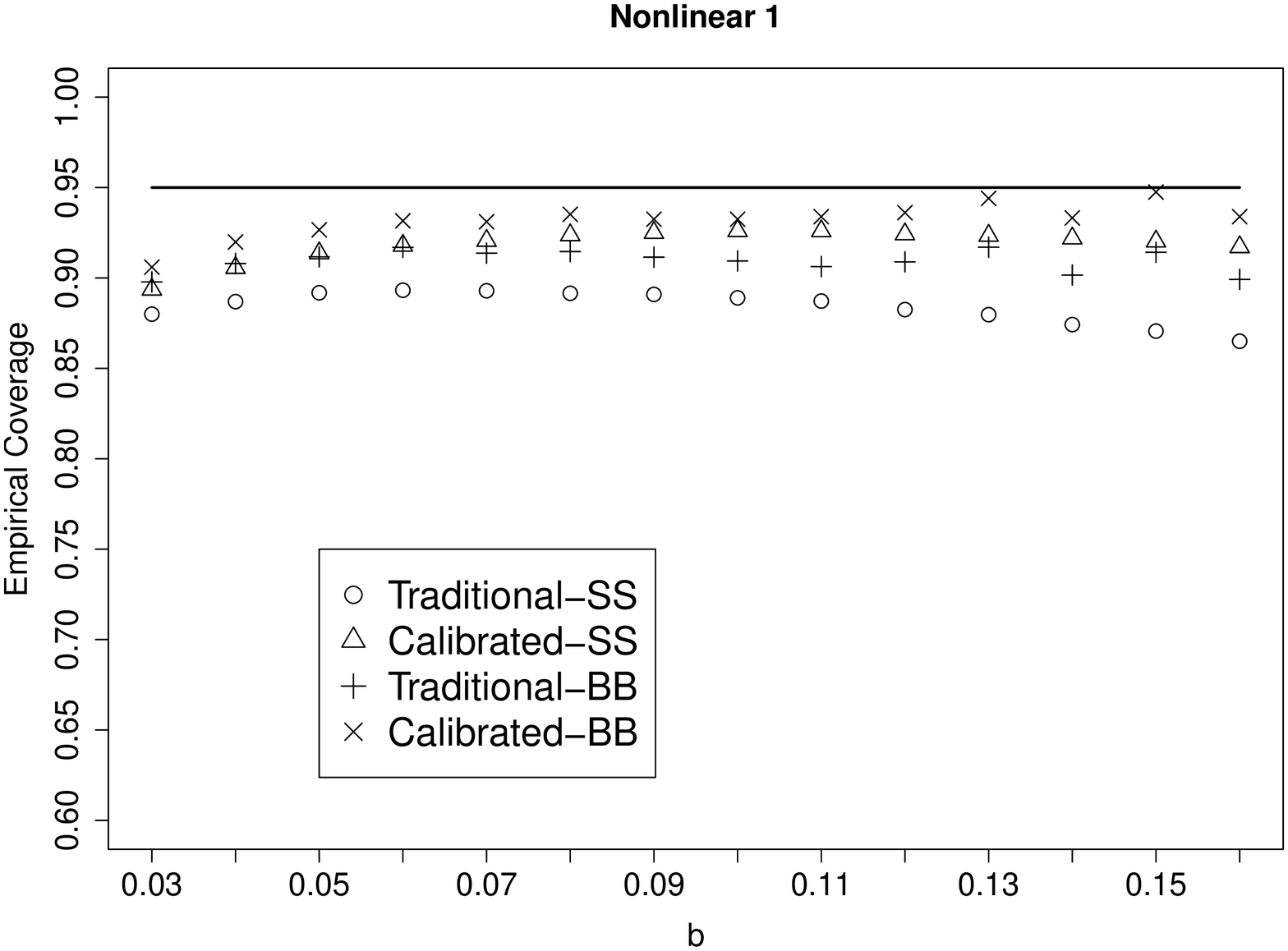}}
{\includegraphics[height=8cm,width=4.5cm,angle=270]{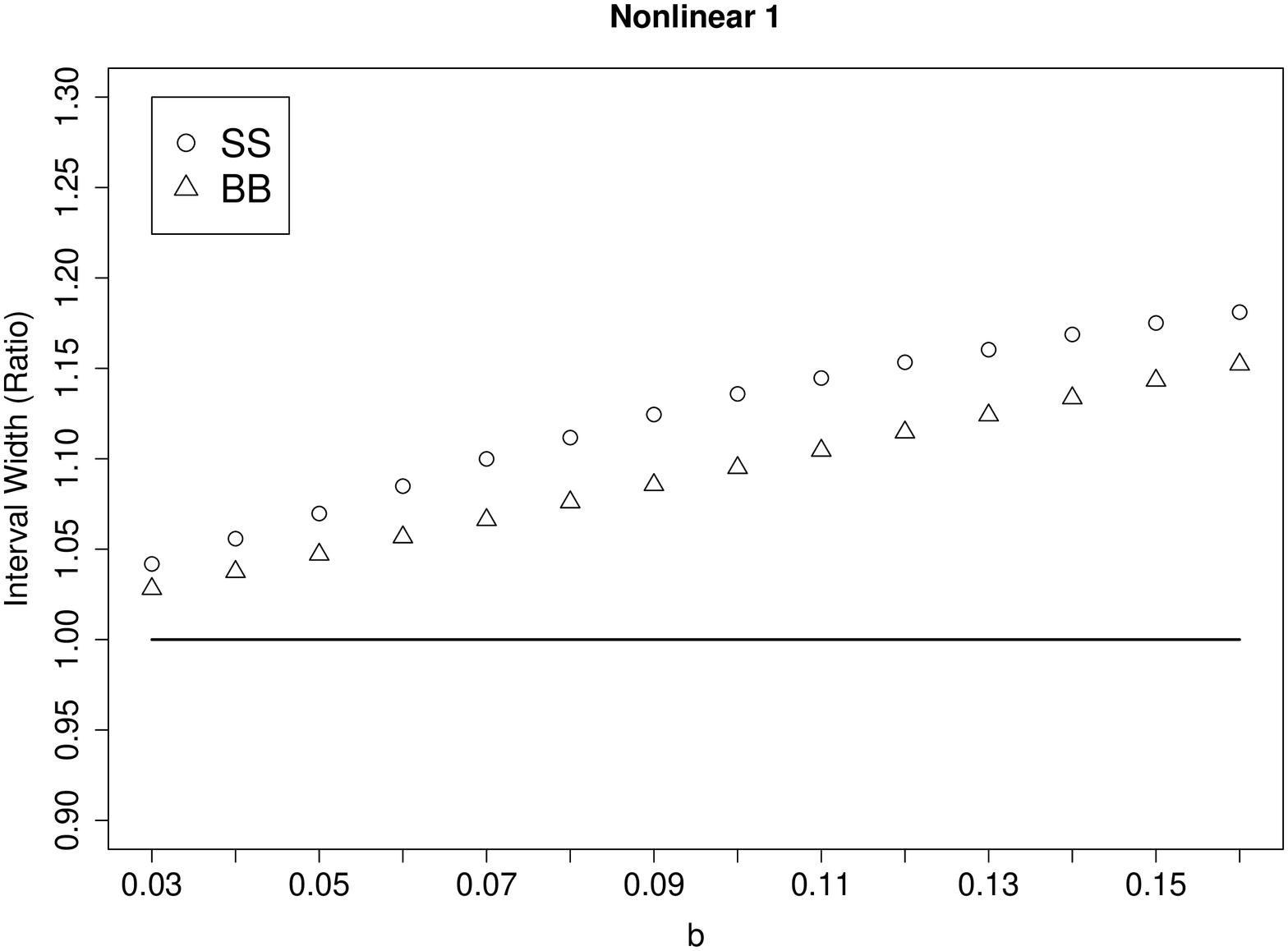}}
{\includegraphics[height=8cm,width=4.5cm,angle=270]{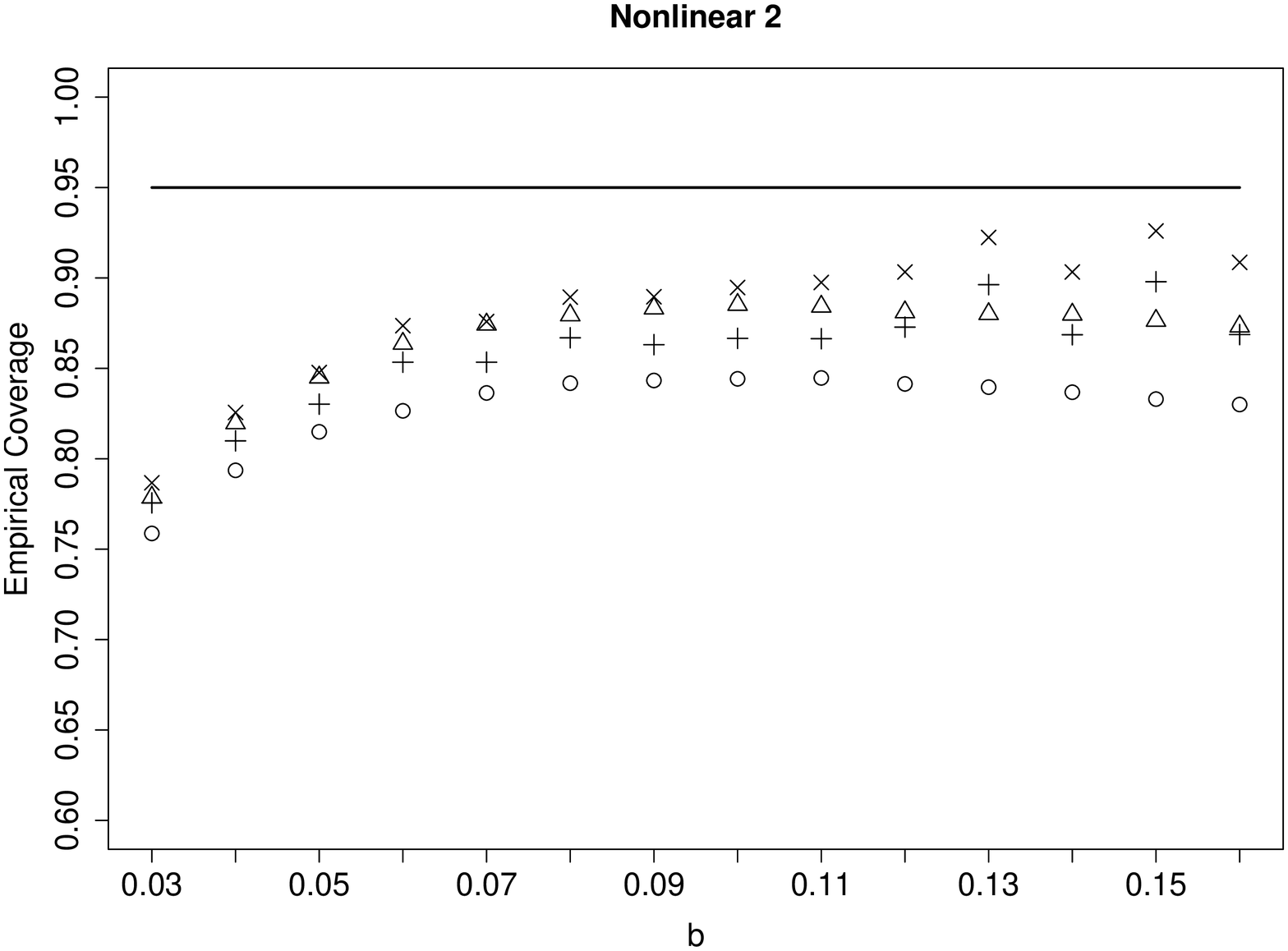}}
{\includegraphics[height=8cm,width=4.5cm,angle=270]{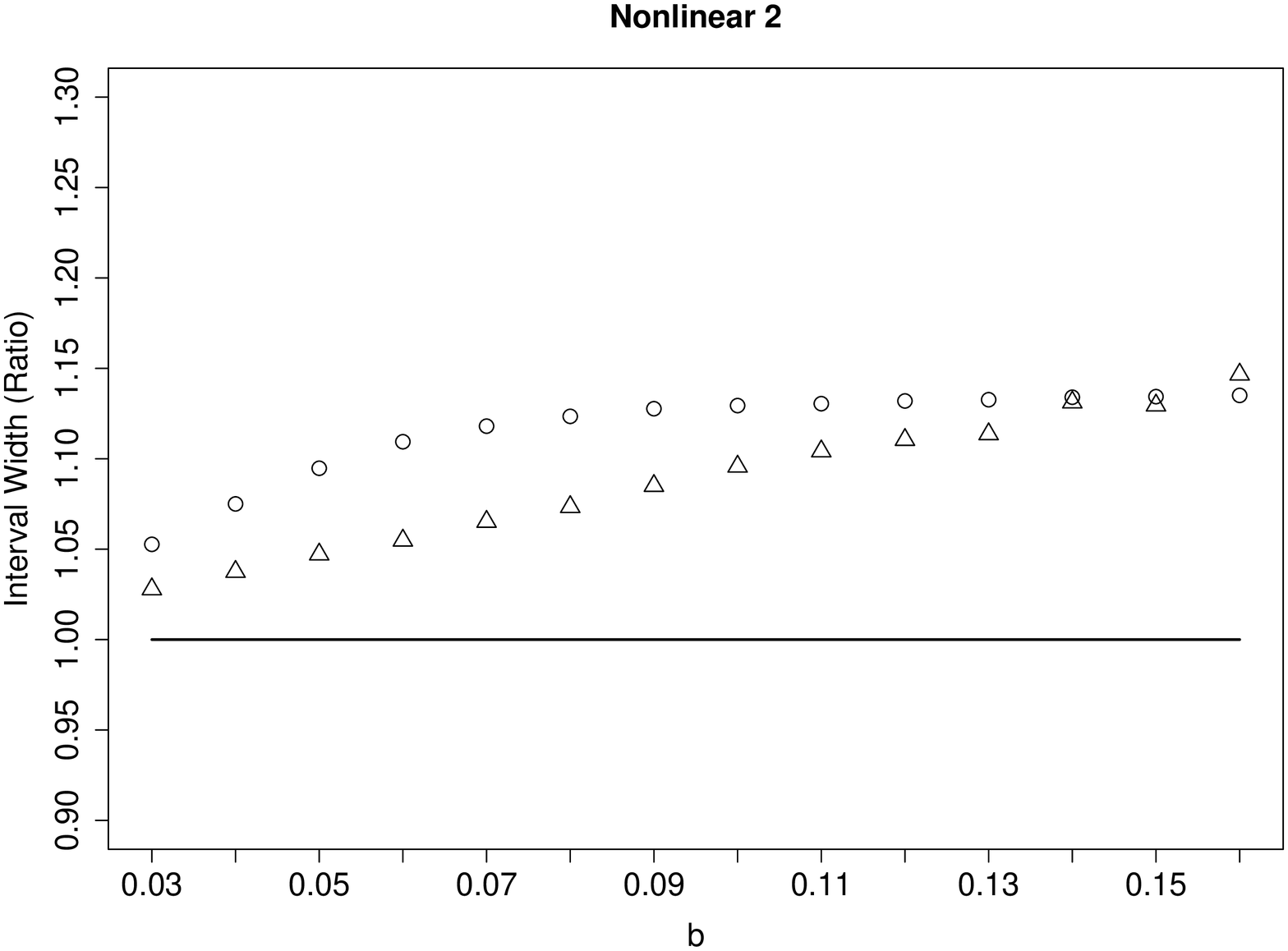}}
{\includegraphics[height=8cm,width=4.5cm,angle=270]{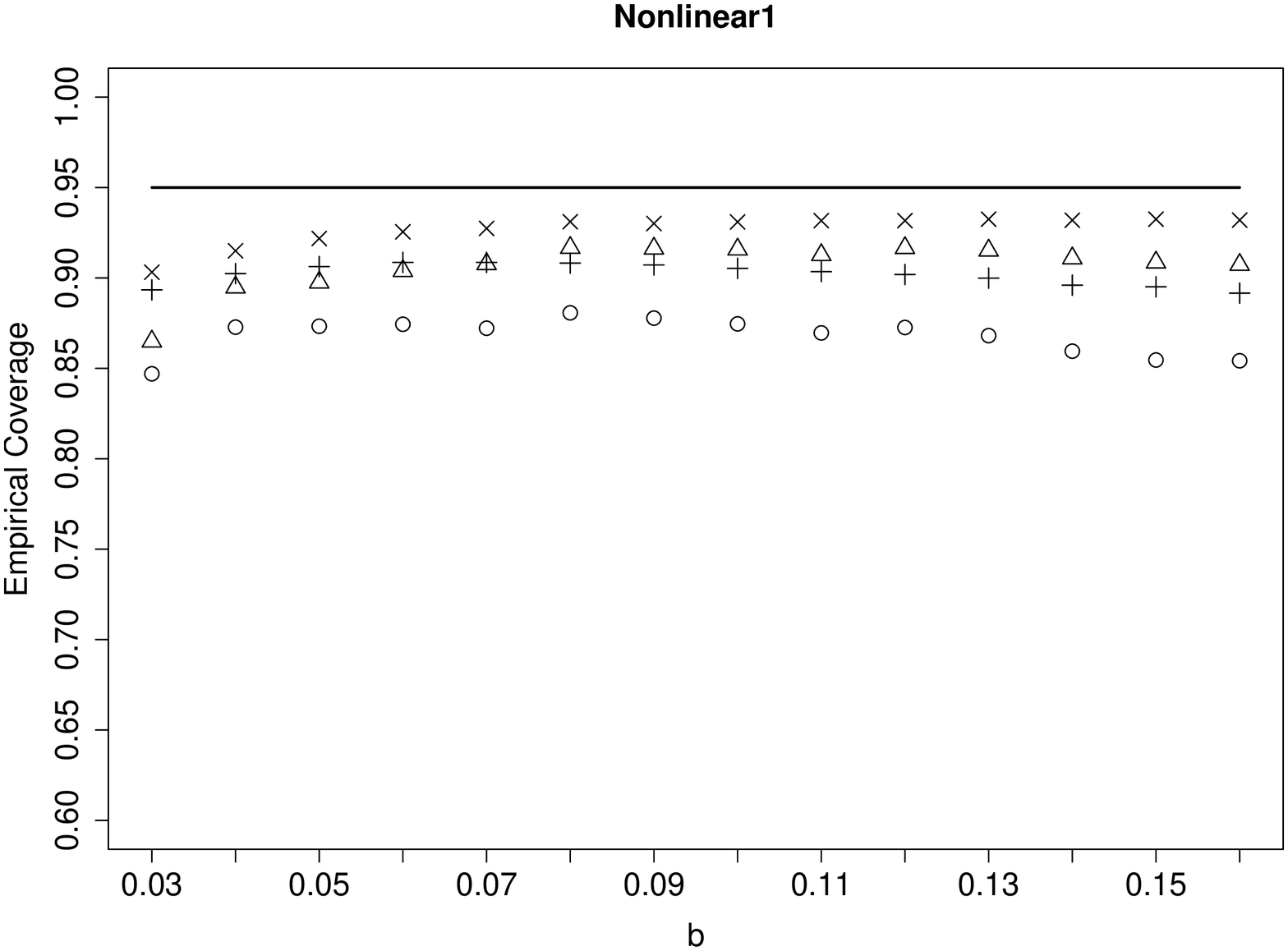}}
{\includegraphics[height=8cm,width=4.5cm,angle=270]{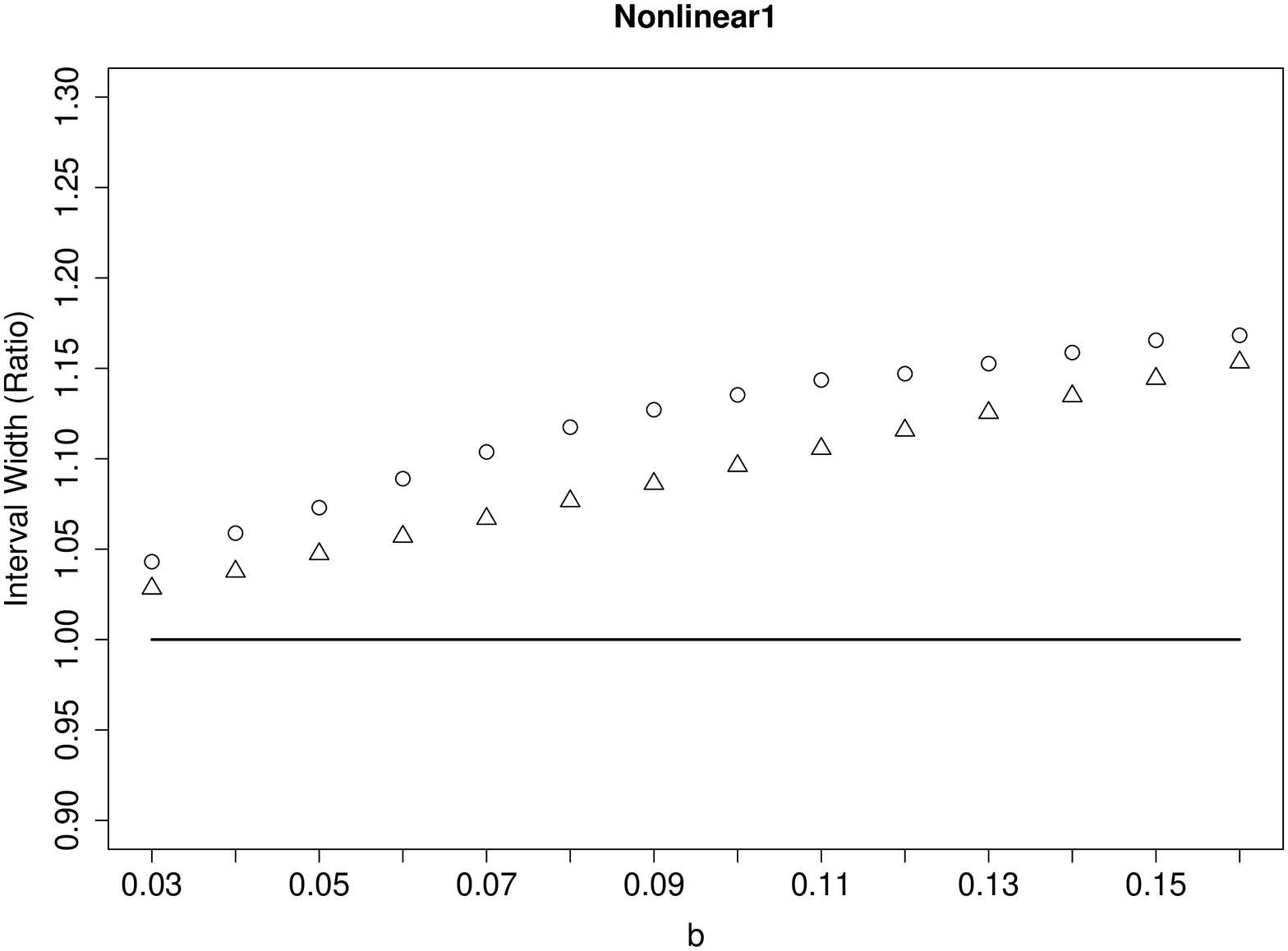}}
{\includegraphics[height=8cm,width=4.5cm,angle=270]{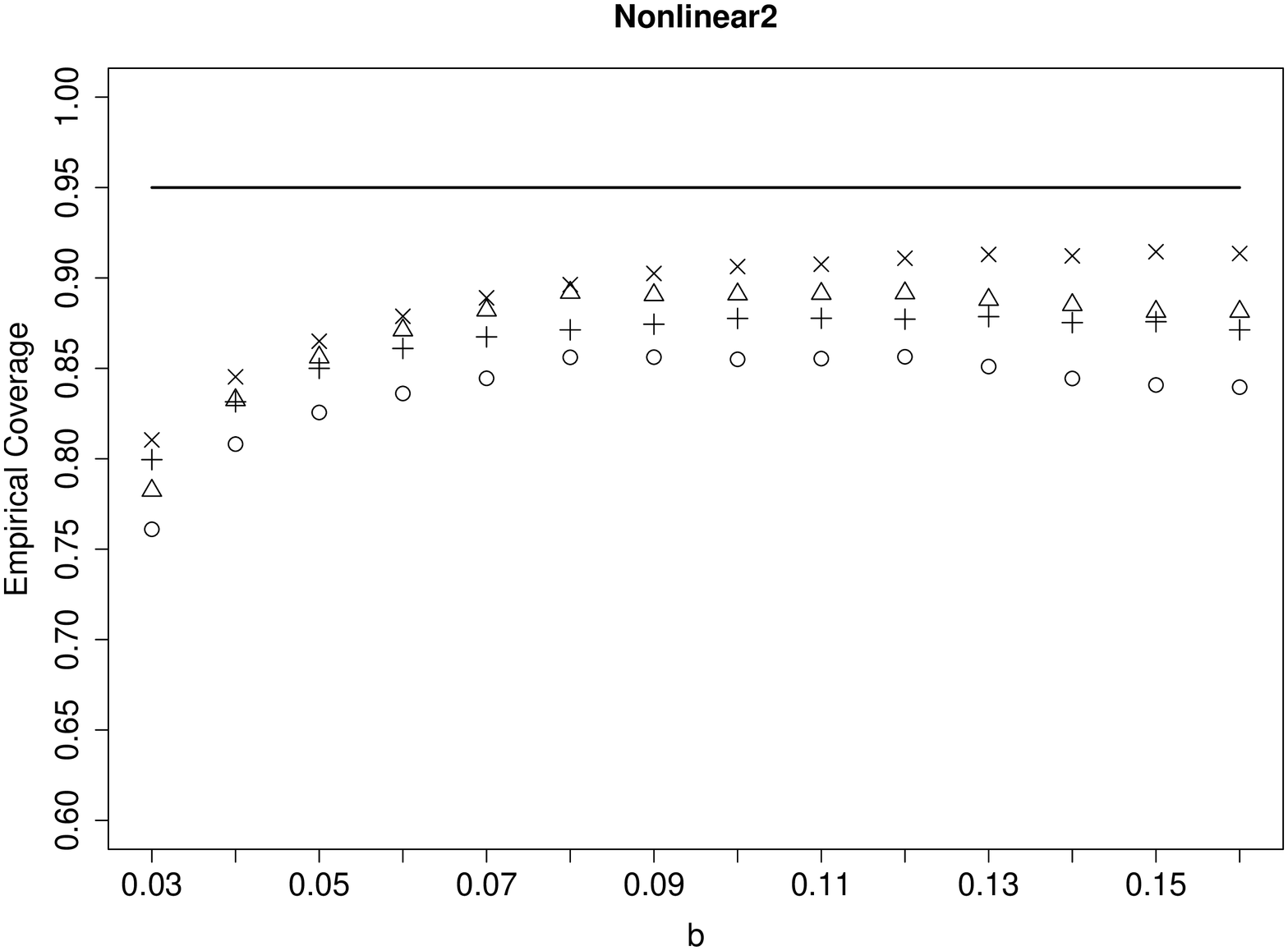}}
{\includegraphics[height=8cm,width=4.5cm,angle=270]{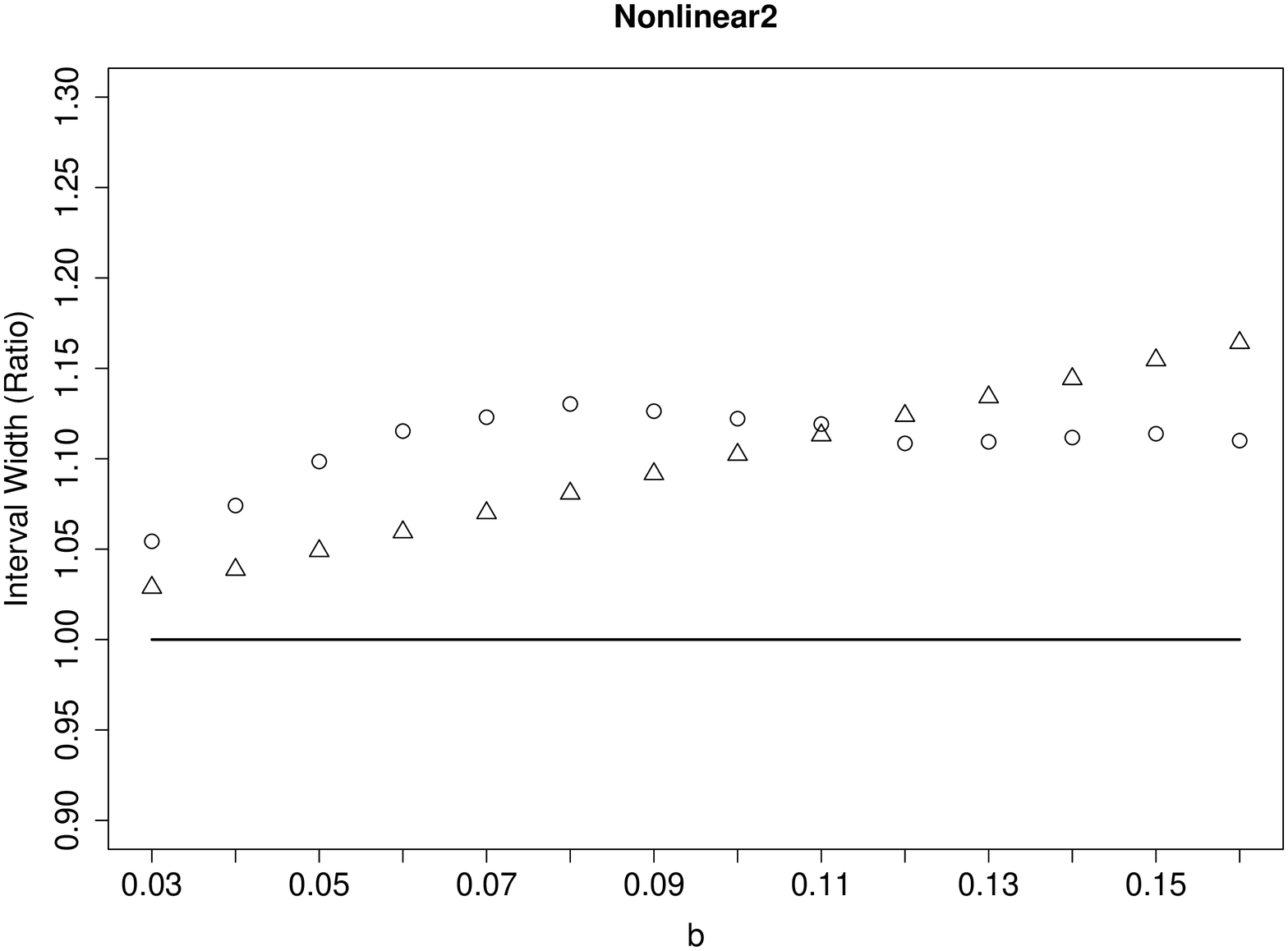}}
\label{fig:nonlinear}
\end{center}

\end{figure}

 The results for the models with exponentially distributed errors as presented in Figure~\ref{fig:meanEXP1} are very similar to the ones for normally distributed errors, suggesting  that the asymmetric shape of exponentially distributed errors has little impact on the finite sample performance. Figures~\ref{fig:trmean1} and~\ref{fig:trmeanEXP1} present the results for the trimmed mean case, which are fairly similar to the results in the mean case. Additionally, the results for the nonlinear models in Figure~\ref{fig:nonlinear} resemble those for AR(1)-N(0,1) models with $\rho=0.5$ in both the mean and the trimmed mean case, indicating that nonlinearity does not affect our results much.

 Due to the duality of confidence interval construction and hypothesis testing, we would expect that the fixed-b approach leads to better size (i.e. size closer to the nominal one) in all the models except MA(1) with $\theta=-0.5$, at the sacrifice of (raw) power. The power loss is expected to be moderate because the ratio of fixed-$b$ based interval width over the small-$b$ based interval width is quite close to $1$. Overall, the simulation results demonstrate that the fixed-$b$ approach delivers more accurate inference for both subsampling and the moving block bootstrap in most situations owing to its more accurate approximation to the finite sample distribution. Of course, we only show the improved accuracy of the fixed-$b$ approximation for a specific $\alpha=0.05$, which is also what Kiefer and Vogelsang did. We also tried $\alpha=10\%$ and qualitatively similar results are obtained.
 It would be interesting to provide some theoretical justifications on the order of the error rejection probability. For the subsampling method,  this boils down to   the order of
 $\sup_{\alpha\in [0,1]} |P(\widetilde{pval}_{n,l}^{SUB}\le \alpha)-P(\widetilde{G}(l/n)\le \alpha)|$ under the fixed-$b$ framework. Note that under the small-$b$ framework, the
error is  $\sup_{\alpha\in [0,1]} |P(\widetilde{pval}_{n,l}^{SUB}\le \alpha)-\alpha|,$
  which is expected to be larger. A formal theoretical proof seems quite challenging and is left for future research.

\subsection{Finite sample performance of confidence regions and confidence bands}
\label{sec:regionband}

In this subsection, we examine the coverage probabilities of confidence regions for the vector parameter of  mean and median, and
 confidence bands for the marginal distribution function $m(\cdot)$ and the
normalized spectral distribution function $\widetilde{F}(\cdot)$. Let $\{X_t\}_{t=1}^{n}$ be generated
from the AR$(1)$ model: $X_t=\rho X_{t-1} + e_t$, where $\rho=-0.6, 0, 0.5, 0.8$, $e_t\sim iid~ N(0,1)$ or $\mbox{EXP}(1)-1$. 
Sample size $n=200$ and
 number of replications is 1000. We use the Euclidean norm in the confidence region construction.
 For both confidence regions and confidence bands, we compared the following three schemes: (1) traditional subsampling-based confidence
 region (band); (2) calibrated subsampling-based
 confidence region (band) with a fixed $n'$, where $n'=15$ for confidence region construction and $n'=30$ for confidence band construction; (3) calibrated subsampling-based confidence
 region (band) with $n'$ chosen in a data driven fashion.  Here we employ a variant of a  block size selection procedure proposed in Bickel and Sakov (2008) for the $m$ out of $n$ bootstrap (also see G\"otze and Ra$\breve{c}$kauskas (2001)), which is closely related to the subsampling method. The use of Bickel and Sakov's automatic bandwidth selection in the subsampling context has been explored in Jach et al. (2011) recently.
 The procedure consists of the following steps (in the confidence band case):
\begin{enumerate}
\item[Step 1] For a predetermined interval $[K_1,K_2]$ and $g\in (0,1)$, we consider a sequence of $n_j$'s of the form
$n_j=\lfloor g^{j-1} K_2\rfloor, ~\mbox{for}~ j=1,2,\cdots, \lfloor\log(K_2/K_1)/\{-\log(g)\}\rfloor$.
\item[Step 2] For each $n_j$, find $J_{n,n_j}$, where $J_{n,n_j}$ is the subsampling-based distribution estimator for the sampling distribution of the p-value.
\item[Step 3] Set $j_0 = \mbox{argmin}_{j=1,\cdots, \lfloor\log(K_2/K_1)/\{-\log(g)\}\rfloor }~\sup_{x\in\R} |J_{n,n_j}(x) - J_{n,n_{j+1}}(x)|$.
 Then the optimal block size is $g^{j_0}K_2$. If the difference is minimized for a few values of $j$, then pick the largest among them.
\end{enumerate}

In our simulation experiment, we set $(K_1,K_2,g)=(5,40,0.75)$ for confidence region construction and  $(K_1,K_2,g)=(10,60,0.75)$
for confidence band construction, which corresponds to a sequence of block lengths as $(40, 30, 22, 16, 12,  9,  7,  5)$
and $(60, 45, 33, 25,18, 14, 10)$, respectively.

\begin{figure}
\caption{The empirical coverage probabilities (left panel) and  the ratios of radii of confidence regions (calibrated fixed-$b$  over traditional small-$b$)
 (right panel) for the vector parameter and for the models with normally distributed errors. Sample size $n=200$ and number of replications is 1000. }
\begin{center}
{\includegraphics[height=8cm,width=4.5cm,angle=270]{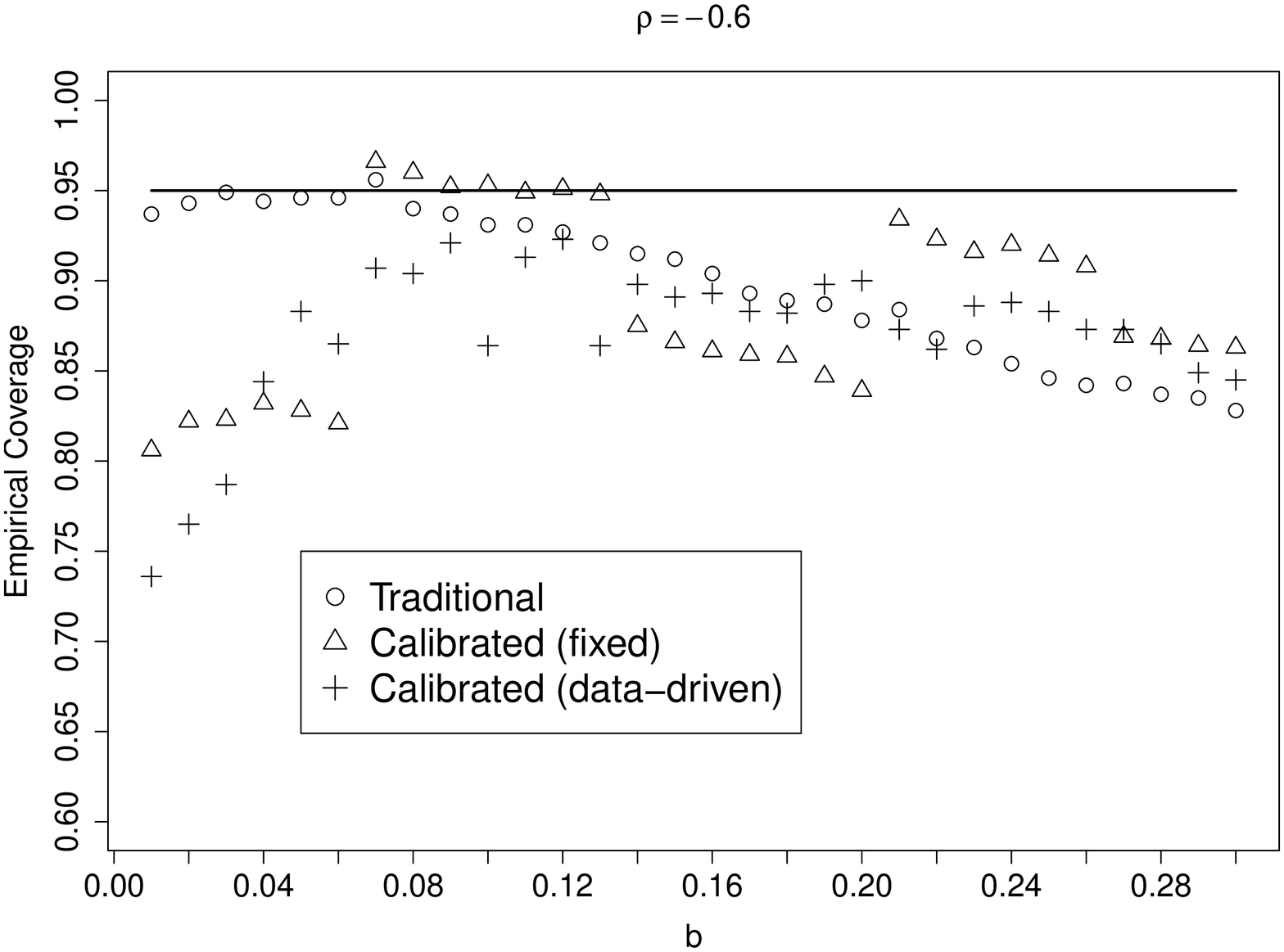}}
{\includegraphics[height=8cm,width=4.5cm,angle=270]{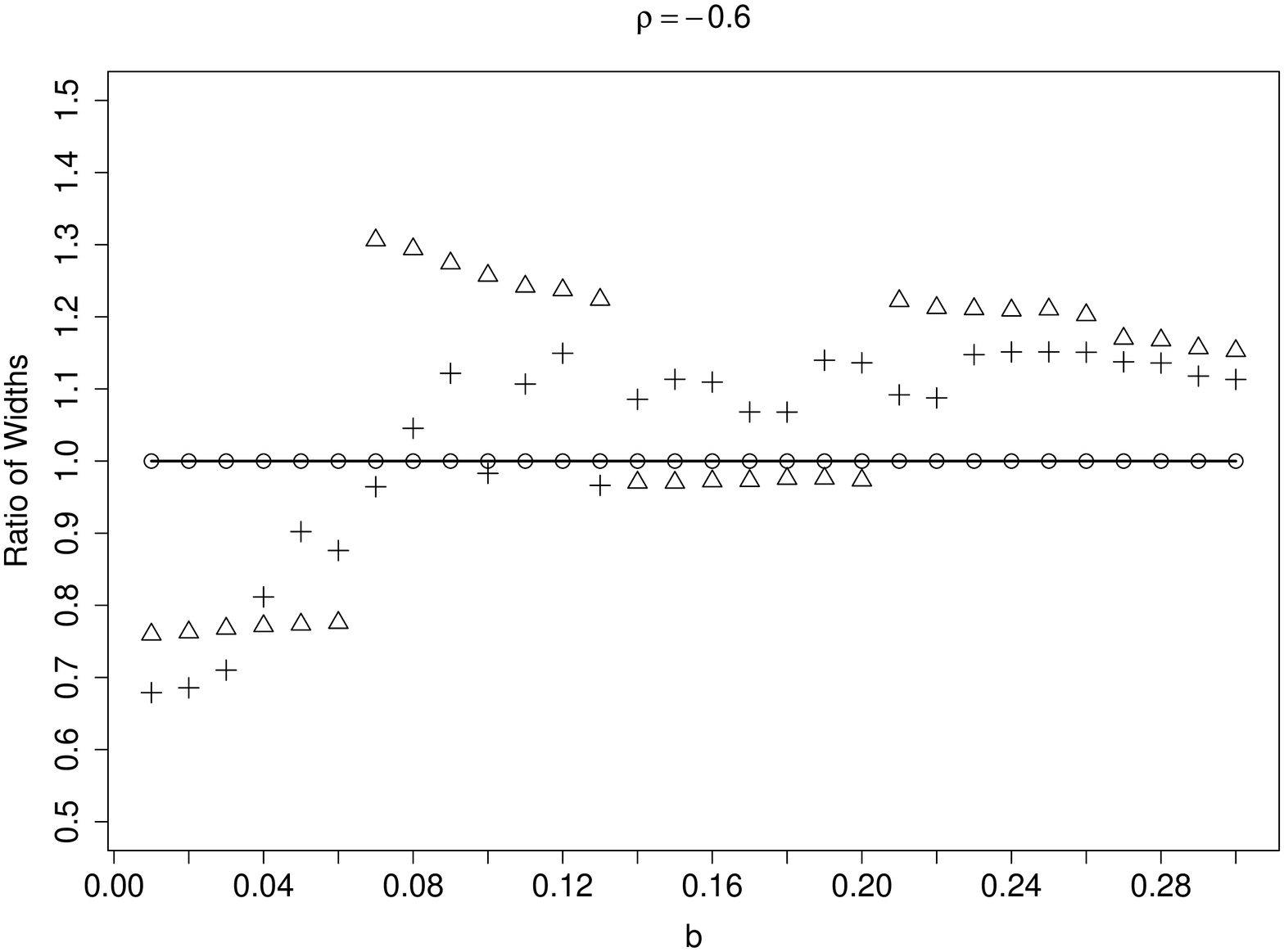}}
{\includegraphics[height=8cm,width=4.5cm,angle=270]{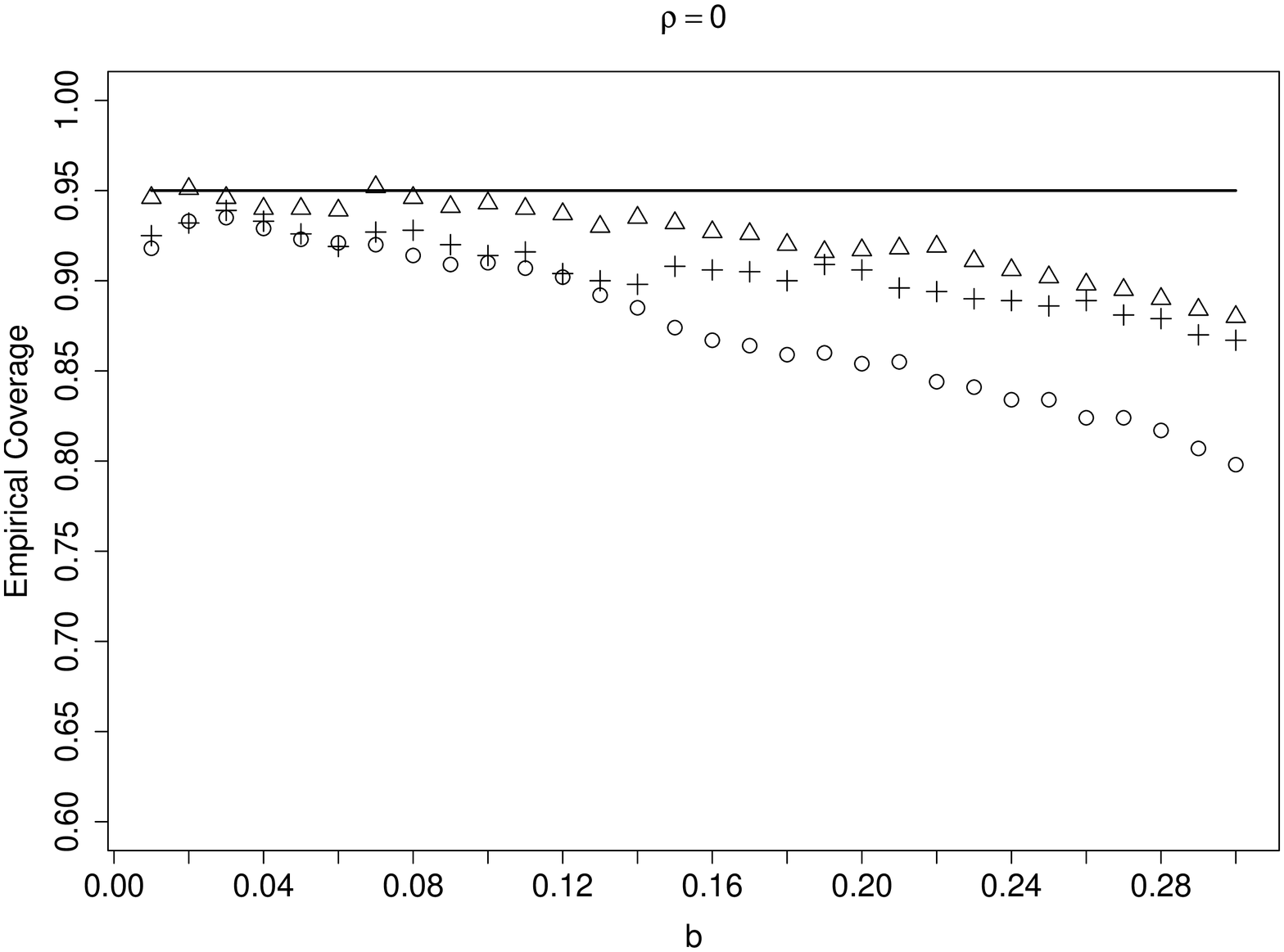}}
{\includegraphics[height=8cm,width=4.5cm,angle=270]{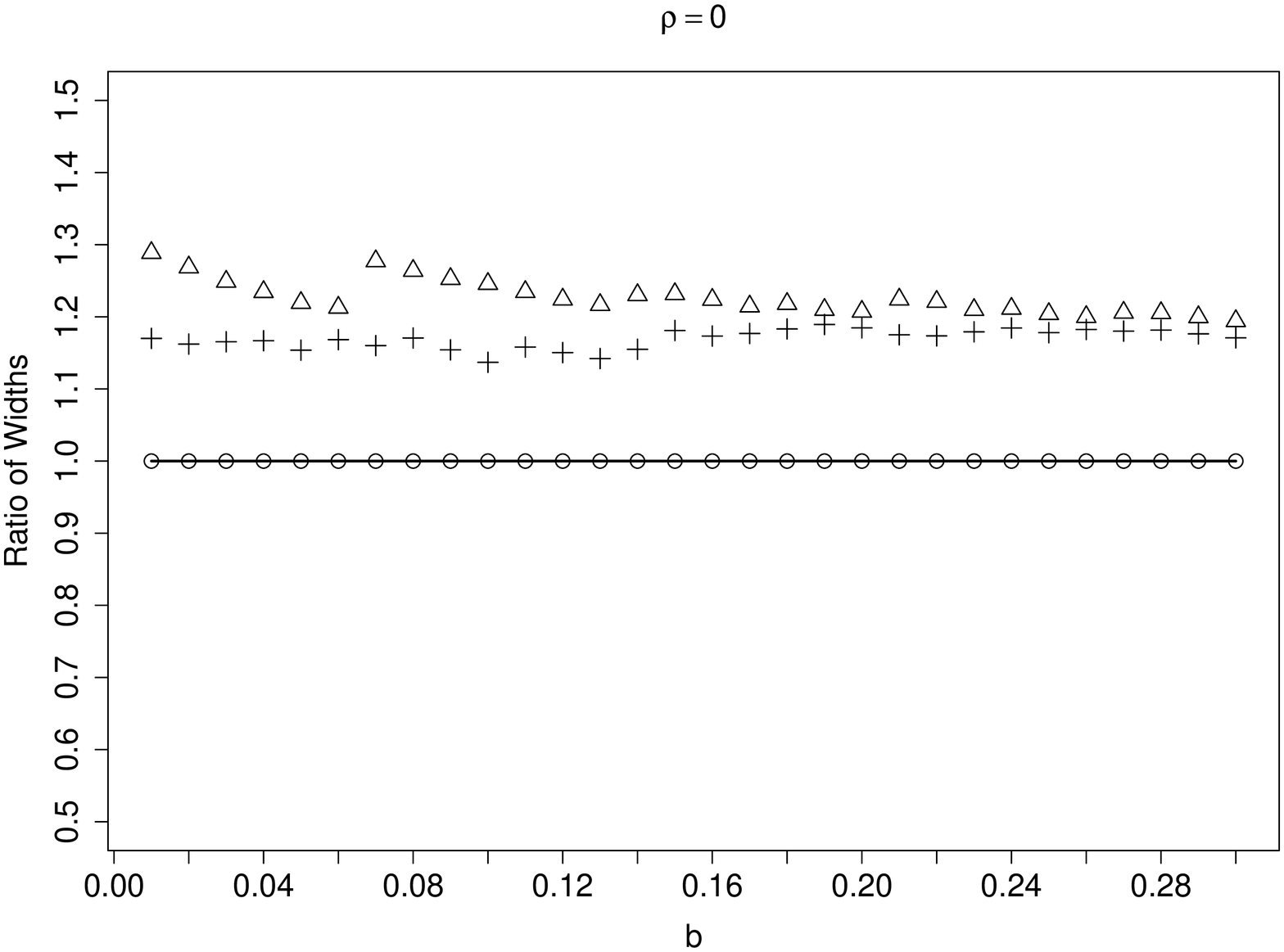}}
{\includegraphics[height=8cm,width=4.5cm,angle=270]{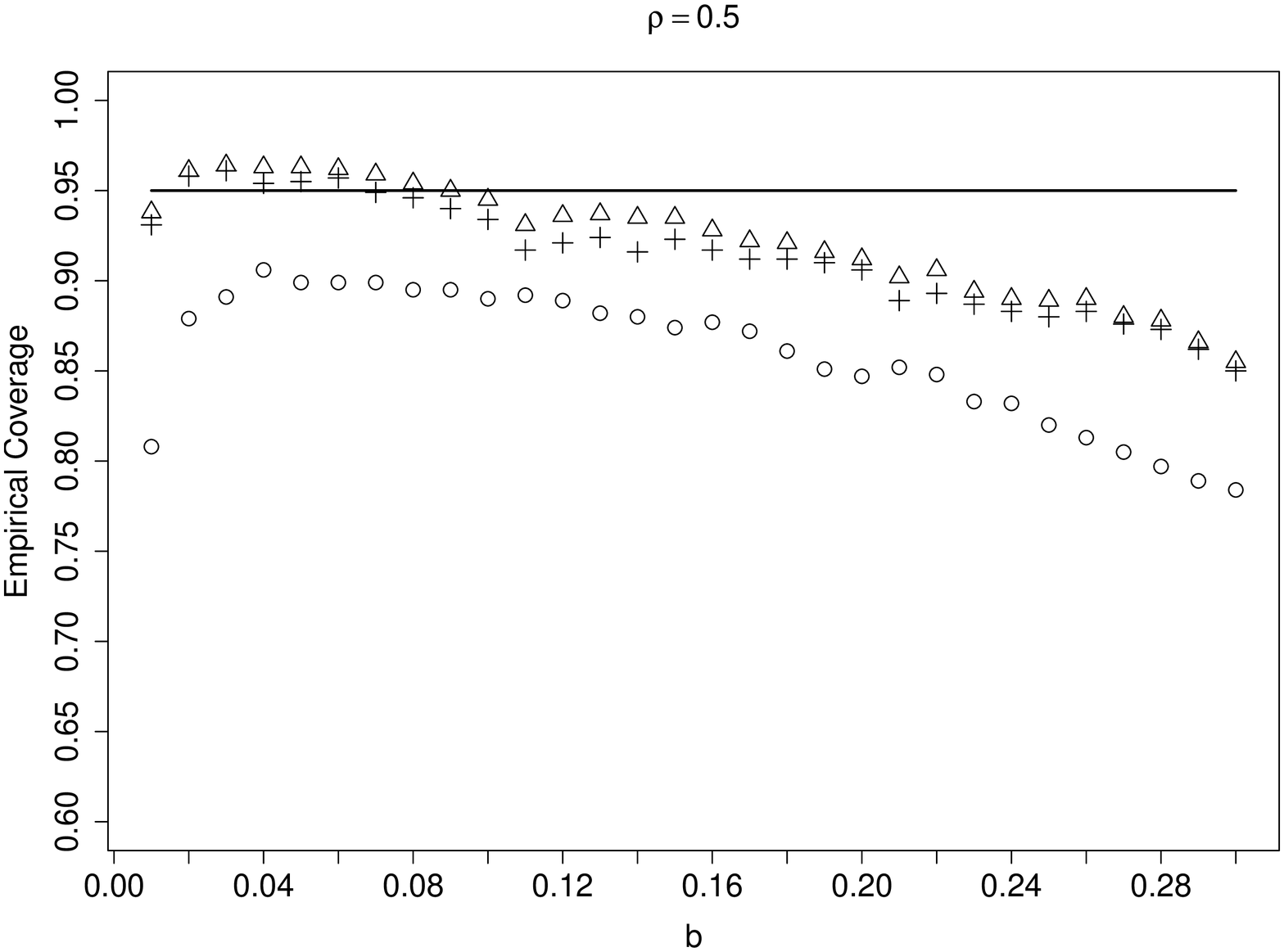}}
{\includegraphics[height=8cm,width=4.5cm,angle=270]{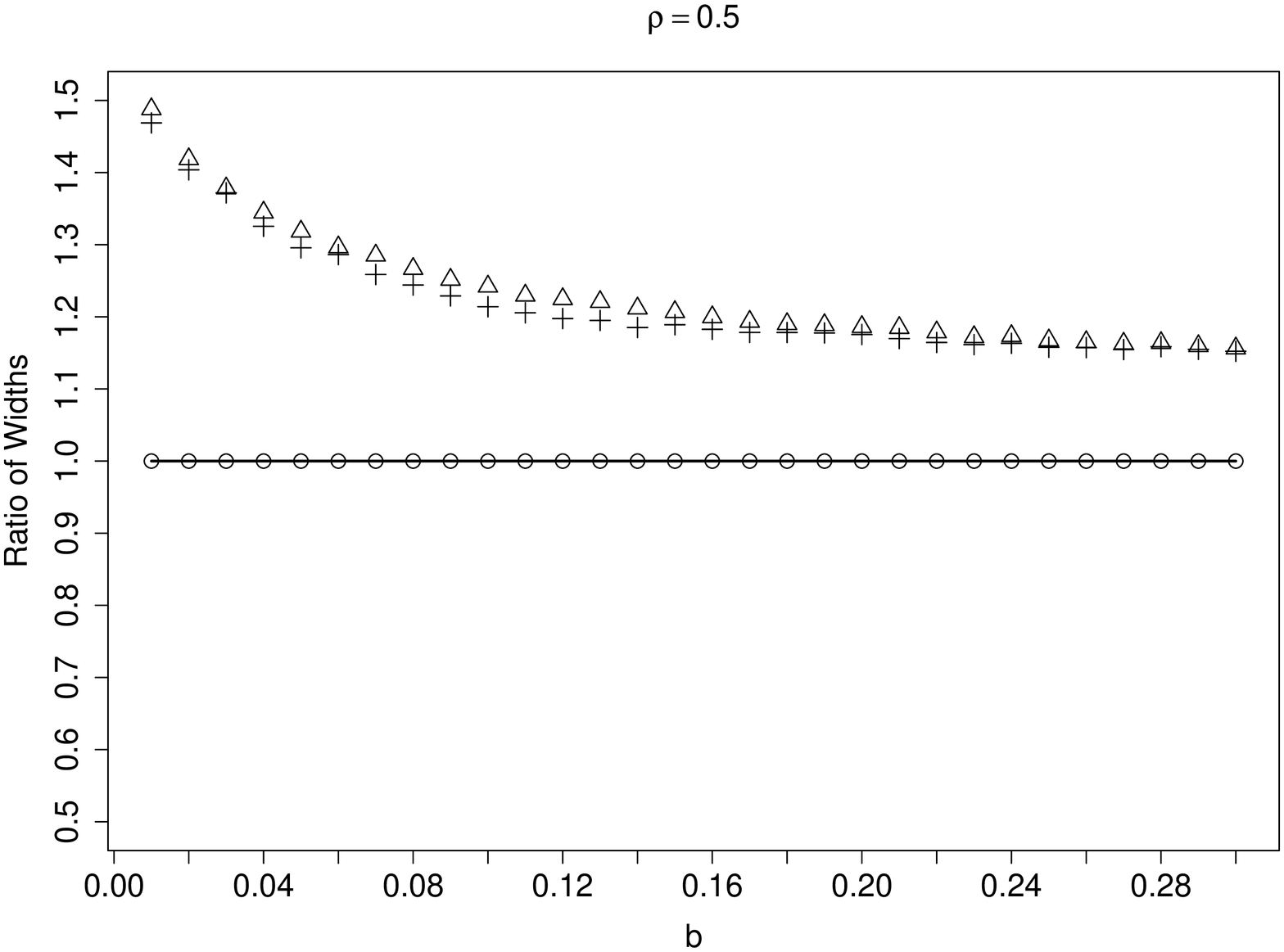}}
{\includegraphics[height=8cm,width=4.5cm,angle=270]{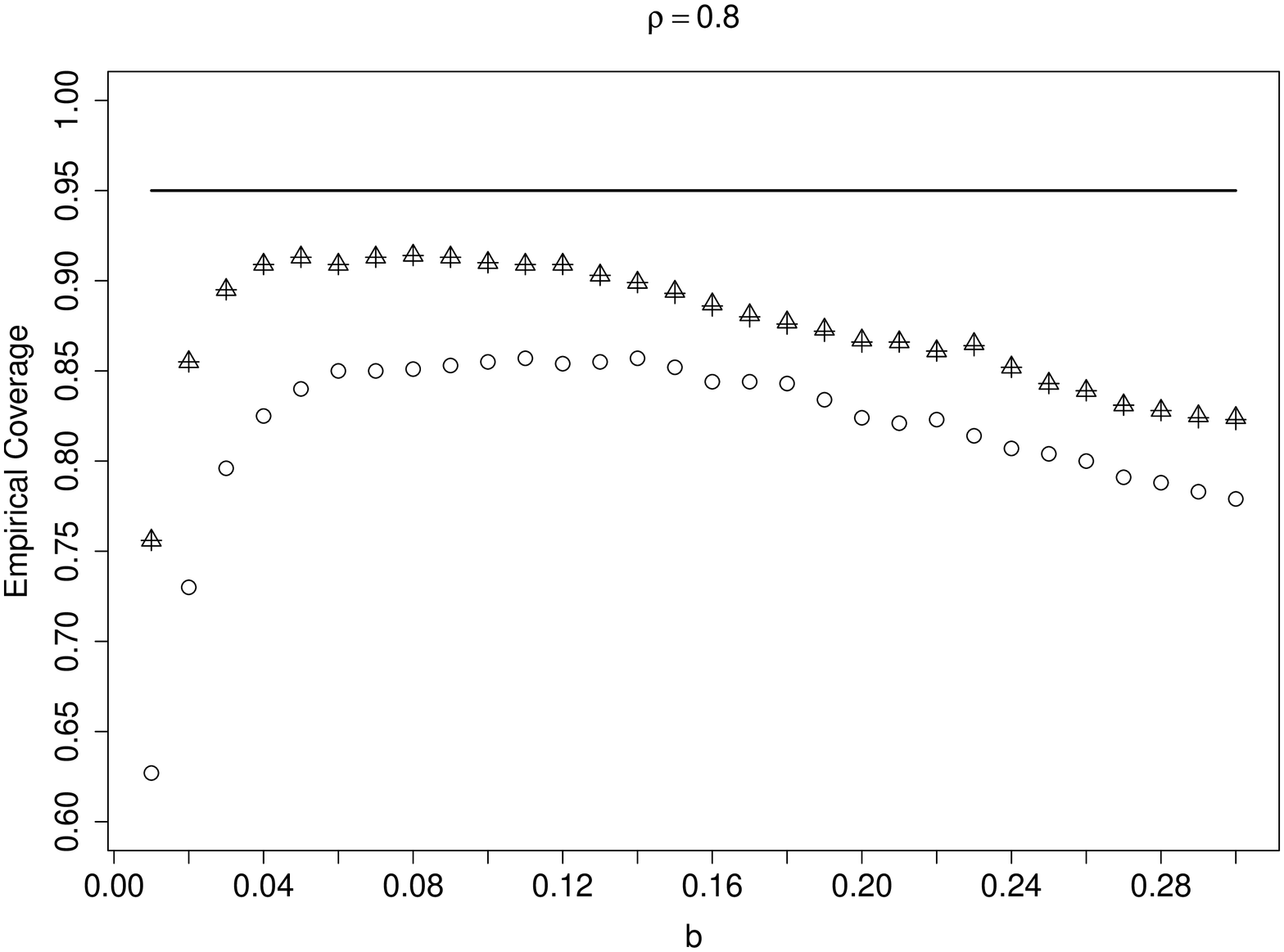}}
{\includegraphics[height=8cm,width=4.5cm,angle=270]{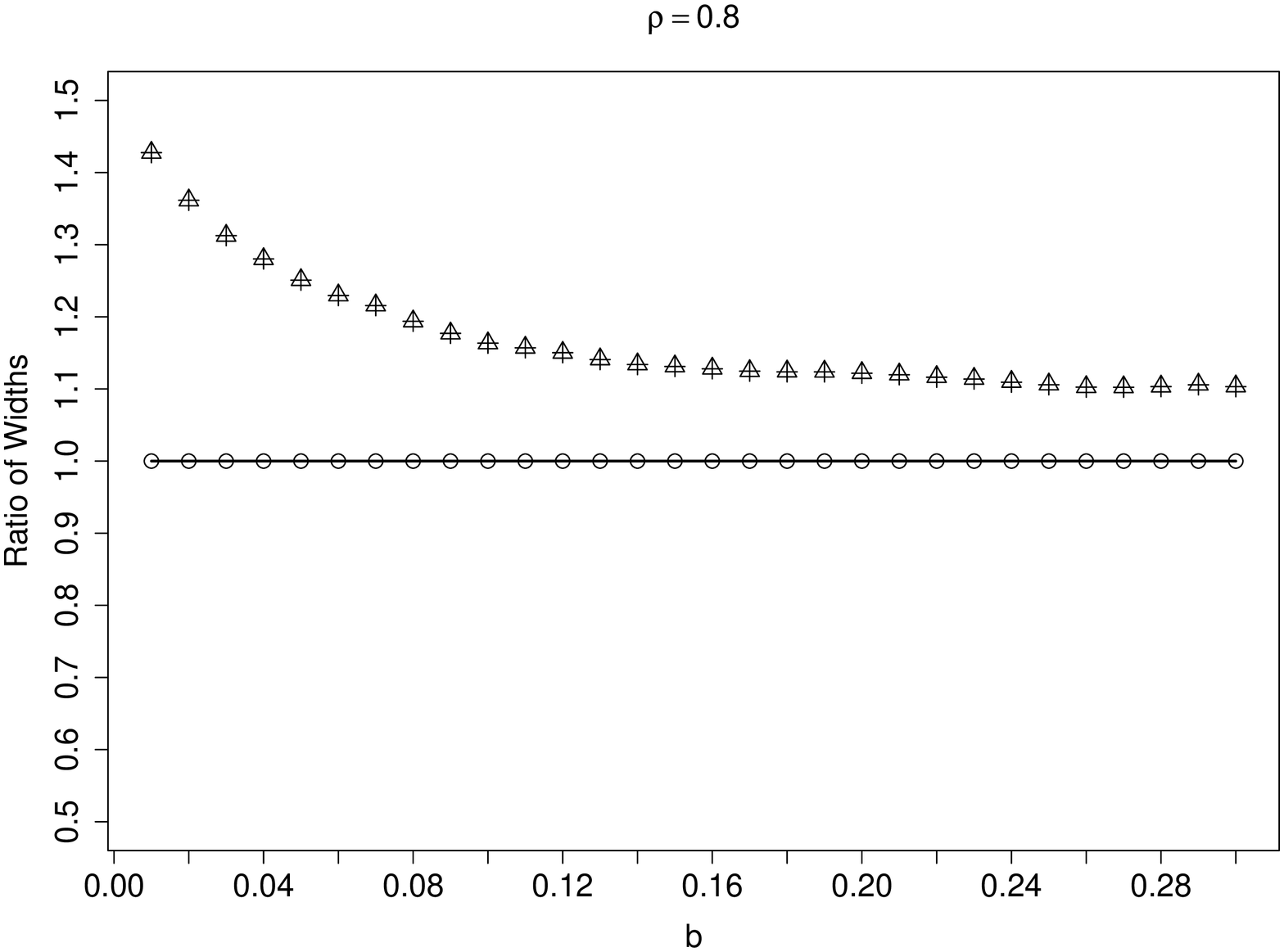}}
\label{fig:region}
\end{center}

\end{figure}

\begin{figure}
\caption{The empirical coverage probabilities (left panel) and  the ratios of radii of confidence regions (calibrated fixed-$b$  over traditional small-$b$)
 (right panel) for the vector parameter and for the models with exponentially distributed errors. Sample size $n=200$ and number of replications is 1000. }
\begin{center}
{\includegraphics[height=8cm,width=4.5cm,angle=270]{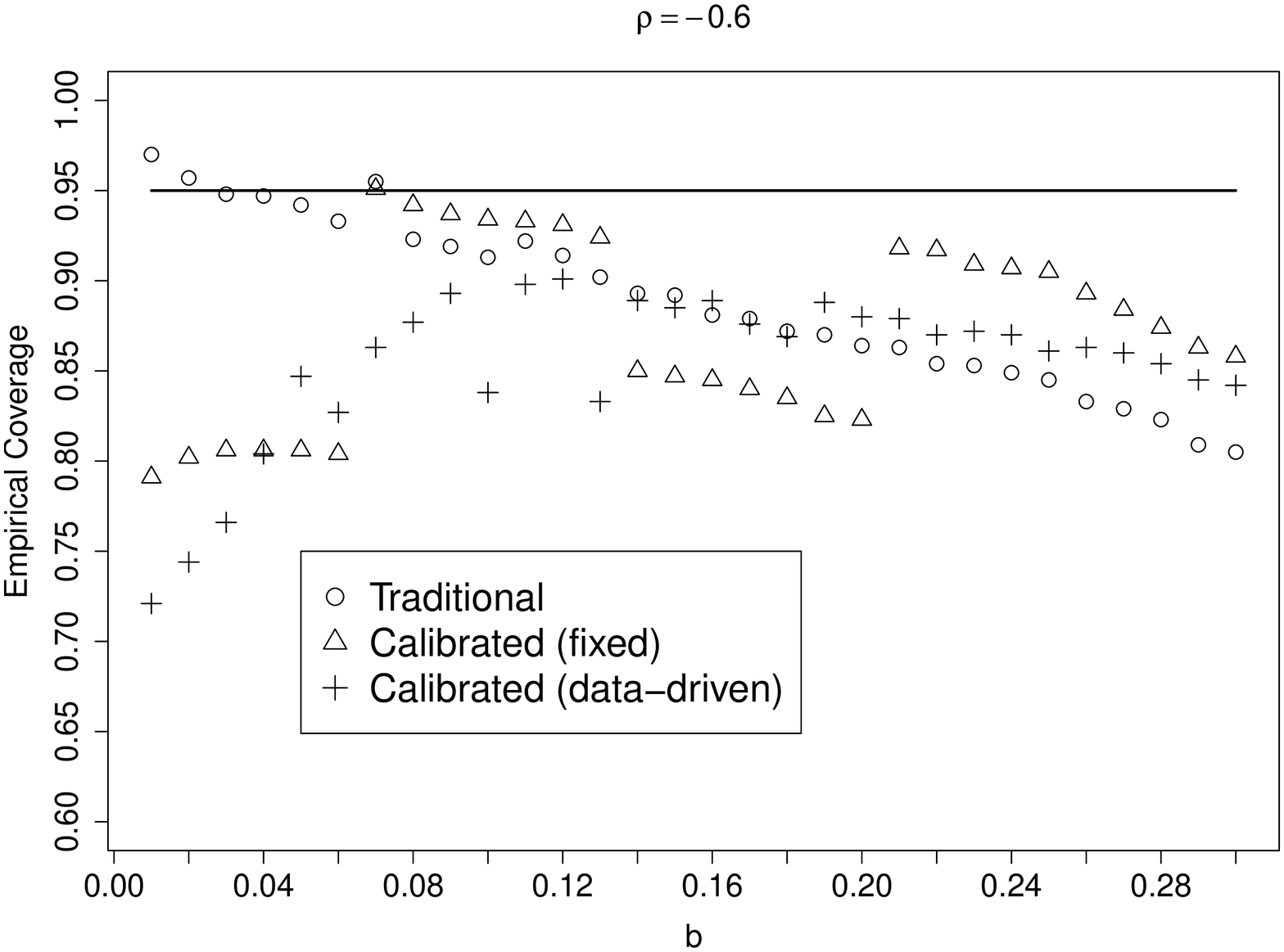}}
{\includegraphics[height=8cm,width=4.5cm,angle=270]{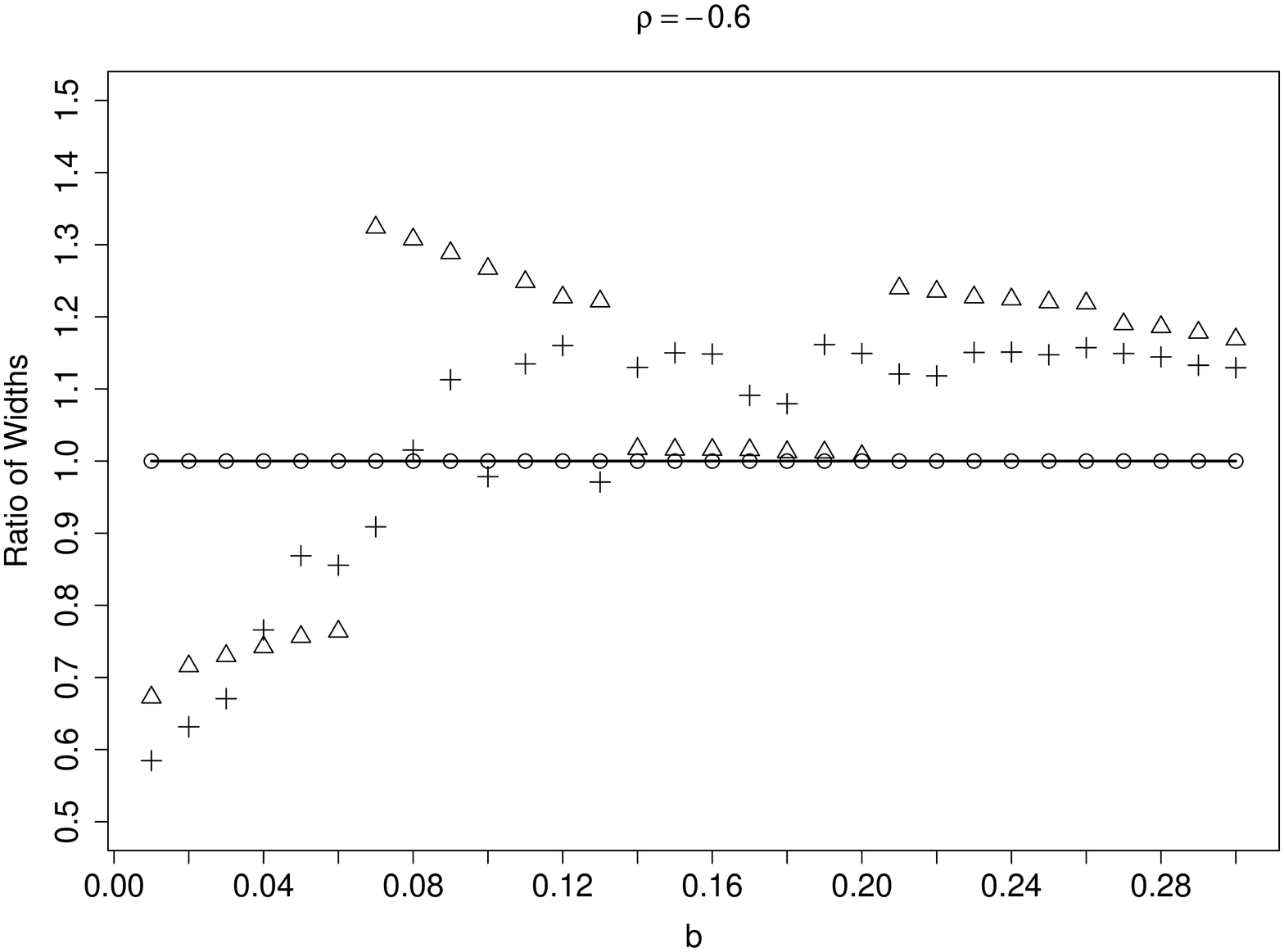}}
{\includegraphics[height=8cm,width=4.5cm,angle=270]{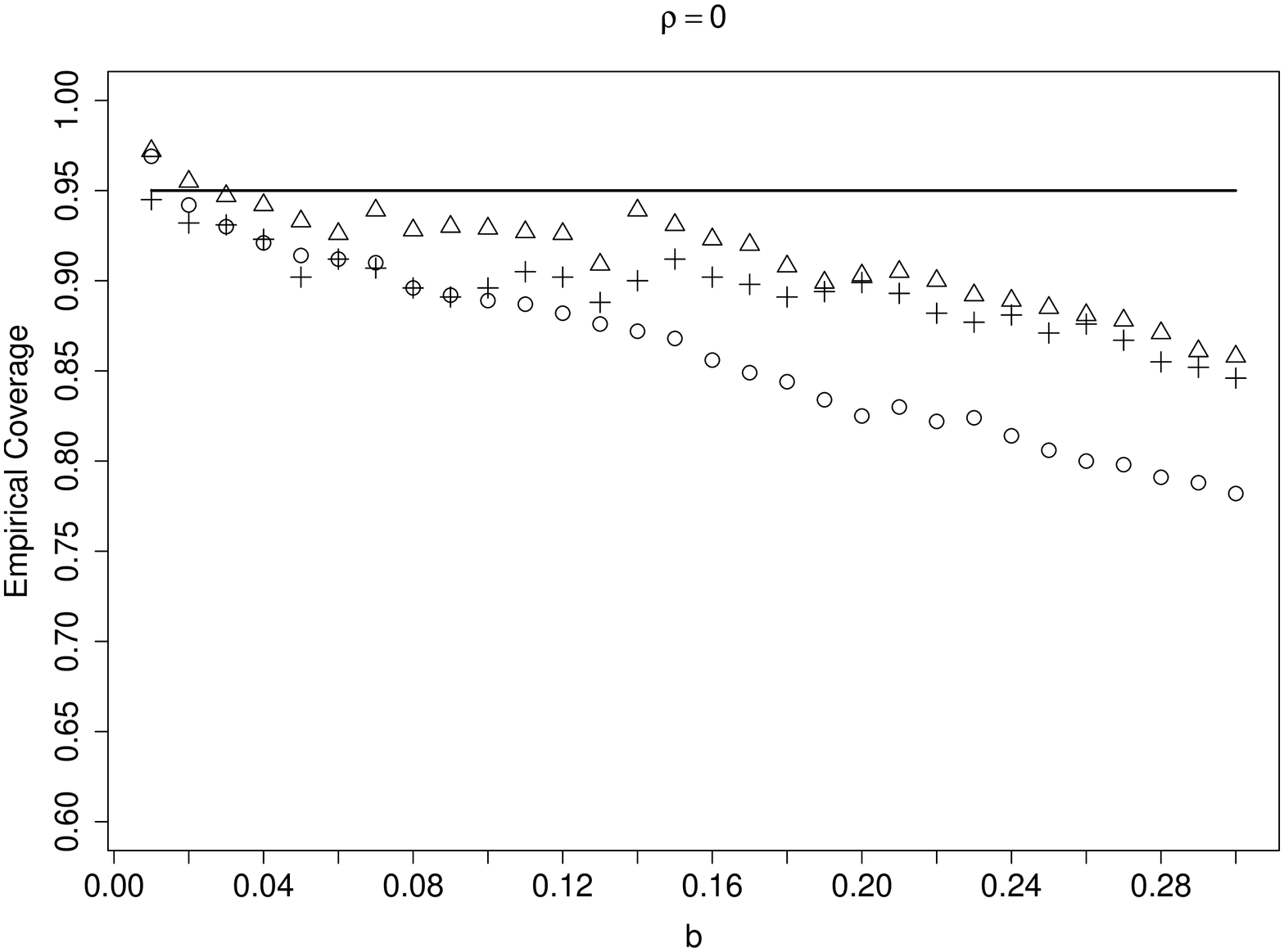}}
{\includegraphics[height=8cm,width=4.5cm,angle=270]{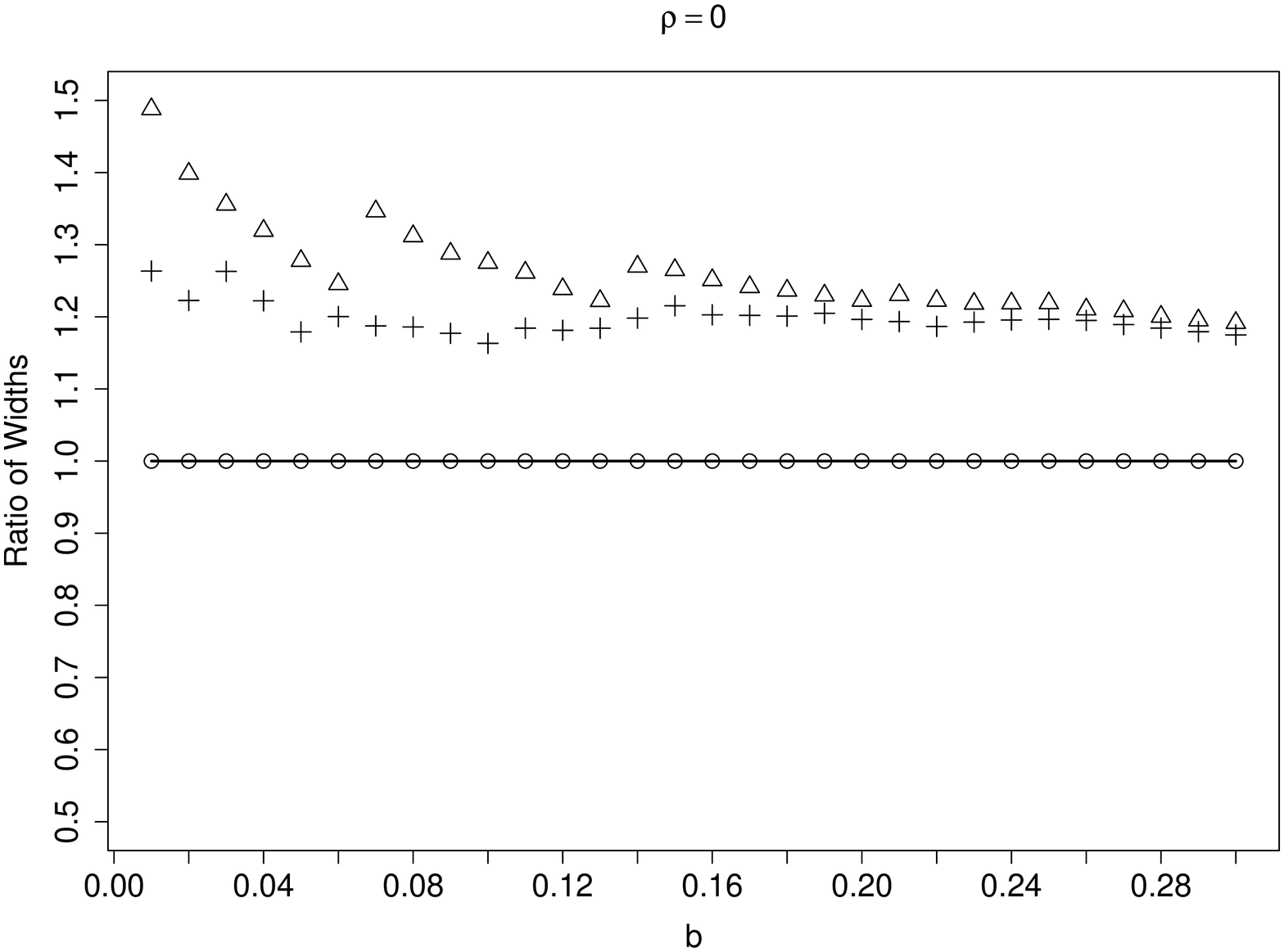}}
{\includegraphics[height=8cm,width=4.5cm,angle=270]{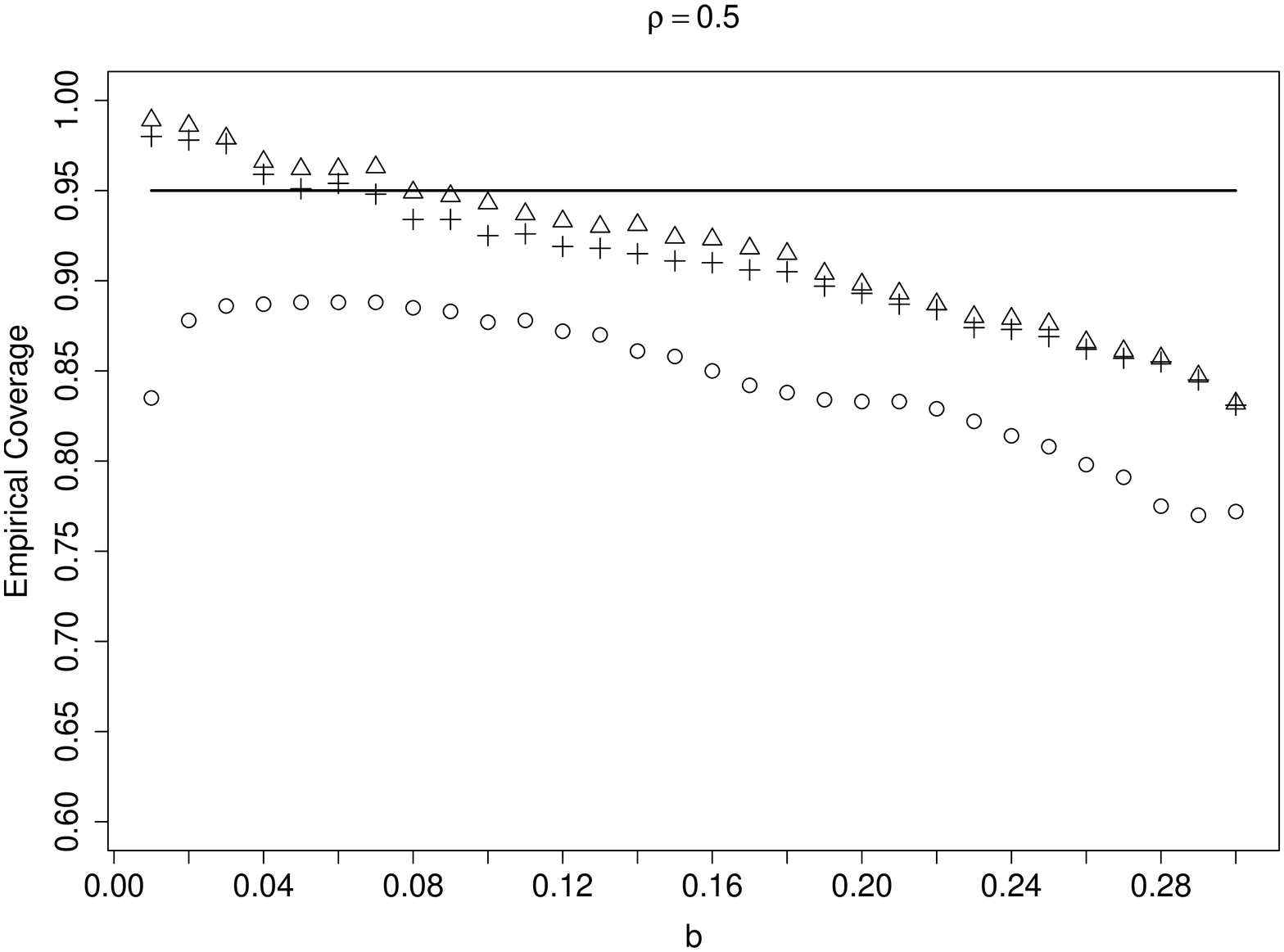}}
{\includegraphics[height=8cm,width=4.5cm,angle=270]{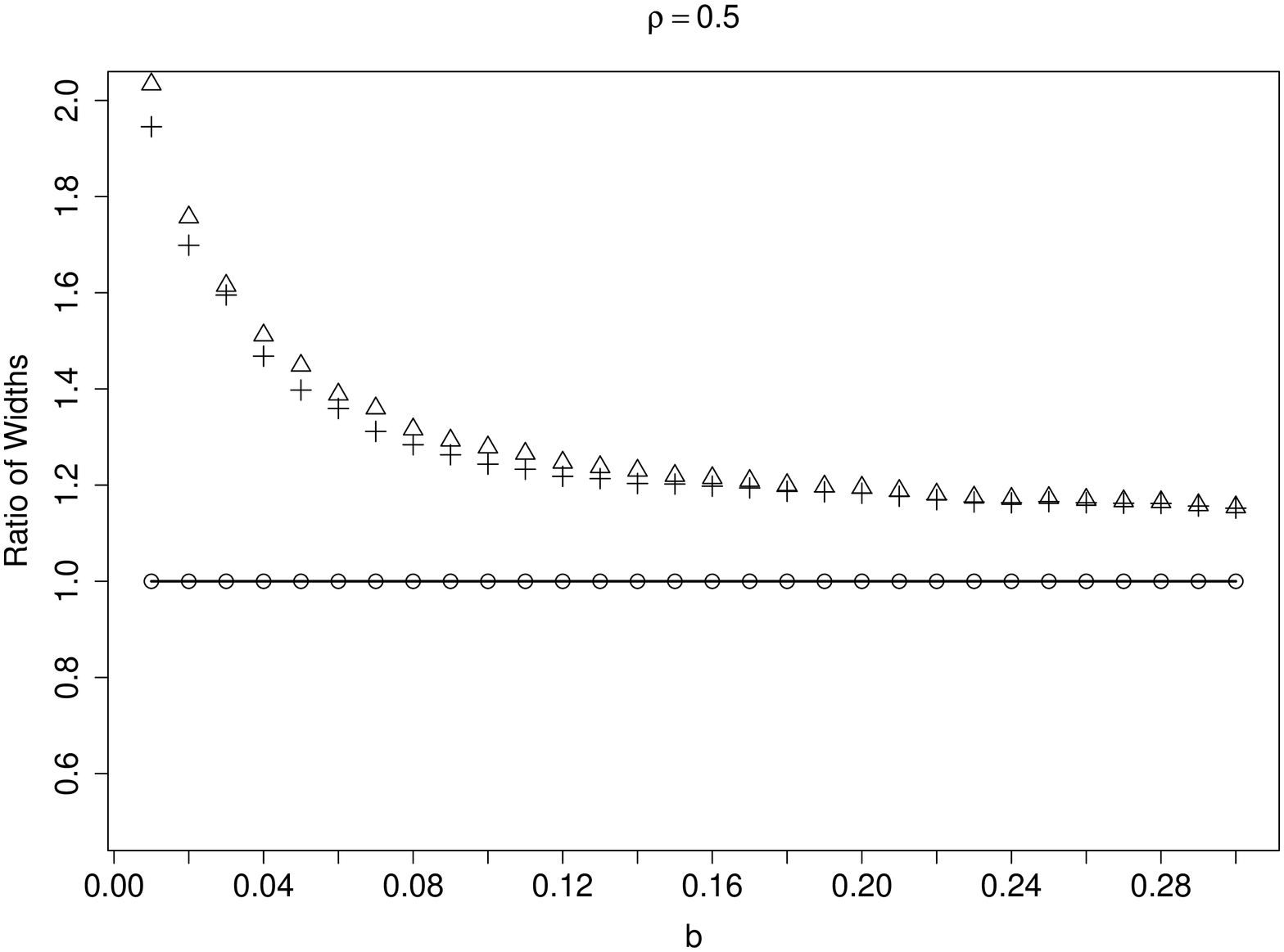}}
{\includegraphics[height=8cm,width=4.5cm,angle=270]{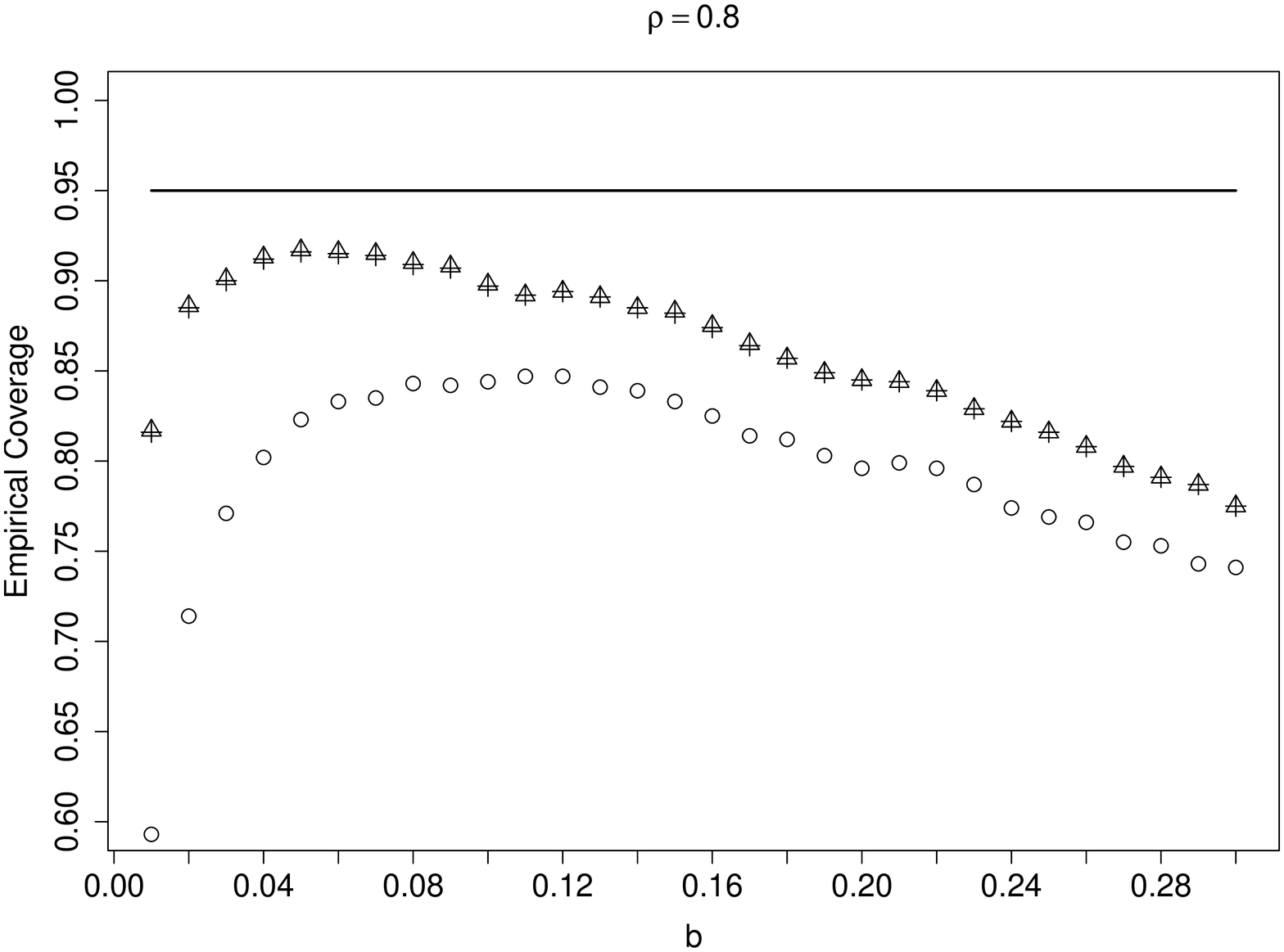}}
{\includegraphics[height=8cm,width=4.5cm,angle=270]{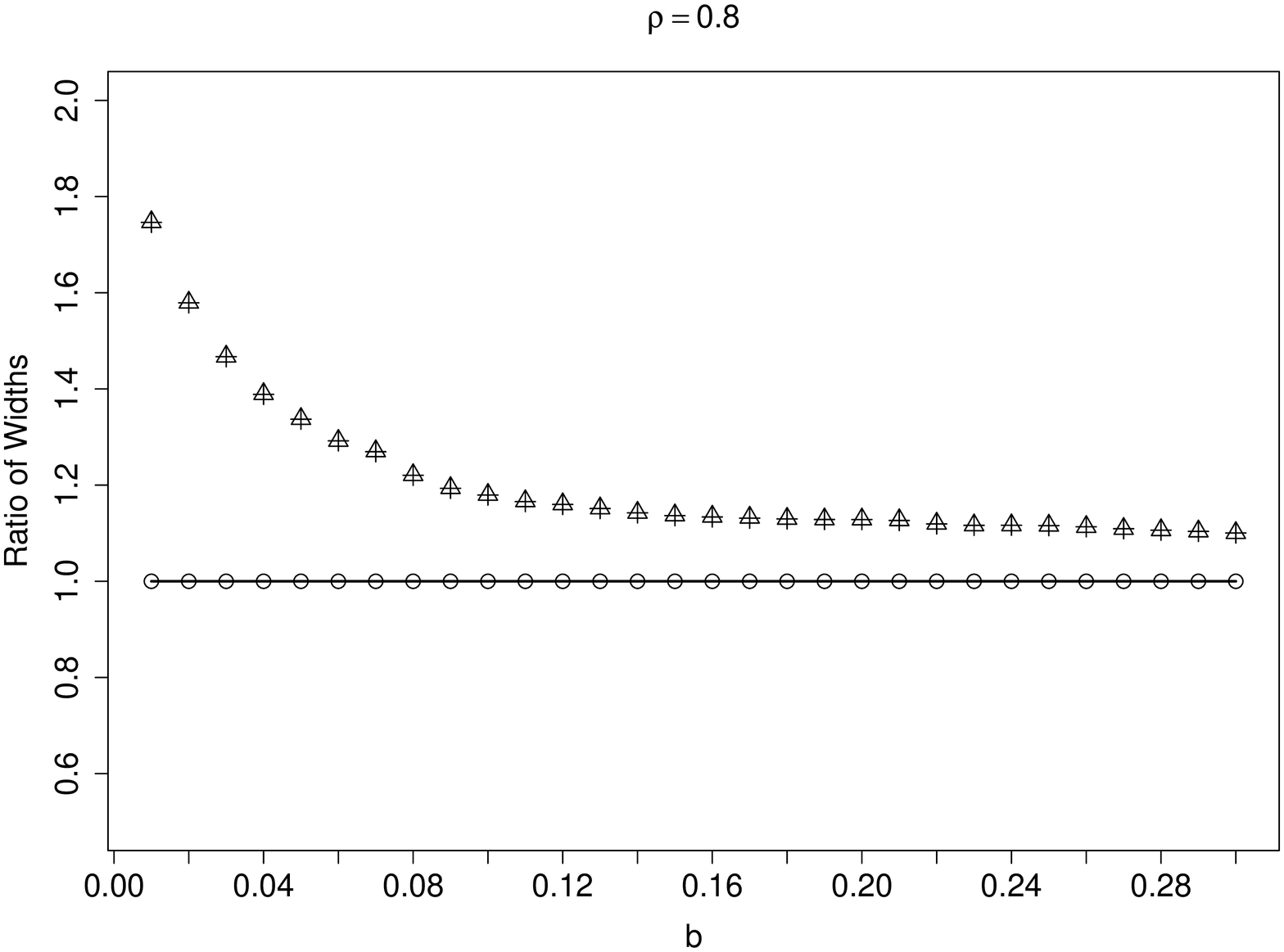}}
\label{fig:regionEXP}
\end{center}

\end{figure}

Figures~\ref{fig:region} and ~\ref{fig:regionEXP} depict the empirical coverages and the ratios of the radii of the confidence regions over that
delivered by the uncalibrated traditional subsampling-based region for the vector parameter and for the models with normally distributed errors and exponentially distributed errors, respectively.
The symbols ``Traditional", ``Calibrated (fixed)" and ``Calibrated (data-driven)" correspond to
  the schemes (1)-(3) described above. As the findings for the normally distributed case and the exponentially distributed case are very close, we shall only describe the results for the normally distributed case.
   When $\rho=0, 0.5,0.8$, there is an undercoverage associated with the traditional subsampling-based approach
  and the coverage errors increase with respect to the magnitude of dependence. The improvement in coverage offered by the calibration is apparent in these cases and it
  holds uniformly over the range of $b$s under examination. On the other hand,    the corresponding radius of the calibrated region is slightly
  larger than that of the uncalibrated counterpart. In the case $\rho=-0.6$, the calibrated region performs worse compared to the traditional counterpart when $b$ is small, but
  still offer some improvement when $b$ is large. It is not fully clear why this phenomenon occurs.   Nevertheless, it suggests that caution has to be exercised in the use of fixed-$b$ based calibration when the autocorrelations of the series have alternating signs.

\begin{figure}
\caption{The empirical coverage probabilities (left panel) and  the ratios of band widths (calibrated fixed-$b$  over traditional small-$b$)
 (right panel) for the marginal distribution function. Sample size $n=200$ and number of replications is 1000. }
\begin{center}
{\includegraphics[height=8cm,width=4.5cm,angle=270]{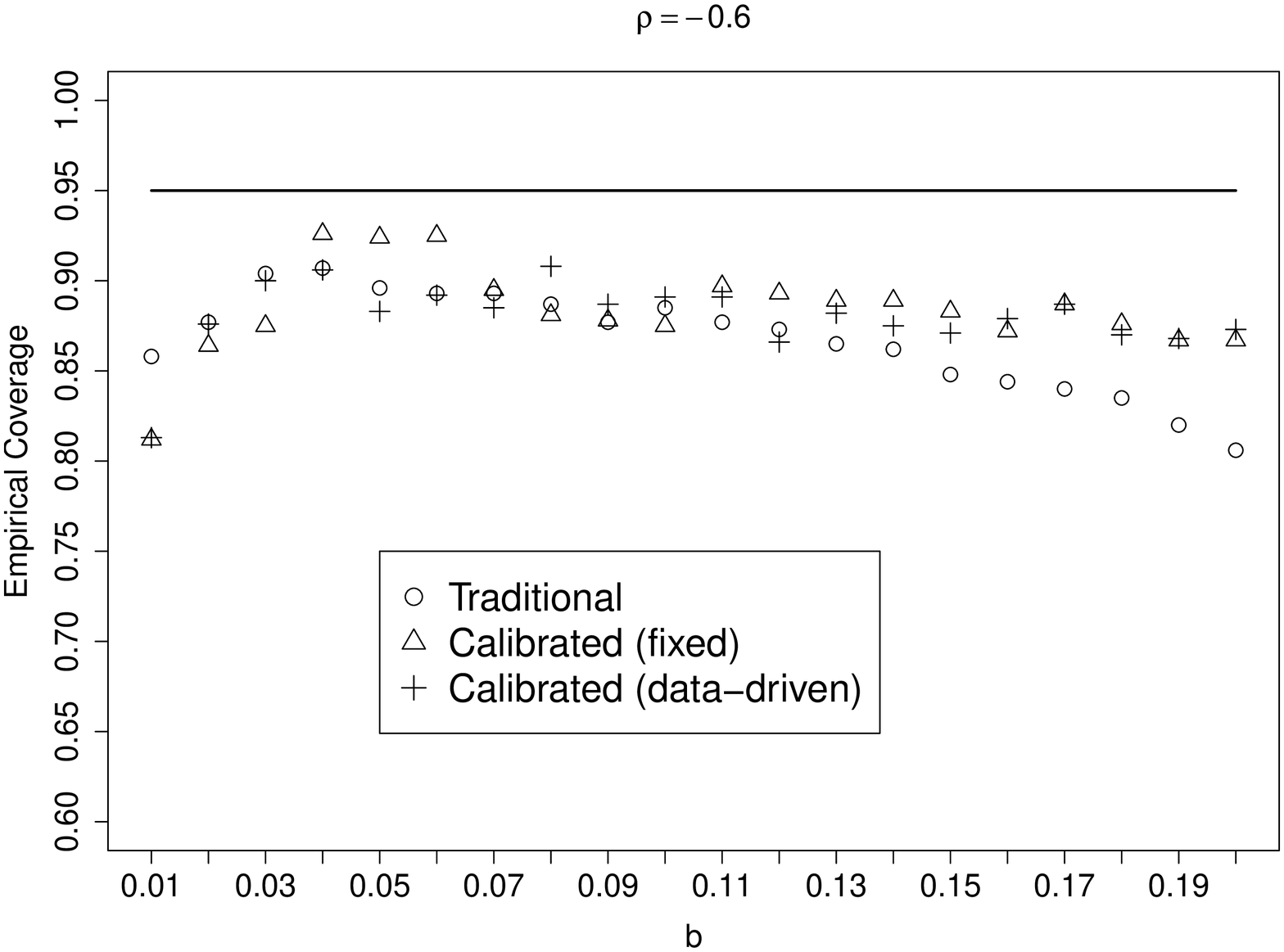}}
{\includegraphics[height=8cm,width=4.5cm,angle=270]{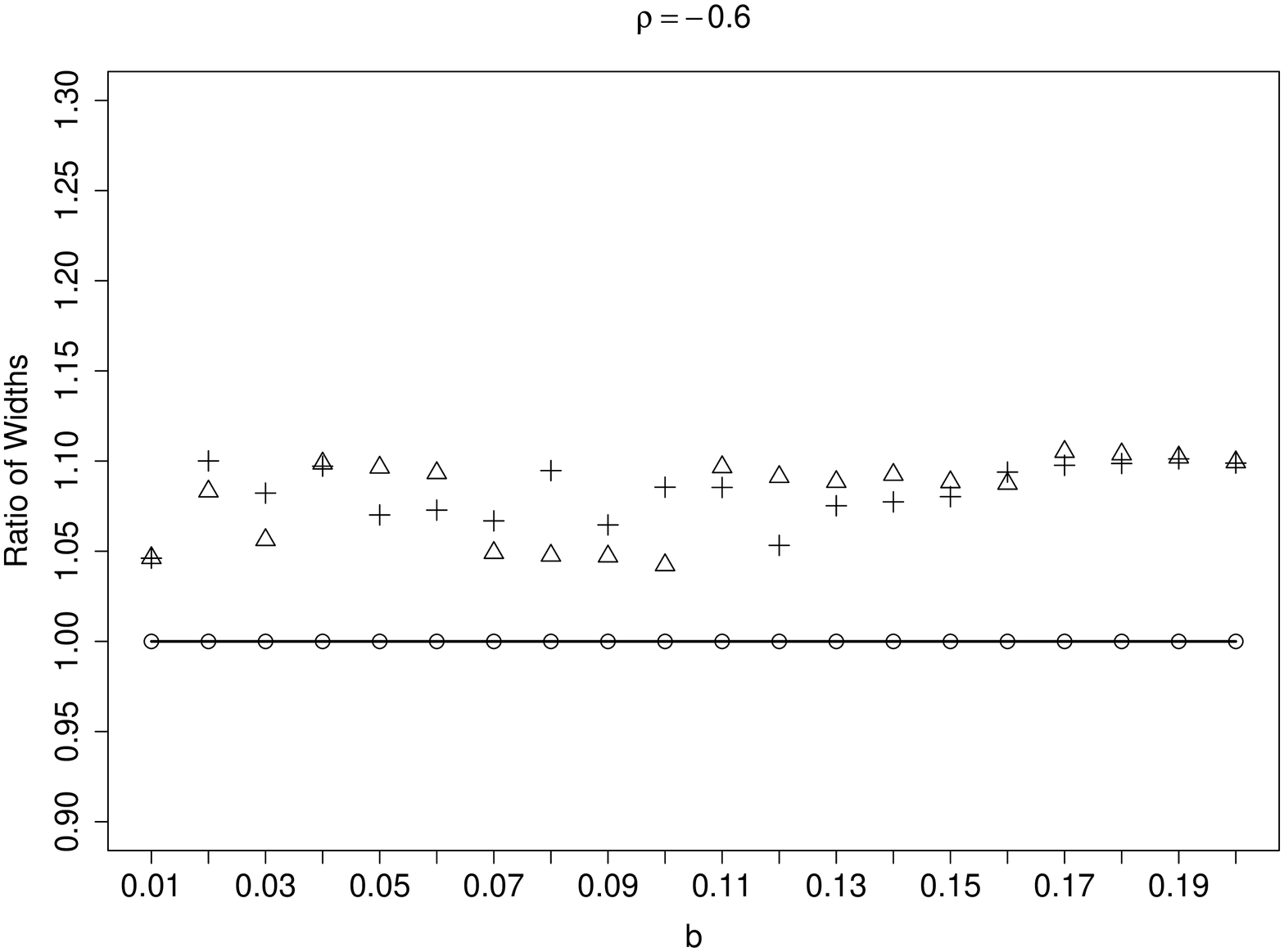}}
{\includegraphics[height=8cm,width=4.5cm,angle=270]{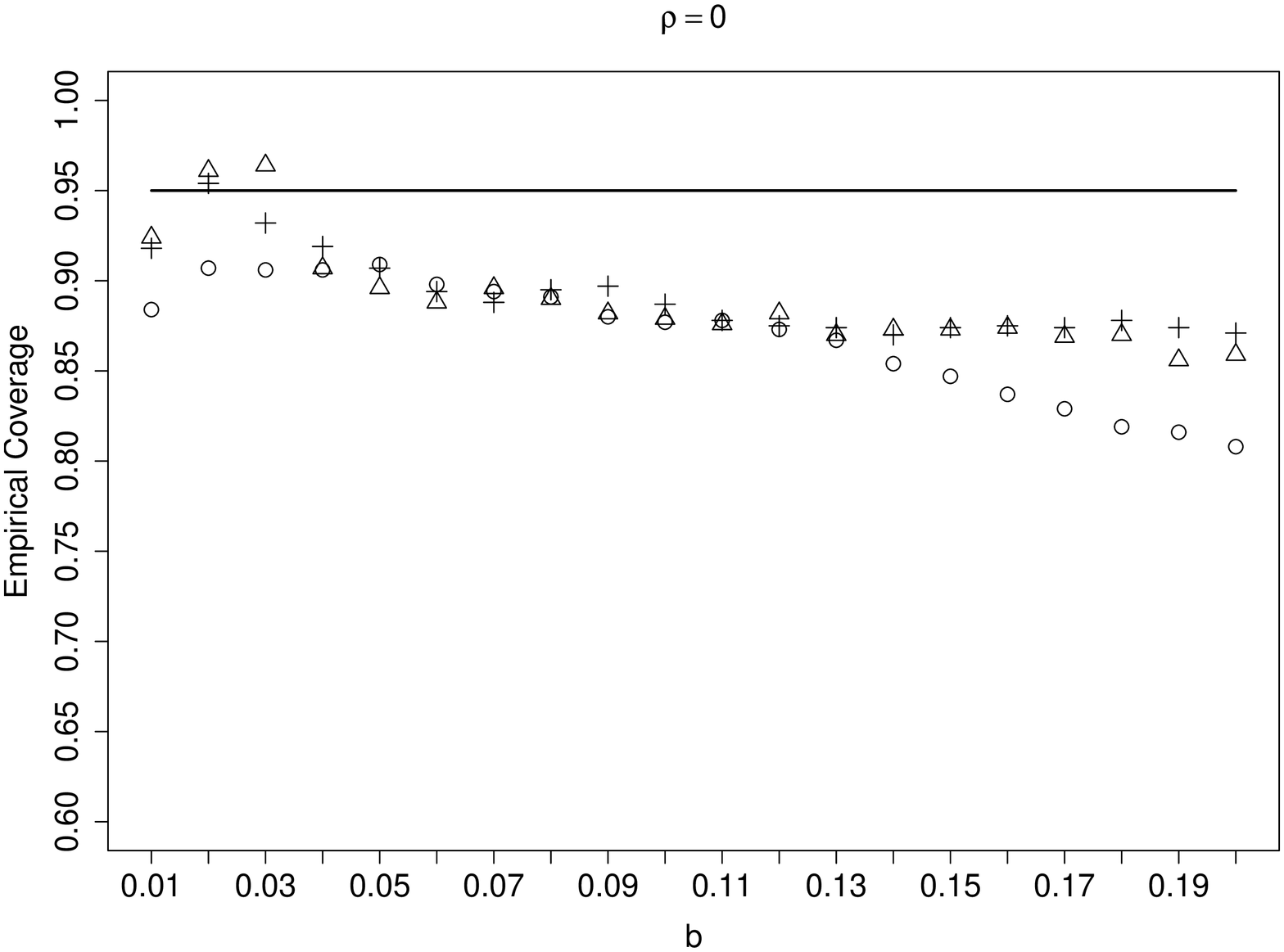}}
{\includegraphics[height=8cm,width=4.5cm,angle=270]{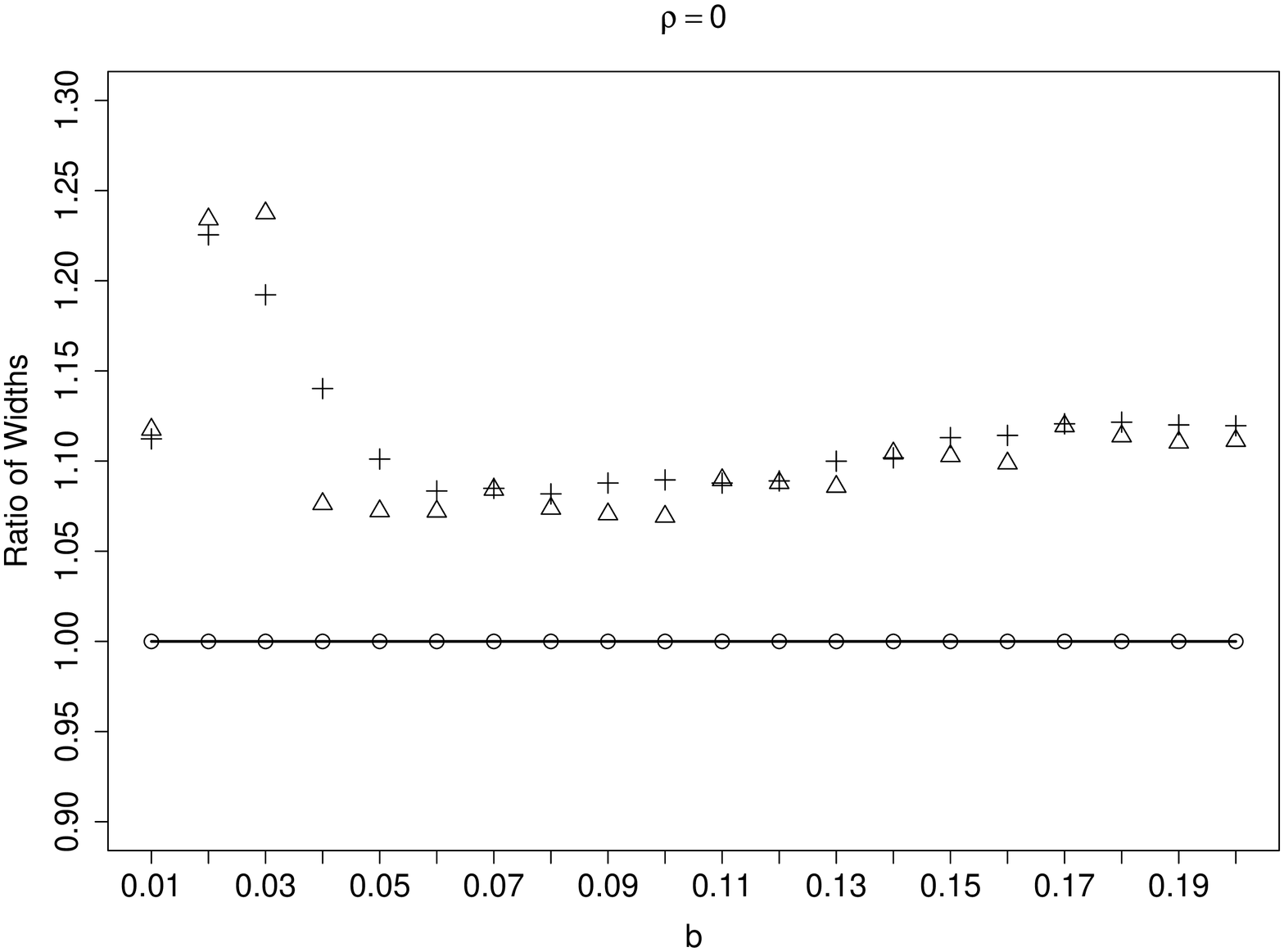}}
{\includegraphics[height=8cm,width=4.5cm,angle=270]{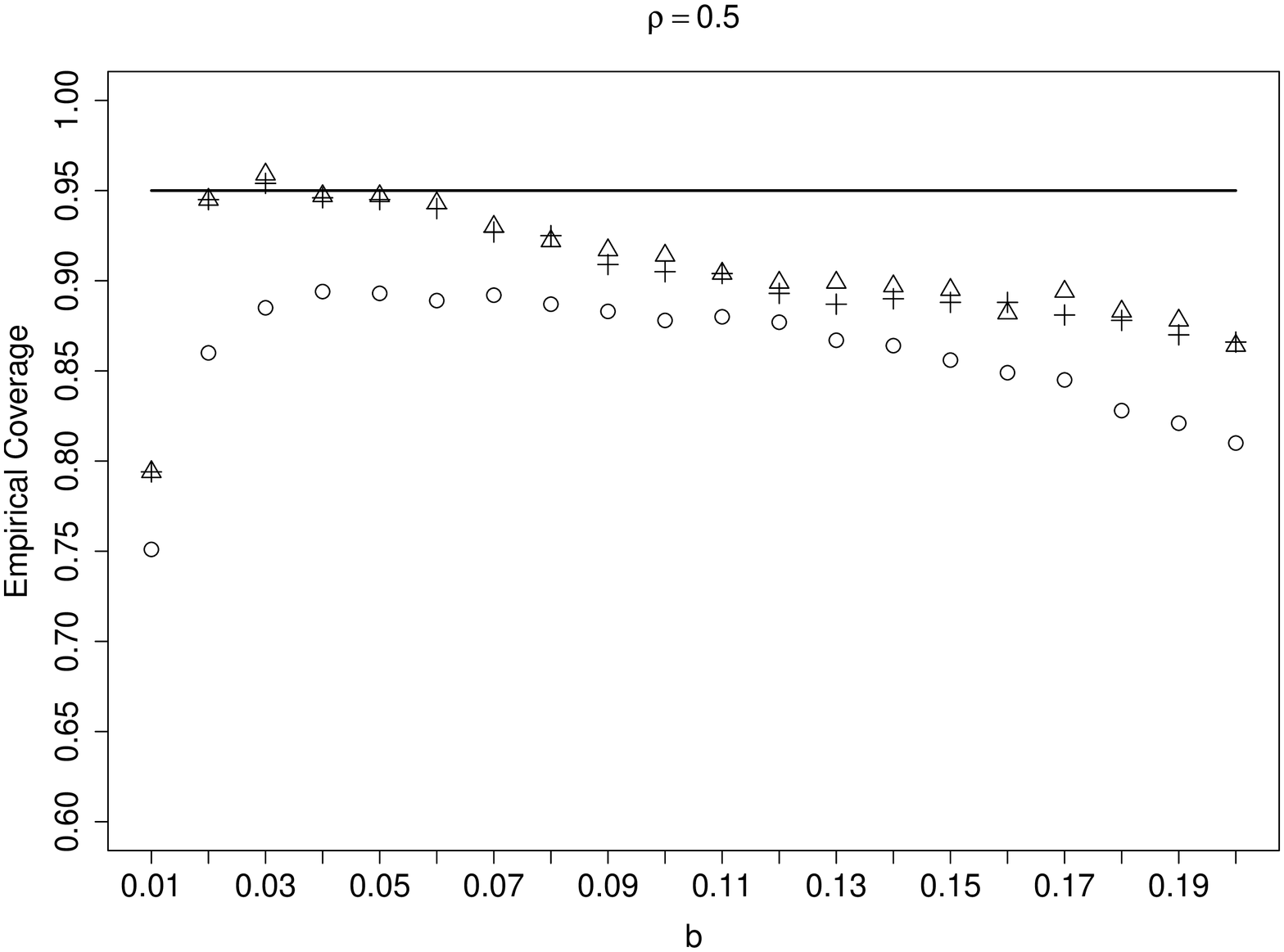}}
{\includegraphics[height=8cm,width=4.5cm,angle=270]{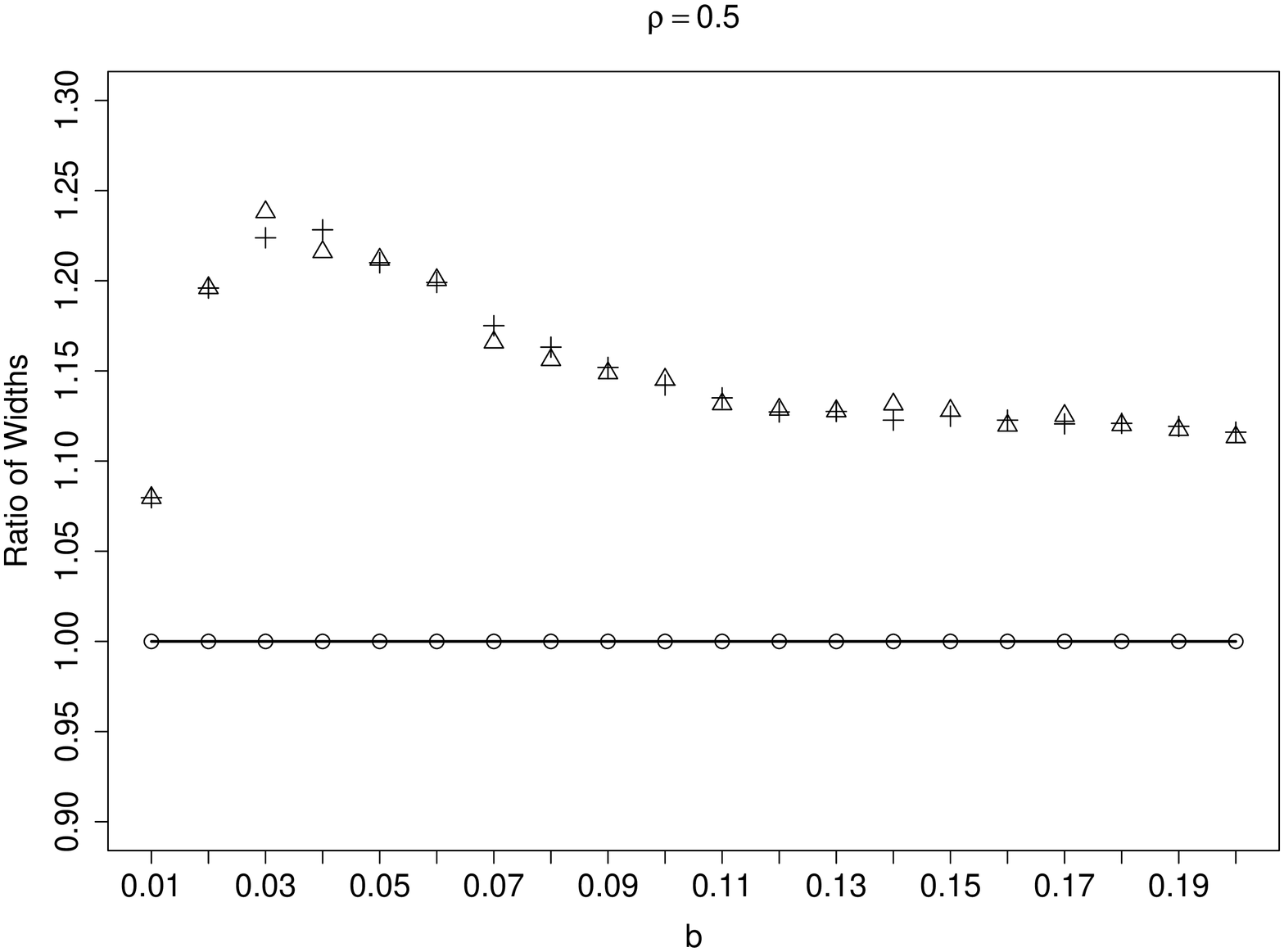}}
{\includegraphics[height=8cm,width=4.5cm,angle=270]{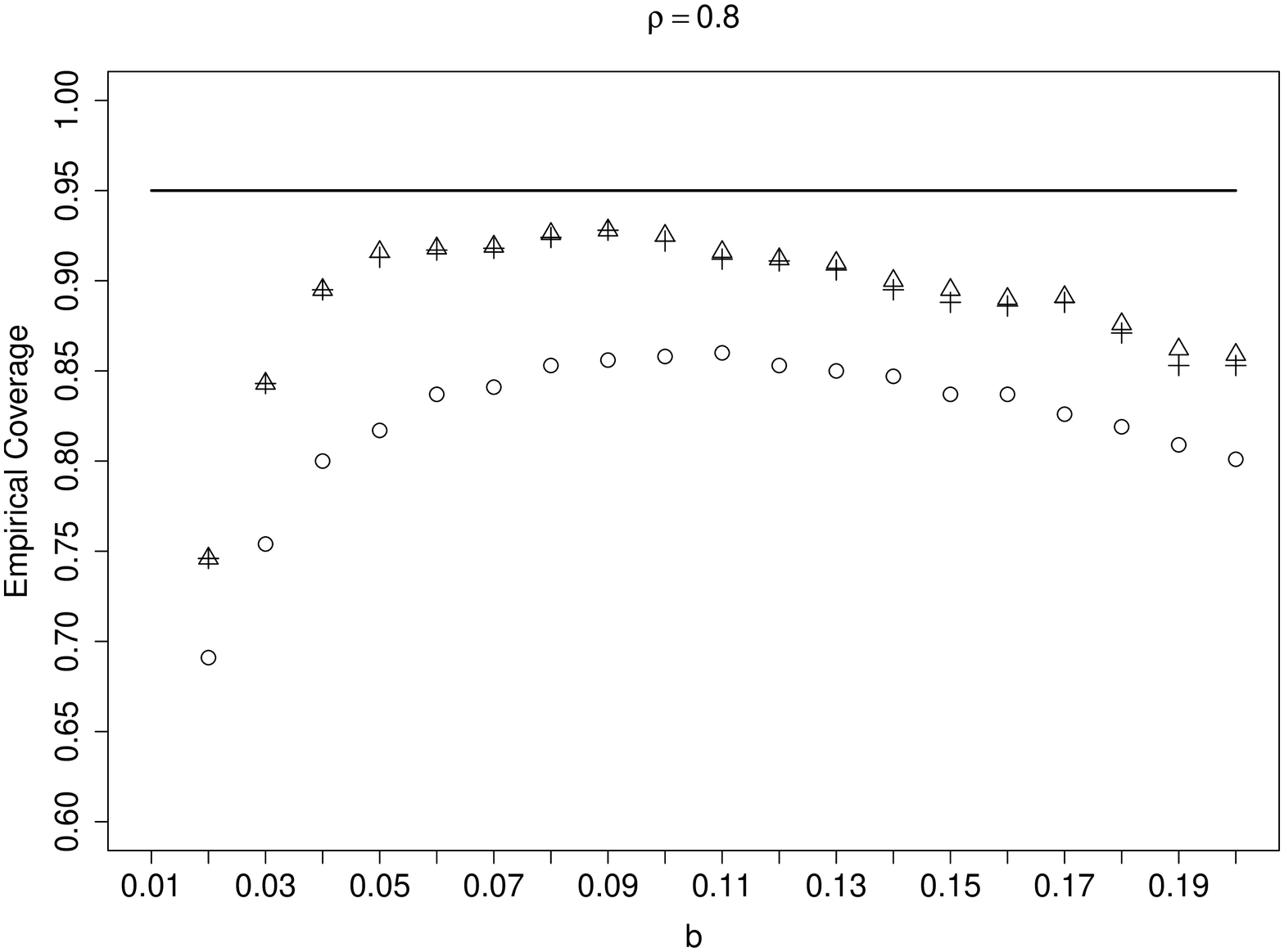}}
{\includegraphics[height=8cm,width=4.5cm,angle=270]{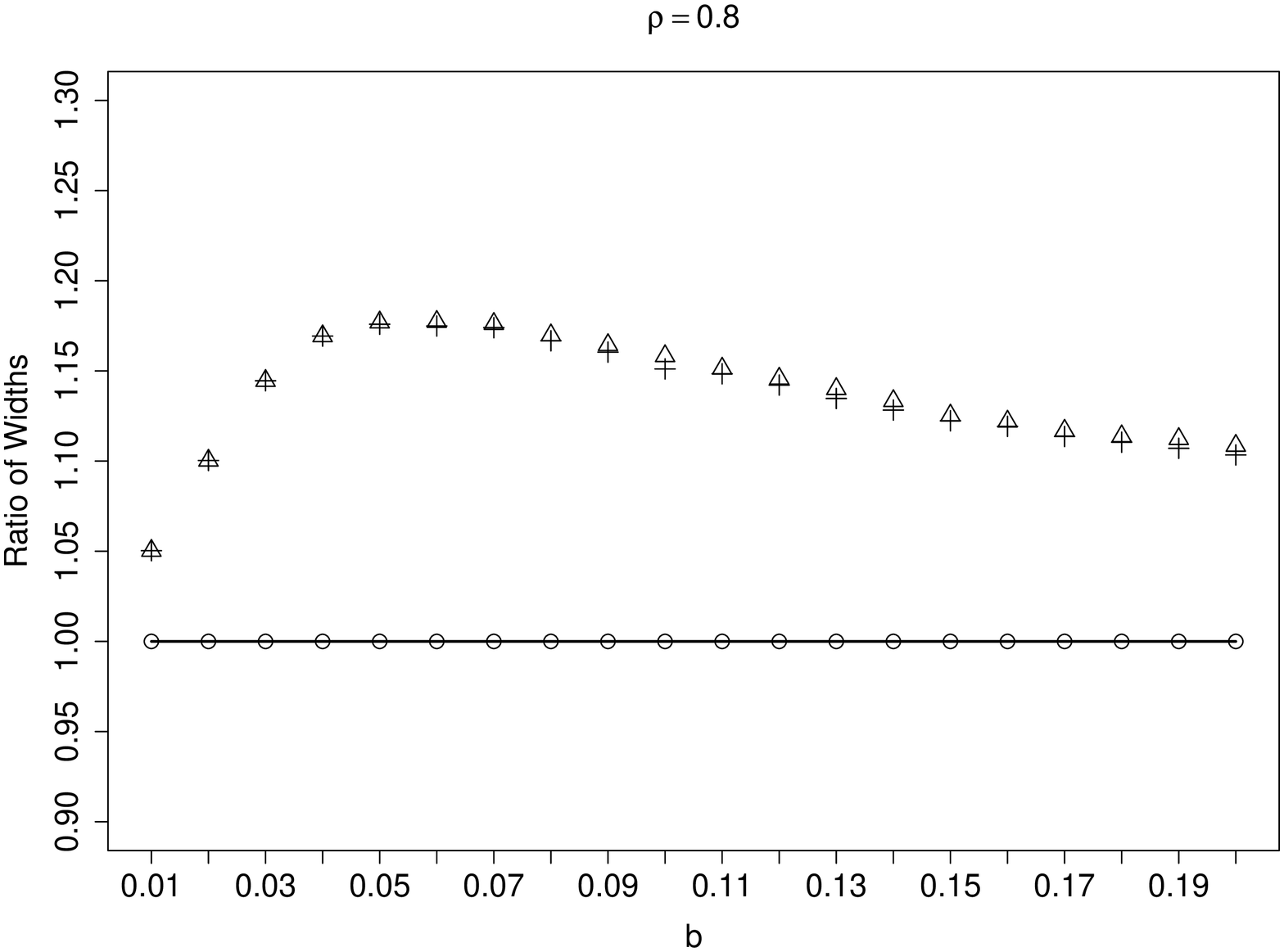}}
\label{fig:band1}
\end{center}

\end{figure}

\begin{figure}
\caption{The empirical coverage probabilities (left panel) and  the ratios of band widths (calibrated fixed-$b$  over traditional small-$b$)
 (right panel) for the normalized spectral distribution function. Sample size $n=200$ and number of replications is 1000. }
\begin{center}
{\includegraphics[height=8cm,width=4.5cm,angle=270]{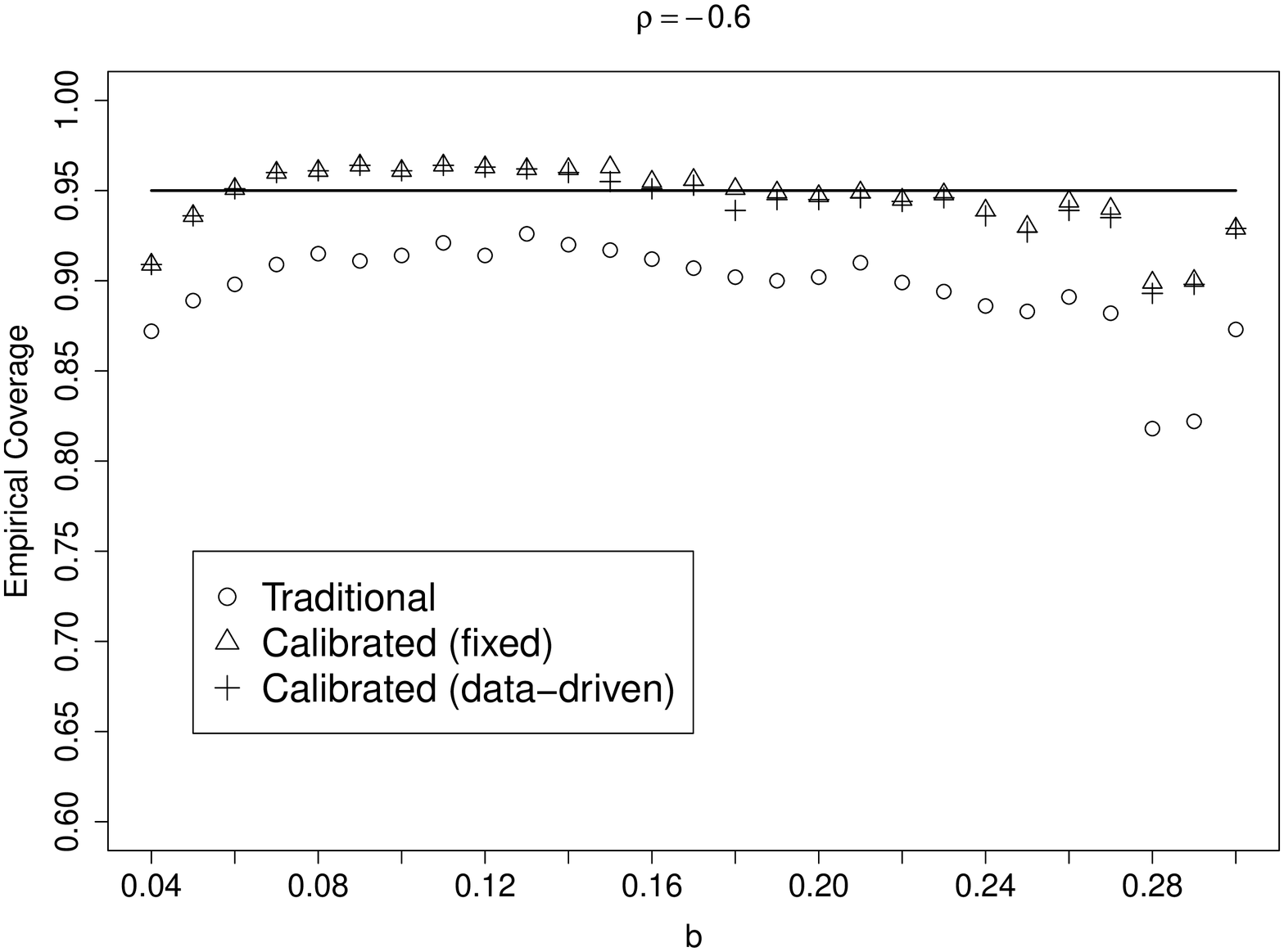}}
{\includegraphics[height=8cm,width=4.5cm,angle=270]{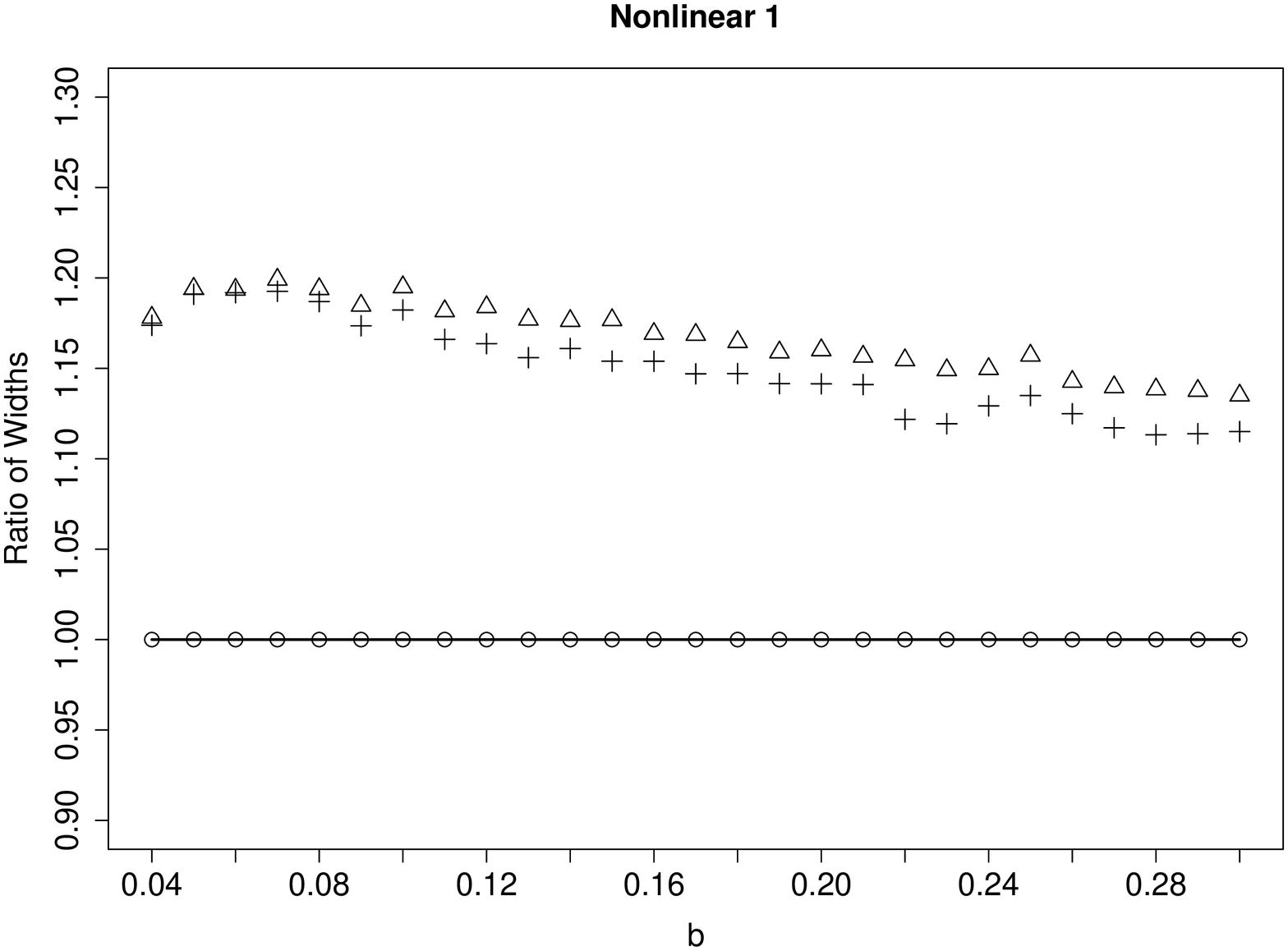}}
{\includegraphics[height=8cm,width=4.5cm,angle=270]{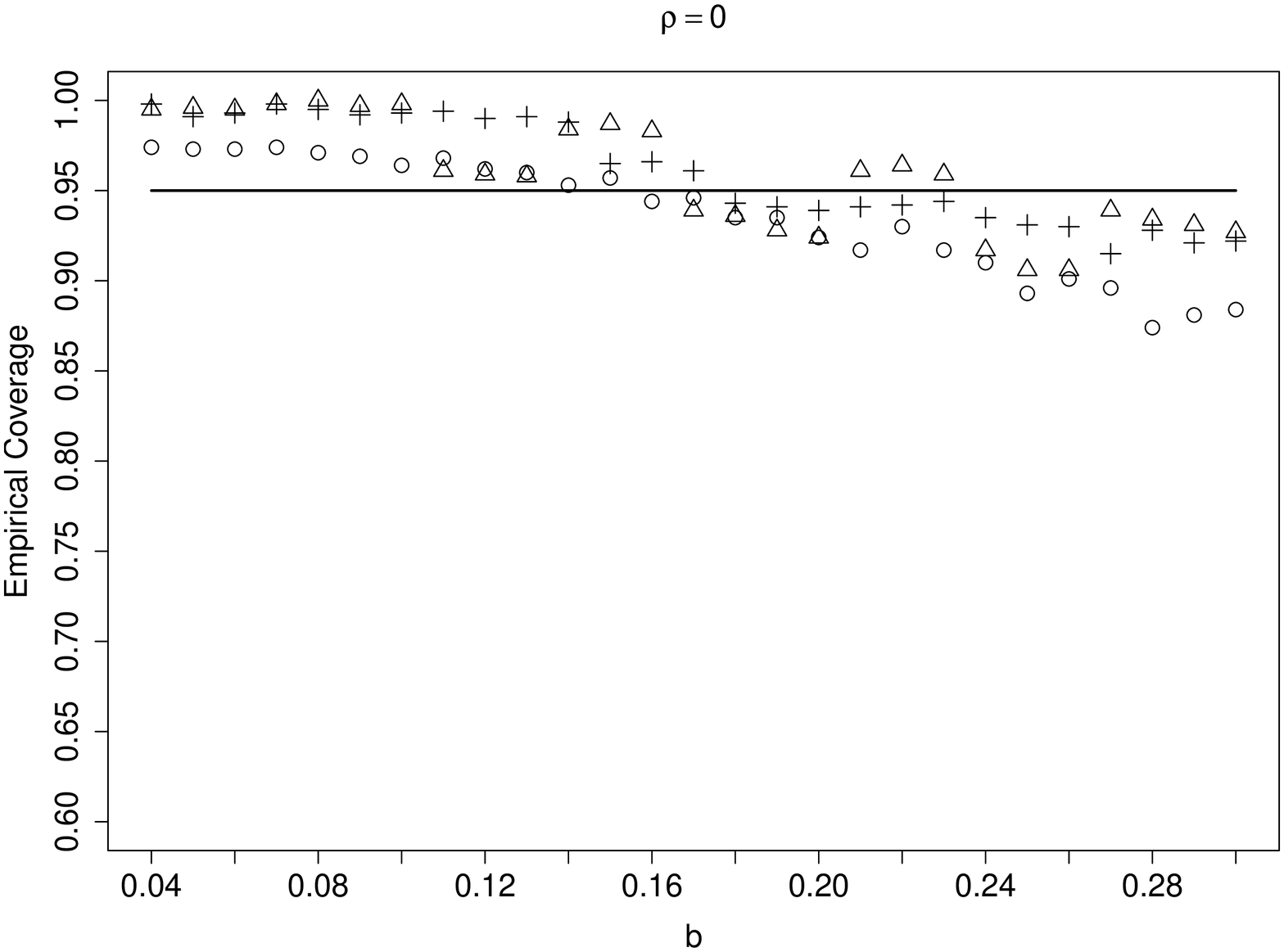}}
{\includegraphics[height=8cm,width=4.5cm,angle=270]{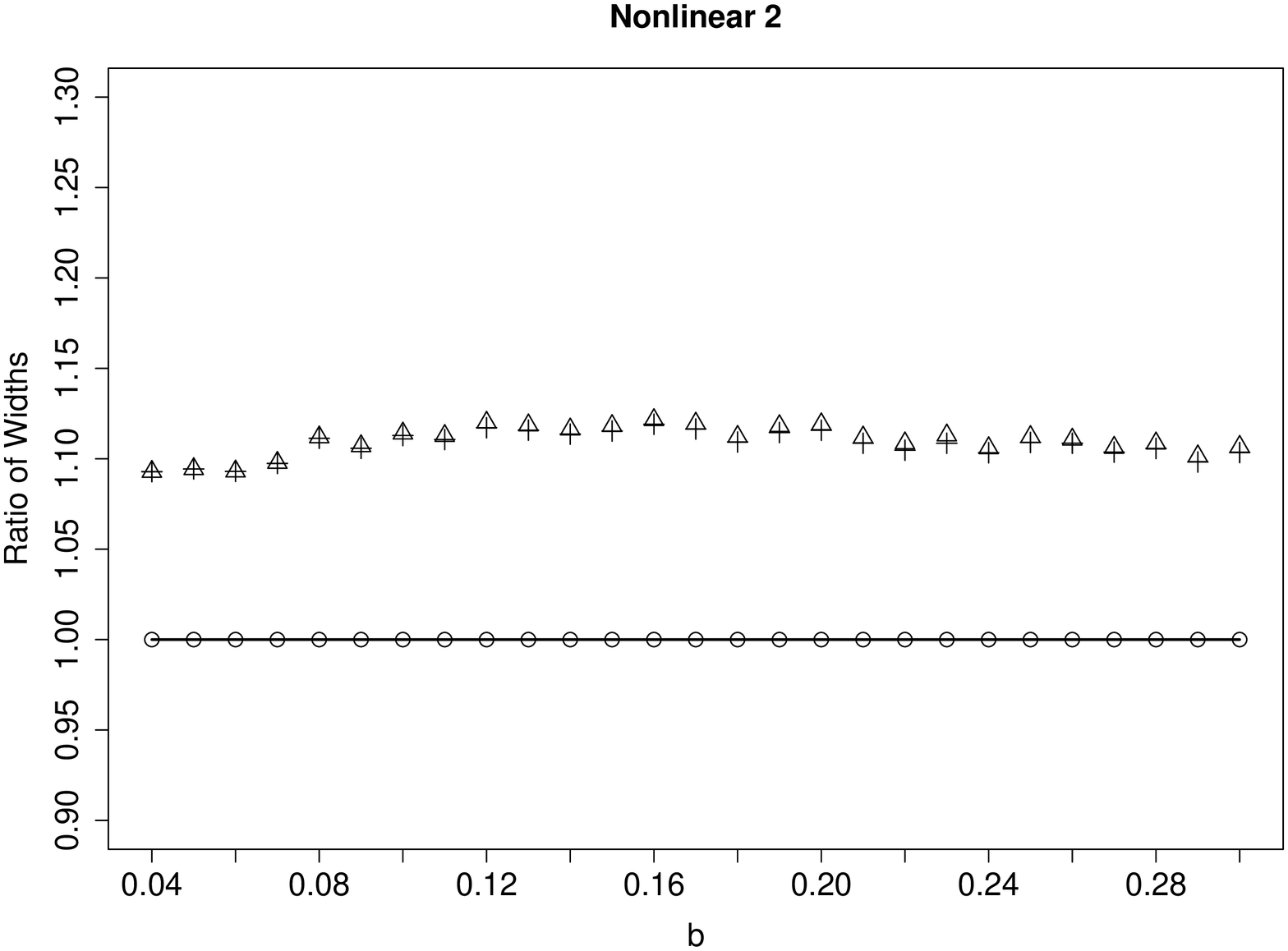}}
{\includegraphics[height=8cm,width=4.5cm,angle=270]{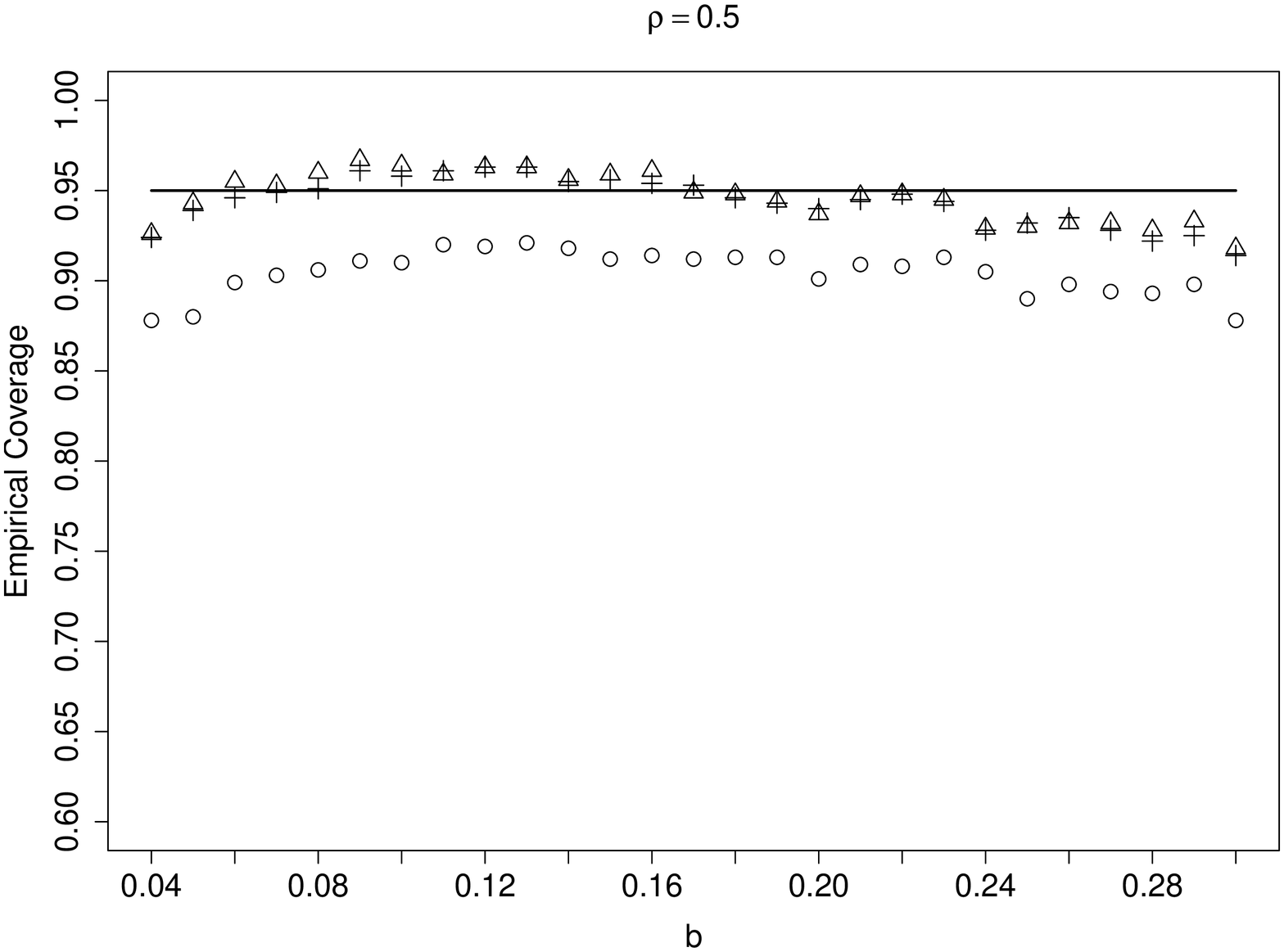}}
{\includegraphics[height=8cm,width=4.5cm,angle=270]{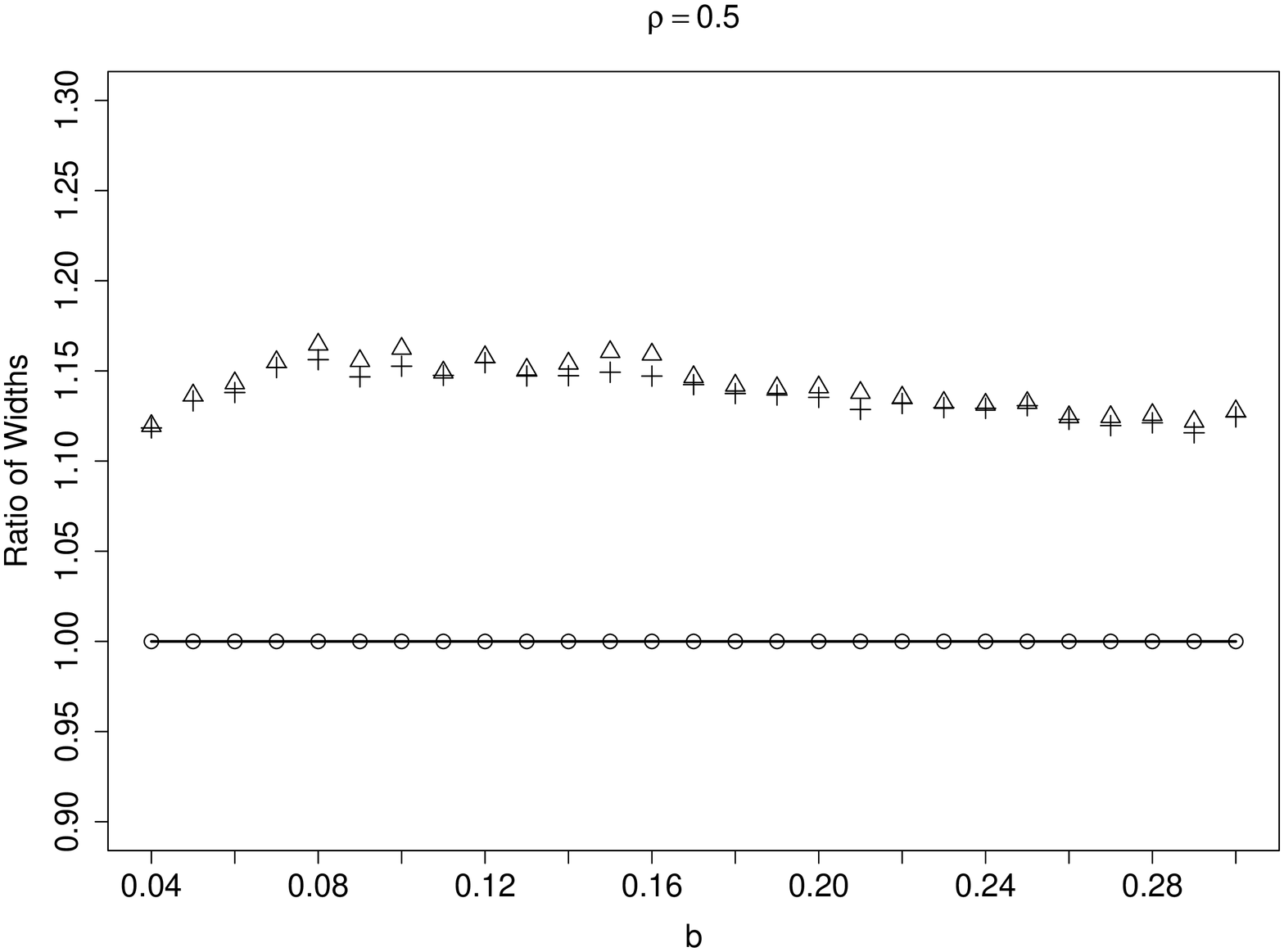}}
{\includegraphics[height=8cm,width=4.5cm,angle=270]{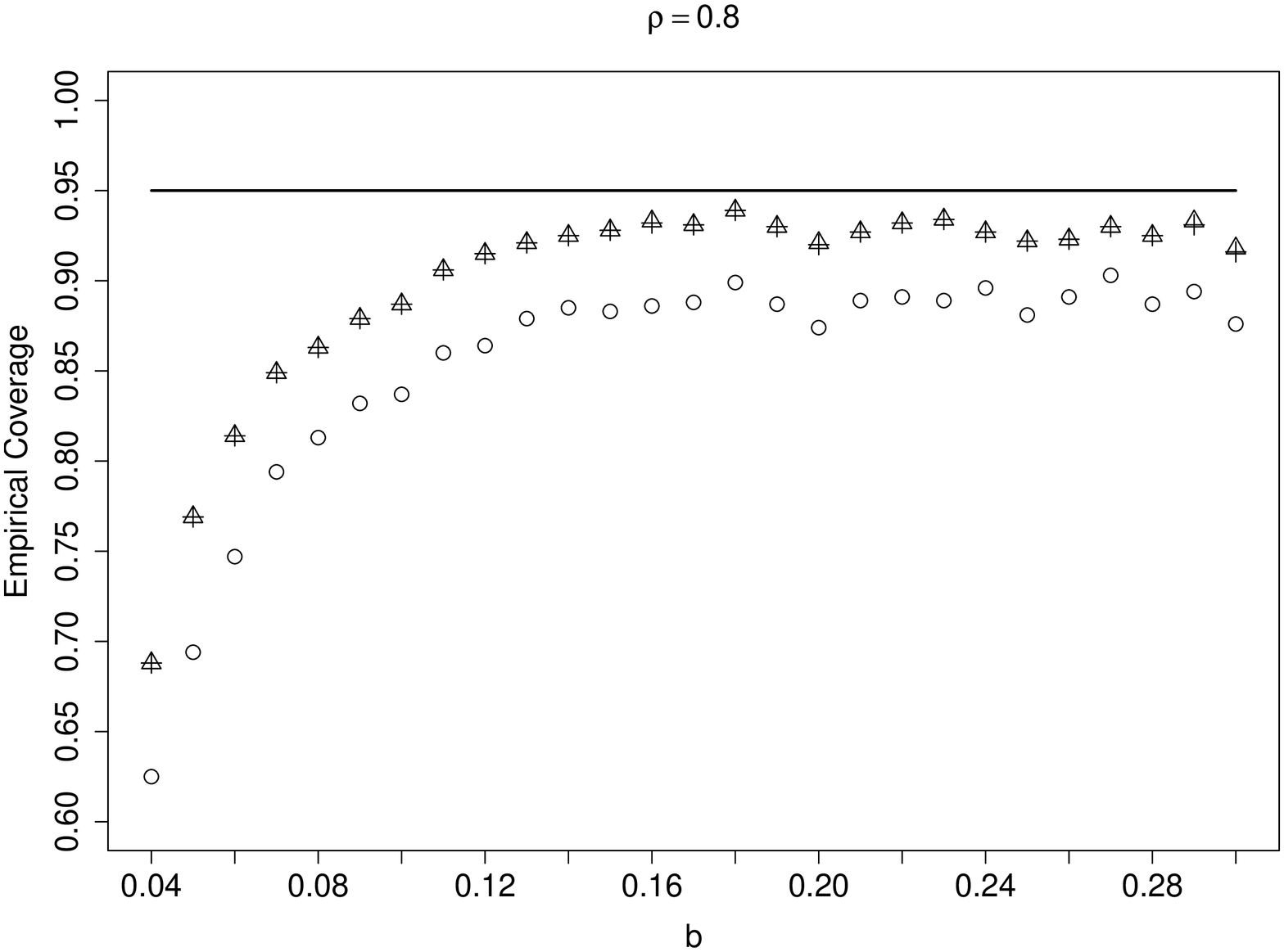}}
{\includegraphics[height=8cm,width=4.5cm,angle=270]{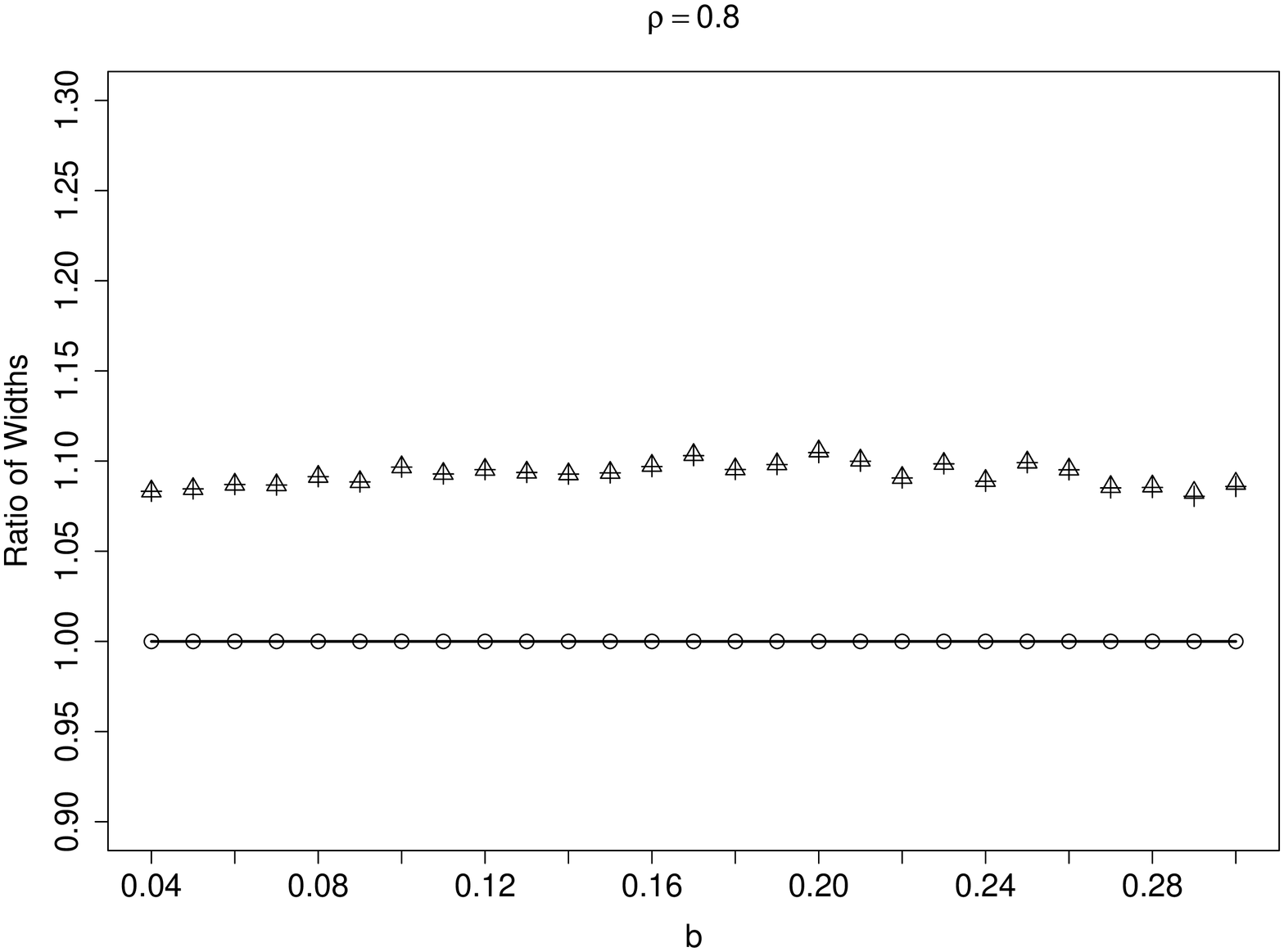}}
\label{fig:band2}
\end{center}

\end{figure}

 Figures~\ref{fig:band1}-\ref{fig:band2} have the same format as Figure~\ref{fig:region} and their right panels show the ratios of the mean band  widths
  over that delivered by the uncalibrated traditional subsampling-based band.   For the marginal distribution function,
 there is an apparent undercoverage for the traditional subsampling-based confidence band in all cases with
   coverage errors increasing  with respect to the magnitude of  dependence (compare the plots for $\rho=0, 0.5, 0.8$), especially at small $b$s.
When $\rho=0,0.5,0.8$ and for almost all $b=0.01,\cdots,0.2$, the coverages delivered by the calibrated bands based on fixed or data driven subsampling width are closer to the nominal level than the traditional counterpart.
 When $\rho=-0.6$,  the calibrated bands based on the fixed or data-dependent bandwidths improve the coverage
  when $b\ge 0.04$, but fails to do so when $b=0.01, 0.02, 0.03$, suggesting that potential improvements can be made about
  the selection of $n'$. In all cases,  the calibrated bands are slightly wider than the uncalibrated counterpart, but the ratios are quite close to $1$. The ``better coverage but wider band" phenomenon is in accordance with the ``better coverage but wider interval" finding in the scalar parameter case. The two calibrated bands  perform similarly in most situations, and their performance is strikingly close when $\rho=0.8$. As seen from Figure~\ref{fig:band2}, which plots the empirical coverages and ratios of mean band widths with respect to $b=0.04,\cdots,0.3$ for the normalized spectral distribution, the improvement of the calibration in terms of coverage error is quite substantial when $\rho=-0.6,0.5,0.8$. In the case $\rho=0$, the calibrated bands
are conservative when $b$ is relatively small, but again provides some improvement when $b$ is close to $0.3$.   Overall the results for the normalized spectral distribution function are qualitatively similar to those for the  marginal distribution function.  Based on the simulation results for confidence intervals reported in  Section~\ref{subsec:interval} and for confidence regions and bands reported in this subsection, it appears that the calibration works very effectively when the series is positively  dependent.

\section{Conclusion}
\label{sec:conclusion}

Subsampling and block-based bootstrap methods have been shown to be widely applicable to many
inference problems in time series analysis. The fixed-$b$ asymptotics developed here explicitly captures the choice of bandwidth parameter in subsampling and the moving block bootstrap, and the resulting  first order approximation is expected to be more accurate than that provided by the small-$b$ asymptotics. As demonstrated in Section~\ref{sec:simulation}, the fixed-$b$ based calibrated confidence intervals (regions, bands) provide  an unambiguous
improvement over the uncalibrated counterparts  in terms of coverage errors in most cases considered. Our calibration method is developed by
estimating the sampling distribution of the p-value, which relates to the prepivoting method by Beran (1987, 1988) and the confidence coefficient calibration method by Loh (1987).  However,
our proposal differs from theirs in two important respects: (i) the limiting null distribution of the p-value is not (necessarily) $U(0,1)$, which is the case for Beran (1987, 1988). In our setting, a pivotal  limiting distribution exists in the scalar parameter case, but not in the case of vector parameter and infinite dimensional parameter, for which the subsampling method is used to provide a good approximation; (ii) Their discussions are limited to the iid setting and inference for finite dimensional parameters. In contrast, our treatment goes substantially beyond their developments by allowing for time series data and the inference of infinite dimensional parameters. Coupled with the recently developed fixed-$b$ approach [Kiefer and Vogelsang (2005)] in econometrics literature, we provide a general recipe for the calibration of the traditional resampling-based inference procedures when smoothing parameters, such as window width in subsampling and block size in the moving block bootstrap  are used to accommodate the dependence.

To conclude the paper, we provide a discussion of   open problems and possible extensions.
(1) Our method can be used as a calibration tool for a properly chosen  smoothing parameter and  it is practically important to choose the smoothing parameter in a sensible way.
The choice of subsampling width and block size for the block-based bootstrap has been discussed in Chapter 9 of Politis et al. (1999a) and Chapter 7 of Lahiri (2003).  It seems natural to ask if it is meaningful to consider the optimal smoothing parameter selection from   a fixed-$b$ based viewpoint, as opposed to the small-$b$ based  approach [see e.g. B\"uhlmann and K\"unsch (1999) and Politis and White (2004)]. A high order expansion of the sampling distribution of the p-value under the null and alternative seems needed to tackle this issue.
(2) The development in this article is confined to time series, although
subsampling and block based bootstrap methods have been extended to spatial settings [see Chapter 5 of Politis et al. (1999a) and
Chapter 12 of Lahiri (2003) and references therein]. An extension of the fixed-$b$ based calibration idea to spatial settings is expected to be
possible but seems nontrivial for irregularly spaced spatial data.
(3) In addition, we impose the weak dependence throughout so the asymptotic normality or functional central limit theorem with
$\sqrt{n}$ convergence rate hold. When the time series is long-range dependent, the subsampling method has been proved to be consistent
 in some situations [see Hall et al. (1998), Nordman and Lahiri (2005)]. It would be interesting to extend the fixed-$b$ approach to calibrate the subsampling based inference in these settings.
(4) A close relative of the block-based bootstraps is the so-called sieve bootstrap [B\"uhlmann (1997)], which also involves
 a bandwidth parameter (i.e., the order of the approximating autoregressive model). It is natural to ask whether it is possible to extend the fixed-$b$ approach
 to calibrate the sieve bootstrap based confidence sets.  We leave these possible extensions for future work.

\section{Appendix}
\label{sec:appendix}

\noindent Proof of Theorem~\ref{th:finite}: For the convenience of notation, let $Y_h=IF(X_h;F)$ and
 $\Delta=\Sigma(P)^{1/2}$. Further let $T_{n,j}={\bf 1}(\|\sqrt{n}(\hat{\theta}_n-\theta_0)\|\le \|\sqrt{l}(\hat{\theta}_{j,j+l-1}-\hat{\theta}_n)\|)$
and $\widetilde{T}_{n,j}={\bf 1}[\|n^{-1/2}\sum_{j=1}^{n} Y_j \|\le \| l^{-1/2} \{\sum_{h=j}^{j+l-1} Y_h - (l/n)\sum_{j=1}^{n} Y_j\}\|]$.
Then $\widetilde{pval}_{n,l}^{SUB}=N^{-1}\sum_{j=1}^{N} T_{n,j}$. Let $D_n(\epsilon)=\{ \|\sqrt{n}R_{1,n}\|<\epsilon, \sup_{j=1,\cdots,N} \|\sqrt{l}R_{j,j+l-1}\| <\epsilon\}$ for any $\epsilon>0$. Then $P\{D_n(\epsilon)\}\rightarrow 1$ as $n\rightarrow\infty$. On $D_n(\epsilon)$, we have that
$$|T_{n,j}-\widetilde{T}_{n,j}|\le {\bf 1}\left[\left|\left\|n^{-1/2}\sum_{j=1}^{n} Y_j \right\| - \left\| l^{-1/2} \left\{\sum_{h=j}^{j+l-1} Y_h - (l/n)\sum_{j=1}^{n} Y_j\right\}\right\|\right|\le 2\epsilon\right].$$
So  the expression $N^{-1}\sum_{j=1}^{N} |T_{n,j}-\widetilde{T}_{n,j}|$ is bounded by
 $$N^{-1}\sum_{j=1}^{N}{\bf 1}\left[\left|\left\|n^{-1/2}\sum_{j=1}^{n} Y_j \right\| - \left\| l^{-1/2} \left\{\sum_{h=j}^{j+l-1} Y_h - (l/n)\sum_{j=1}^{n} Y_j\right\}\right\|\right|\le 2\epsilon\right],$$
  which, by the continuous mapping theorem, converges in distribution to $I(b,\epsilon)$,
where
$$I(b,\epsilon):=(1-b)^{-1}\int_0^{1-b} {\bf 1}\left[ \left|\|\Delta W_k(1)\| - \left\|\Delta \{W_k(b+t)-W_k(t)-bW_k(1)\}/\sqrt{b}\right\| \right|\le 2\epsilon \right]dt.$$
It is not hard to see that for each $t\in [0,1-b]$, the integrand in $I(b,\epsilon)\downarrow 0$ almost surely as $\epsilon\downarrow 0$, which implies that $\lim_{\epsilon\downarrow 0} I(b,\epsilon) = 0$ almost surely by the Lebesgue dominated convergence theorem.
Since $N^{-1}\sum_{j=1}^{N}\widetilde{T}_{n,j}\rightarrow_{D} \widetilde{G}(b;k) $ by the continuous mapping theorem, the conclusion follows by letting
$\epsilon\downarrow 0$ and $n\rightarrow\infty$.

We provide a justification for the use of the continuous mapping theorem above.
For any $x\in D^k[0,1]$, define the functional
$$f_1(x)=(1-b)^{-1}\int_0^{1-b} {\bf 1}\left[ \left|\|x(1)\| - \left\|\{x(b+t)-x(t)-b x(1)\}/\sqrt{b}\right\| \right|\le 2\epsilon \right]dt.$$
and
$$f_2(x)=(1-b)^{-1}\int_0^{1-b} {\bf 1}\left[ \|x(1)\| \le \left\|\{x(b+t)-x(t)-b x(1)\}/\sqrt{b}\right\|  \right]dt.$$

To use the continuous mapping theorem, we need to show that both $f_1$ and $f_2$ are $\Delta W_k(\cdot)-$continuous almost surely. We shall focus on $f_1$ and
the same argument applies to $f_2$.  Define $D_{f_1}=\{x: f_1~\mbox{ is not continuous at}~ x\}$. Then
$$D_{f_1}\subset \widetilde{D}_{f_1}=\{x: \lambda\{t\in [0,1-b]: \|x(1)\| - \left\|\{x(b+t)-x(t)-b x(1)\}/\sqrt{b}\right\| = \pm 2\epsilon\}>0\},$$
where $\lambda$ stands for Lebesgue measure. It is enough to show $P(\Delta W_k(\cdot)\in \widetilde{D}_{f_1})=0$. To
this end, we note that
 \begin{eqnarray*}
&&\E \int _0^{1-b} {\bf 1}(\|\Delta W_k(1)\| - \left\|\{\Delta W_k(b+t)-\Delta W_k(t)-b \Delta W_k(1)\}/\sqrt{b}\right\| = \pm 2\epsilon)dt
\\
&&\hspace{1cm}=\int_0^{1-b} P(\|\Delta W_k(1)\| - \left\|\{\Delta W_k(b+t)-\Delta W_k(t)-b \Delta W_k(1)\}/\sqrt{b}\right\| = \pm 2\epsilon)dt=0
 \end{eqnarray*}
 where we have used the fact that for each $t\in [0,1-b]$,
 \begin{eqnarray}
 \label{eq:fact1}
 P(\|\Delta W_k(1)\| - \left\|\{\Delta W_k(b+t)-\Delta W_k(t)-b \Delta W_k(1)\}/\sqrt{b}\right\| = \pm 2\epsilon)=0
 \end{eqnarray}
 The fact (\ref{eq:fact1}) can be easily shown by noticing that the joint distribution of $(\Delta W_k(1),\{\Delta W_k(b+t)-\Delta W_k(t)-b \Delta W_k(1)\}/\sqrt{b})$ is multivariate normal with
 a positive definite covariance matrix.
 So $P(\Delta W_k(\cdot)\in\widetilde{D}_{f_1})=0$ holds and the use of the continuous mapping theorem is justified. The proof is thus complete.
 \qed

\noindent Proof of Theorem~\ref{th:theoremregion}: The proof is similar to that of Theorem~\ref{th:theorem2}, so we omit the details.
\qed

\noindent Proof of Theorem~\ref{th:theorem2}:
(a) The proof follows from the use of the continuous mapping theorem. Here
the mapping $f: D([-\infty,\infty]\times [0,1])\rightarrow \R$ is defined as
\[f(x)=(1-b)^{-1}\int_0^{1-b} {\bf 1}\left\{ \sup_{s\in\R} |x(s,r+b) - x(s,r) - bx(s,1)|/\sqrt{b}  \ge \sup_{s\in \R} |x(s,1)|\right\} dr.\]
Following the argument in the proof of Theorem~\ref{th:finite}, we can show that if $\widetilde{V}_{b}(r,0)=0$ for every $r\in [0,1-b]$, then
the mapping $f$ is $K$-continuous, i.e., the probability that the Gaussian process $K(\cdot,\cdot)$ falls into the discontinuity set of $f$ is zero.
This completes the proof.

(b) In view of the continuity assumption of $J(b)$
and the monotonicity of $J_{n,n'}(x)$, it suffices to show $J_{n,n'}(x)=P(J(b)\le x)+o_p(1)$ for each $x\in \R$.
 Let
\[\widehat{h}_{n',t}=(N')^{-1}\sum_{j=t}^{t+N'-1} {\bf 1}\left\{ \sqrt{l'} \|m_{j,j+l'-1}-m_{t,t+n'-1}\|_{\infty} \ge \sqrt{n'}\|m_{t,t+n'-1}-m\|_{\infty} \right\}\]
 for $t=1,\cdots,n-n'+1$ and  $\widehat{J}_{n,n'}(x)=(n-n'+1)^{-1}\sum_{t=1}^{n-n'+1} {\bf 1}(\widehat{h}_{n',t}\le x)$.  Note that
  \begin{eqnarray*}
  &&\sqrt{n'}\|m_{t,t+n'-1}-m\|_{\infty} - \sqrt{n'} \|m_{n}-m\|_{\infty} \le \sqrt{n'}\|m_{t,t+n'-1}-m_n\|_{\infty}\\
  &&\hspace{2cm} \le  \sqrt{n'}\|m_{t,t+n'-1}-m\|_{\infty} + \sqrt{n'}\|m_{n}-m\|_{\infty}.
  \end{eqnarray*}
For any $\epsilon>0$,  let $E_n(\epsilon)=\{\sqrt{n'} \|m_{n}-m\|_{\infty} \le \epsilon\}$ and
$$V_{b}(r,\epsilon):=P\left\{\left|\sup_{s\in\R} |K(s,r+b) - K(s,r) - bK(s,1)|/\sqrt{b}  - \sup_{s\in \R} |K(s,1)|\right|\le \epsilon\right\}.$$
 Then $P\{E_n(\epsilon)\}\rightarrow 1$ as $n\rightarrow\infty$.
  On $E_n(\epsilon)$, we have that for each $t=1,\cdots, n-n'+1$, $|h_{n',t}-\widehat{h}_{n',t}| \le W_n(t;\epsilon)$, where
\[W_n(t;\epsilon) :=   (N')^{-1}\sum_{j=t}^{t+N'-1} {\bf 1}\left\{ \left|\sqrt{l'} \|m_{j,j+l'-1}-m_{t,t+n'-1}\|_{\infty} - \sqrt{n'}\|m_{t,t+n'-1}-m\|_{\infty} \right| \le \epsilon \right\}.\]
By stationarity, we have that
  \begin{eqnarray*}
  \E|h_{n',t}-\widehat{h}_{n',t}| {\bf 1}\{E_n(\epsilon)\} &\le& (N')^{-1}\sum_{j=1}^{N'}  P\left\{\left| \sqrt{l'} \|m_{j,j+l'-1}-m_{1,n'}\|_{\infty}-\sqrt{n'}\|m_{1,n'}-m\|_{\infty} \right| \le \epsilon
  \right\}\\
  &\rightarrow& L(b,\epsilon):=(1-b)^{-1}\int_0^{1-b} V_b(r,\epsilon) dr
  \end{eqnarray*}
 The above convergence  follows from Theorem 3 of Ferguson (1996)
  and  the fact that $W_n(1;\epsilon)\rightarrow_{D} J(b,\epsilon)$, where
 \[ J(b,\epsilon):=(1-b)^{-1}\int_0^{1-b} {\bf 1}\left\{ \left|\sup_{s\in\R} |K(s,r+b) - K(s,r) - bK(s,1)|/\sqrt{b}  - \sup_{s\in \R} |K(s,1)|\right|\le \epsilon\right\} dr \]
Again the continuous mapping theorem is invoked to derive the weak convergence of $W_n(1,\epsilon)$
 and following the argument in the proof of Theorem~\ref{th:finite}, its use can be justified under the assumption that  $\widetilde{V}_b(r,\epsilon)=0$ for
 each $r\in [0,1-b]$ and $\epsilon\ge 0$.

Next it is not hard  to see that
   $\lim_{\epsilon\downarrow 0} L(b,\epsilon)=0$ since $V_b(r,\epsilon)\downarrow V_b(r,0)=0$ as $\epsilon\downarrow 0$ for every $r\in [0,1-b]$. Thus
   $\sup_{t=1,\cdots,n-n'+1}\E|h_{n',t}-\widehat{h}_{n',t}|\le \E|h_{n',1}-\widehat{h}_{n',1}|{\bf 1}\{E_n(\epsilon)\}+2P(E_n(\epsilon)^c)\le \epsilon$
    for large enough $n$.   Furthermore,
 \begin{eqnarray*}
 &&\widehat{J}_{n,n'}(x-\sqrt{\epsilon})-(n-n'+1)^{-1}\sum_{t=1}^{n-n'+1} {\bf 1} \{ |\widehat{h}_{n',t}-h_{n',t}|\ge \sqrt{\epsilon}\} \le J_{n,n'}(x)\\
 &&\hspace{2cm}\le \widehat{J}_{n,n'}(x+\sqrt{\epsilon})+(n-n'+1)^{-1}\sum_{t=1}^{n-n'+1} {\bf 1} \{ |\widehat{h}_{n',t}-h_{n',t}|\ge \sqrt{\epsilon}\}.
 \end{eqnarray*}
 By the Markov inequality, $(n-n'+1)^{-1}\sum_{t=1}^{n-n'+1} P\{ |\widehat{h}_{n',t}-h_{n',t}|\ge \sqrt{\epsilon}\}\le \sqrt{\epsilon}$.
 Using the same argument in the proof of Theorem 3.2.1. of Politis et al. (1999a),
 we can show that $\widehat{J}_{n,n'}(x)-P\{J(b)\le x\}=o_p(1)$, which follows from the stationarity and strong mixing properties of $X_t$ and
 the boundness of $\{\widehat{h}_{n',t}\}_{t=1}^{n-n'+1}$. The conclusion then follows from an elementary argument.

\qed

\bigskip

\centerline{\bf\large References}



\par\noindent\hangindent2.3em\hangafter 1
Anderson, T., (1993) Goodness of fit tests for spectral distributions. {\it Annals of Statistics}, 21, 830-847.

\par\noindent\hangindent2.3em\hangafter 1
Aza\"is, J. M. and Wschebor, M. (2001) On the regularity of the distribution of the maximum of one-parameter Gaussian processes.
{\it Probability Theory and Related Fields}, 119, 70-98.

\par\noindent\hangindent2.3em\hangafter 1
Beran, R. (1987) Previoting to reduce level error of confidence sets. {\it Biometrika}, 74, 457-468.

\par\noindent\hangindent2.3em\hangafter 1
Beran, R. (1988) Prepivoting test statistics: A bootstrap view of asymptotic refinements.  {\it Journal of the American Statistical Association}, 83, 687-697.

\par\noindent\hangindent2.3em\hangafter 1
Berg, A.,  McMurry, T. and  Politis, D. N. (2010) Subsampling p-values. {\it Statistics and Probability Letters}, 1358-1364.

\par\noindent\hangindent2.3em\hangafter 1
Berkes, I., H\"ormann, S. and Schauer, J. (2009) Asymptotic results for the empirical process of stationary sequences.
{\it Stochastic Processes and their Applications}, 119, 1298-1324.

\par\noindent\hangindent2.3em\hangafter 1
Bickel, P. and Sakov, A. (2008) On the choice of $m$ in the $m$ out of $n$ bootstrap and confidence bounds for extrema.
{\it Statistica Sinica}, 18,  967-985.

\par\noindent\hangindent2.3em\hangafter 1
Billingsley, P. (1968) {\it Convergence of Probability Measures}, Wiley.

\par\noindent\hangindent2.3em\hangafter 1
Brillinger, D. (1975) {\it Time Series: Data Analysis and Theory}, Holt, Rinehart and Winston, New York.

\par\noindent\hangindent2.3em\hangafter 1
B\"uhlmann, P. (1994) Blockwise bootstrapped empirical process
for stationary sequences. {\it Annals of Statistics}, 22, 995-1012.

\par\noindent\hangindent2.3em\hangafter 1
B\"uhlmann, P. (1997) Sieve bootstrap for time series. {\it Bernoulli}, 3, 123-148.

\par\noindent\hangindent2.3em\hangafter 1
B\"uhlmann, P. and K\"unsch, H. R. (1999). Block length selection in the bootstrap for time series.
{\it Computational Statistics and Data Analysis}, 31, 295-310.

\par\noindent\hangindent2.3em\hangafter 1
Dahlhaus, R. (1985a) On the asymptotic distribution of Bartlett's $U_p$-statistic. {\it Journal of Time Series Analysis}, 6(4), 213-227.

\par\noindent\hangindent2.3em\hangafter 1
Dahlhaus, R. (1985b) A functional central limit theorem for tapered empirical distribution functions.
{\it Stochastic Processes and its Applications}, 19, 135-149.

\par\noindent\hangindent2.3em\hangafter 1
Diebolt, J. and Posse, C. (1996) On the density of the maximum of smooth Gaussian processes.
{\it Annals of Probability}, 24, 1104-1129.

\par\noindent\hangindent2.3em\hangafter 1
Ferguson, T. S. (1996) {\it A Course in Large Sample Theory}. Chapman\&Hall.

\par\noindent\hangindent2.3em\hangafter 1
Goncalves, S. and Vogelsang, T. J. (2011) Block bootstrap HAC robust tests: sophistication of the naive bootstrap.
{\it Econometric Theory}, 27, 745-791.

\par\noindent\hangindent2.3em\hangafter 1
G\"otze, F. and Ra$\breve{c}$kauskas, A. (2001) Adaptive choice of bootstrap sample sizes. In {\it State of the
Art in Probability and Statistics}, IMS Lecture Notes Monogr. Ser., 36 (ed. Aad van der Vaart
Mathisca de Gunst Chris Klaassen), pp. 286-309. Cambridge University Press.

\par\noindent\hangindent2.3em\hangafter 1
Hall, P. (1992) {\it The Bootstrap and Edgeworth Expansion}, New York: Springer.

\par\noindent\hangindent2.3em\hangafter 1
Hall, P., Lahiri, S. N., and Jing, B.-Y. (1998) On the subsampling window method for long-range dependent
data. {\it Statistica Sinica}, 8, 1189-1204.

\par\noindent\hangindent2.3em\hangafter 1
Hampel, F., Ronchetti, E., Rousseeuw, P. and Stahel, W. (1986) {\it Robust Statistics: The Approach Based on Influence Functions}, New
York: John Wiley.

\par\noindent\hangindent2.3em\hangafter 1
 Hashimzade, N. and Vogelsang, T. J. (2008) Fixed-b asymptotic approximation of the sampling behavior of nonparametric spectral density estimators. {\it Journal of Time Series Analysis}, 29, 142-162.

\par\noindent\hangindent2.3em\hangafter 1
Jach, A., McElroy, A. and Politis, D. N.  (2011) Subsampling inference for the mean
of heavy-tailed long memory time series. to appear in {\it Journal of  Time Series Analysis}.

\par\noindent\hangindent2.3em\hangafter 1
{Jansson, M.} (2004) The error rejection probability of simple autocorrelation
 robust test. {\it Econometrica}, 72, 937-946.

\par\noindent\hangindent2.3em\hangafter 1
{Kiefer, N. M.}, and {Vogelsang, T. J.} (2005) A new asymptotic theory for heteroskedasticity-autocorrelation robust tests.
 {\it Econometric Theory}, 21, 1130-1164.

\par\noindent\hangindent2.3em\hangafter 1
Kiefer, N. M., Vogelsang, T. J. and Bunzel, H. (2000) Simple robust testing of regression hypotheses.
{\it Econometrica}, 68, 695-714.

\par\noindent\hangindent2.3em\hangafter 1
{K\"unsch, H.} (1989). The jackknife and the bootstrap for general stationary observations. {\it Annals of Statistics}
 {\bf 17}, 1217-1241.

\par\noindent\hangindent2.3em\hangafter 1
Lahiri, S. N. (2001) Effects of block lengths on the validity of block resampling methods.
{\it Probability Theory and Related Fields}, 121, 73-97.

\par\noindent\hangindent2.3em\hangafter 1
Lahiri, S. N. (2003) {\it Resampling Methods for Dependent Data}, New York:
Springer.

\par\noindent\hangindent2.3em\hangafter 1
Lee, S. M. S. and Lai, P. Y. (2009) Double block bootstrap confidence intervals for dependent data. {\it Biometrika},
96(2), 427-443.

\par\noindent\hangindent2.3em\hangafter 1
{Liu, R. Y.} and {Singh, K.} (1992). Moving blocks jackknife and bootstrap capture weak dependence. In {\it Exploring
the Limits of Bootstrap}, (Ed. R. LePage and L.Billard), 225-248.  John Wiley, New York.


\par\noindent\hangindent2.3em\hangafter 1
Loh, W.-Y. (1987) Calibrating confidence coefficients. {\it Journal of the American Statistical Association}, 82, 155-162.

\par\noindent\hangindent2.3em\hangafter 1
Loh, W.-Y. (1991) Bootstrap calibration for confidence interval construction and selection. {\it Statistica Sinica},
1, 477-491.

\par\noindent\hangindent2.3em\hangafter 1
Naik-Nimbalkar, U. and Rajarshi, M. (1994) Validity of blockwise bootstrap for empirical processes
with stationary observations. {\it Annals of Statistics}, 22, 980-994.

\par\noindent\hangindent2.3em\hangafter 1
Nordman, D. and Lahiri, S. N. (2005) Validity of sampling window method for linear long-range dependent
processes. {\it Econometric Theory}, 21, 1087-1111.

\par\noindent\hangindent2.3em\hangafter 1
{Paparoditis, E.} and {Politis, D. N.} (2001). Tapered block bootstrap. {\it Biometrika} {\bf 88}, 1105-1119.

\par\noindent\hangindent2.3em\hangafter 1
{Paparoditis, E.} and {Politis, D. N.} (2002). The tapered block bootstrap for general statistics from stationary sequences. {\it Econometrics Journal} {\bf 5}, 131-148.

\par\noindent\hangindent2.3em\hangafter 1
Politis, D. N. (2011)  Higher-order accurate, positive semi-definite estimation of large-sample covariance and spectral density matrices. {\it  Econometric Theory} {\bf 27}, 703-744.

\par\noindent\hangindent2.3em\hangafter 1
Politis, D. N. and Romano, J. P. (1992) A circular block-resampling procedure for stationary data.  {\it Exploring the Limits of Bootstrap}, (Raoul LePage and Lynne Billard, eds.), John Wiley,  p263-270.

\par\noindent\hangindent2.3em\hangafter 1
Politis, D. N. and Romano, J. P. (1994) Large sample confidence regions based on subsamples under minimal assumptions. {\it  Annals of Statistics},  22,  2031-2050.


\par\noindent\hangindent2.3em\hangafter 1
{Politis, D. N.}, {Romano, J. P.} and {Wolf, M.} (1999a), {\it
Subsampling}, Springer-Verlag, New York.

\par\noindent\hangindent2.3em\hangafter 1
{Politis, D. N.}, {Romano, J. P.} and {You, L.} (1993) Uniform confidence bands for the spectrum based on subsamples, in Computing Science and Statistics, Proceedings of the 25th Symposium on the Interface, San Diego, California. (M. Tarter and M. Lock, eds.), The Interface Foundation of North America,  346-351.

\par\noindent\hangindent2.3em\hangafter 1
{Politis, D. N.}, {Romano, J. P.} and {Wolf, M.} (1999b) Weak convergence
of dependent empirical measures with application to subsampling in function
spaces. {\it Journal of Statistical Planning and Inference}, 79, 179-190.

\par\noindent\hangindent2.3em\hangafter 1
{Politis, D. N.} and {White, H.} (2004) Automatic block-length selection for the dependent bootstrap. {\it Econometric Reviews},
 23,  53-70.

\par\noindent\hangindent2.3em\hangafter 1
Sayginsoy, O. and Vogelsang, T. J. (2011) Testing for a shift in trend at an unknown date: a fixed-b analysis of heteroskedasticity autocorrelation robust OLS based tests. {\it Econometric Theory}, 27, 992-1025.

\par\noindent\hangindent2.3em\hangafter 1
Shao, X. (2009) Confidence intervals for spectral mean and ratio statistics. {\it Biometrika}, 96, 107-117.

\par\noindent\hangindent2.3em\hangafter 1
Shao, X. (2010a)  A self-normalized approach to confidence interval construction in time series.  {\it Journal of the Royal Statistical Society, Series B},  72(3), 343-366.  Corrigendum: 2010, 72(5), 695-696.

\par\noindent\hangindent2.3em\hangafter 1
Shao, X. (2010b) Extended tapered block bootstrap. {\it Statistica Sinica}, 20, 807-821.

\par\noindent\hangindent2.3em\hangafter 1
Shao, X. (2010c)  The dependent wild bootstrap. {\it Journal of the American Statistical Association}, 105, 218-235.

\par\noindent\hangindent2.3em\hangafter 1
Shao, X. and Politis, D. (2011) Fixed-$b$  subsampling and block bootstrap: improved confidence sets based on p-value calibration. Technical Report, Department of Statistics,
University of Illinois at Urbana-Champaign. Available at ArXiv

\par\noindent\hangindent2.3em\hangafter 1
Sun, Y., Phillips, P. C. B. and Jin, S. (2008) Optimal bandwidth selection in heteroscedasticity-autocorrelation robust testing.
 {\it Econometrica}, 76, 175-194.

\newpage

\begin{table}[!h]
\caption{Simulated values of $G_{\alpha}(b)$, $\widetilde{G}_{\alpha}(b)$, $H_{\alpha}(b)$ and $\widetilde{H}_{\alpha}(b)$ when fitted with a quadratic polynomial $cv(b)=a_0+a_1b+a_2b^2$, $b\in (0,0.2]$.  $\alpha=0.05,0.1$. The simulated values are based on $n=5000$ and 50000 replications. In the moving block bootstrap case, we use $50000$ bootstrap replications.   }
\label{tb:cv}
\begin{center}
\begin{tabular}{c|cccc}
\hline
&$a_0$&$a_1$&$a_2$&$R^2$\\
\hline
$G_{0.05}(b)$&0.05&-0.2289  &-0.1325&0.9980
\\
$G_{0.1}(b)$&0.1&-0.1039 & -0.8407& 0.9997
\\
$\widetilde{G}_{0.05}(b)$&0.05&-0.3929 &  0.6394&0.9978
\\
$\widetilde{G}_{0.1}(b)$&0.1&-0.3285 & -0.4088 &0.9994
\\
$H_{0.05}(b)$&0.05&-0.3431 &  0.5766&   0.9868
\\
$H_{0.1}(b)$&0.1&-0.4079 &   0.2256&  0.9681
\\
$\widetilde{H}_{0.05}(b)$&0.05&-0.2121  & 0.2624 & 0.9610
\\
$\widetilde{H}_{0.1}(b)$&0.1&-0.2461 &  0.1174  & 0.9584
\\
\hline
\end{tabular}
\end{center}
\end{table}

\newpage

\newpage

\newpage

\newpage


\end{document}